\documentclass[10pt]{article}

\usepackage{amsmath,amssymb}
\usepackage{bm}
\usepackage{graphicx}
\usepackage[english]{babel}
\usepackage{epsfig}
\usepackage{float}
\usepackage{cancel}
\usepackage{ulem}
\usepackage{caption}
\usepackage{subcaption}

\usepackage{color}

\usepackage[margin=1.2in]{geometry}
\usepackage{verbatim}
\usepackage{varwidth}
\usepackage[colorlinks=true]{hyperref}
\hypersetup{urlcolor=blue, citecolor=red}
\newcommand\beq[1]{ \begin{equation}\label{#1} }
\newcommand{\eeq}{ \end{equation} }

\newcommand\beqa[1]{ \begin{eqnarray} \label{#1}}
\newcommand{\eeqa}{ \end{eqnarray} }
\newcommand{\beqano}{ \begin{eqnarray*} }
\newcommand{\eeqano}{ \end{eqnarray*} }

\newcommand{\implica}{\ \Longrightarrow\ }

\newcommand{\R}{\mathbb{R}}     
\newcommand{\N}{\mathbb{N}}     
     
\newcommand{\CC}{{\rm C}}   
   \newcommand{\HH}{{\rm H}}
\newcommand{\KK}{\rm{K}}   
   
\newcommand{\RR}{\mathcal{R}}   \newcommand{\WW}{\mathcal{W}}
   
   \newcommand{\UU}{{\rm U}}
\newcommand{\EE}{{\rm E}}

\newcommand{\II}{{\rm I}}

\newcommand{\bx}{{\mathbf x}}
\newcommand{\by}{{\mathbf y}}
\newcommand{\bC}{{\mathbf C}}
\newcommand{\bG}{{\mathbf G}}
\newcommand{\bi}{{\mathbf i}}
\newcommand{\bj}{{\mathbf j}}
\newcommand{\bk}{{\mathbf k}}
\newcommand{\hh}{{\rm h}}

\newtheorem{theorem}{Theorem}[section]

\newtheorem{definition}[theorem]{Definition}
\newtheorem{proposition}{Proposition}[section]
\newtheorem{conjecture}{Conjecture}[section]
\newtheorem{numevid}{Numerical Evidence}[section]
\newtheorem{remark}{Remark}[section]
\newtheorem{corollary}{Corollary}[section]

\begin{document}

\title{\bf Euler integral as a  source of chaos in the three--body problem} 
\author{Sara Di Ruzza and Gabriella Pinzari} 
\date{\today}
\maketitle
\tableofcontents
\newpage

\begin{abstract}
In this paper we address,  from a purely numerical point of view, the question, raised in \cite{pinzari19, pinzari20b}, and partly considered in \cite{pinzari20a, diruzzaDP20, chenPi2021}, whether a certain function, referred  to as ``Euler Integral'', is a quasi--integral along the trajectories of the three--body problem.  Differently from our previous investigations, here we focus on the region of the ``unperturbed separatrix'', which turns to be complicated by a collision singularity.
Concretely, we reduce the Hamiltonian to two degrees of freedom and, after fixing some energy level, we discuss in detail the resulting three--dimensional phase space around an elliptic and an hyperbolic periodic orbit. After measuring the strength of variation of the Euler Integral (which are in fact small), we detect the existence of chaos closely to the unperturbed separatrix. The latter result is obtained through a careful use of the machinery of covering relations, developed in \cite{GierzkiewiczZ19, ZgliczynskiG2004, WilczakZ2003}.
\end{abstract}


\section{Purpose of the paper}
This paper is a numerical study on the three--body problem. It is to be specified that we deal with (a suitably simplified version of) the Hamiltonian of the {\it full} three--body problem, where ``full'' is used here as opposed to the so--called ``restricted'' problem -- maybe more known to non specialists -- to which much of the arguments discussed here also can also be applied. The full three--body (in general, many--body) problem inherits much of its reputation --especially in Hamiltonian mechanics -- after the breakthrough paper by V. I. Arnold~\cite{arnold63} which will be recalled below. In fact, this paper is motivated by previous research~\cite{pinzari19, diruzzaDP20, pinzari20a, pinzari20b}, which here we briefly recall, in order to keep the paper self--contained.

\noindent
We fix a reference frame $(\bi, \bj, \bk)$ in the Euclidean space, which we identify with $\mathbb R^3$. In such a space we consider three masses $1$, $\mu$ and $\kappa$, with $\mu$, $\kappa<1$, interacting through gravity only. We reduce the translation symmetry relating the positions of two (out of three) masses  to the position of the third one, as described in \cite[\S 5]{giorgilli}. Contrarily to the usual practice, we choose $\mu$ as reference mass (usually, the unit mass is chosen). With such choice, the Hamiltonian governing the motions  of the masses $1$ and $\kappa$ is 
\\
\begin{equation*}
	\HH_{\rm 3b}(\by',\by,\bx',\bx)= \frac{\kappa+\mu}{\kappa\mu} \frac{\|\by\|^2}{2}-\frac{\kappa\mu}{\|\bx\|} + \frac{\mu+1}{\mu} \frac{\|\by'\|^2}{2} - \frac{\mu}{\|\bx'\|}-\frac{\kappa }{\|\bx-\bx'\|} + \frac{1}{\mu}\by\cdot \by' \, .
\end{equation*}
\begin{figure}[H]
 \centering
 \includegraphics[width=10cm,height=8cm,draft=false]{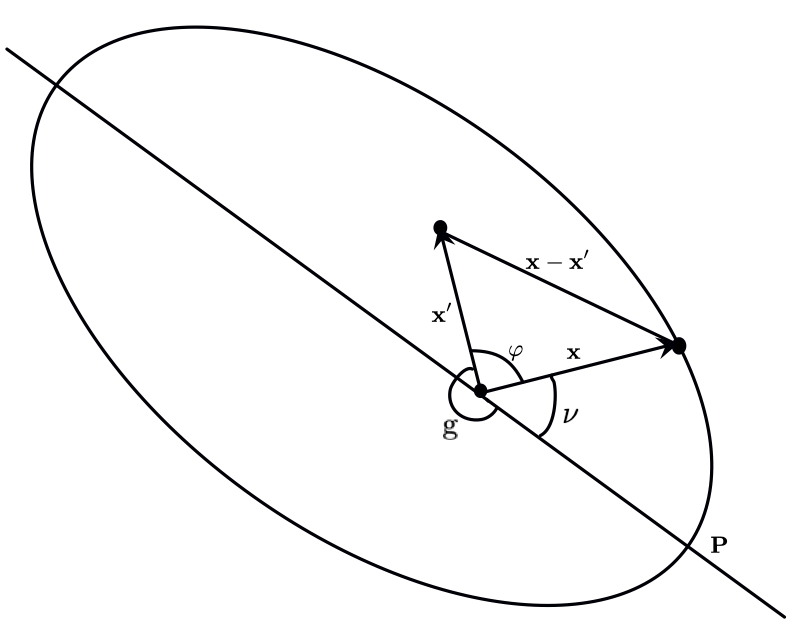}
 \caption{The  three--body problem. }\label{model}
 \end{figure}

 \noindent
 where $\bx'=(x'_1, x'_2, x'_3)$, $\bx=(x_1, x_2, x_3)$ are the position coordinates of $1$ and $k$; $\by'=(y'_1, y'_2, y'_3)$, $\by=(y_1, y_2, y_3)$ are their respective linear momenta; $\|\cdot\|$ denotes the Euclidean distance and, finally, the gravity constant has been conventionally fixed to one.
\\ The rescaling 
	\begin{equation*}
		(\by',\by) \rightarrow \frac{\mu^2\kappa^2}{\kappa+\mu} (\by',\by)\, , \qquad (\bx',\bx) \rightarrow \frac{\kappa+\mu}{\mu^2\kappa^2}(\bx',\bx) \,,\qquad t\to \frac{\mu^3\kappa^3}{\kappa+\mu}t\end{equation*}
(with $t$ denoting the time) does not alter the motion equations, provided that $\HH_{\rm 3b}$ is changed to
\begin{equation}
\label{ham_4dof}
	\HH_{\rm 3b}(\by',\by,\bx',\bx) =  \frac{\|\by\|^2}{2}-\frac{1}{\|\bx\|} + \delta \biggl(  \frac{\|\by'\|^2}{2} -\frac{\alpha }{\|\bx-\bx'\|} -\frac{\beta}{\|\bx'\|} + \gamma \, \by\cdot \by' \biggr) \, ,
\end{equation}
with
\\
\begin{equation}\label{newmasses}
	\alpha:= \frac{\kappa+\mu}{\kappa\mu(\mu+1)} \, , \qquad 
	\beta:= \frac{\kappa+\mu}{\kappa^2(\mu+1)} \, , \qquad
	\gamma:= \frac{1}{\mu+1} \, , \qquad
	\delta:= \frac{\kappa(\mu+1)}{\kappa+\mu} \, .
\end{equation}
As only two parameters among~\eqref{newmasses} can be regarded as independent, from this point on, we choose $\alpha$ and $\beta$. This will simplify later analysis (compare Equation~\eqref{ham} below).  We restrict our attention to the so called ``planar problem'', which corresponds to take the respective third components of position and momentum coordinates identically  vanishing: $x_3=x_3'=y_3=y_3'=0$.
In such a case, $\HH_{\rm 3b}$ in~\eqref{ham_4dof} has four degrees of freedom. We now describe  a procedure which will  reduce the number of degrees of freedom to {\it two}. One degree of freedom can be eliminated   exploiting the ``rotations invariance'', namely the fact that the Hamiltnian $\HH_{\rm 3b}$ remain unchanged   under the group of transformations
\beqa{rotinv}(\by',\by)\rightarrow (\RR \by', \RR\by)\qquad (\bx',\bx)\rightarrow (\RR \bx', \RR\bx) \eeqa
where $\RR$ is any constant orthogonal matrix, i.e, verifying $$\RR \RR^{\rm t}=\II=\RR^{\rm t}\RR$$ with the superscript ``t'' denoting transpose, and $\II$ being the identity matrix. The existence of such  group of diffeomeorphisms 
is  caused by the conservation of the components of the  ``angular momentum'' vector $\bC=(\CC_1, \CC_2, \CC_3)$ given by
\begin{eqnarray*}
	\bC = \bx \times \by + \bx' \times \by'
\end{eqnarray*}
along the trajectories of $\HH_{\rm 3b}$.
Clearly, rotation invariance is not specific of the planar problem. In the planar case, it  allows for the reduction of one\footnote{Incidentally, in the general case, the number of degrees of freedom is lowered by {\it two} units, due to the fact that the components of $\bC$ are not pairwise commuting. See~\cite{chierchiaPi11b} for a case study.} degree of freedom, as
$\bC$ has one only non--trivial coordinate $\CC_3$, which, from now on, we shall simply denote as $\CC$.   One further degree of freedom can be eliminated under the assumption that the ``Keplerian term'' outside parentheses in~\eqref{ham_4dof}, namely,
 \begin{align}\label{Kep}\frac{\|\mathbf y\|^2}{2}-\frac{1}{\|\mathbf x\|}\end{align} 
takes negative values and is  ``leading''  in the  Hamiltonian~\eqref{ham_4dof}. 
 To better specify this assumption, we need to describe canonical coordinates  explicitly performing the reduction of~\eqref{rotinv}  and, simultaneously,
 integrating~\eqref{Kep}. For the planar case, such coordinates are easy to be produced. We proceed as follows.

 \noindent
On a 6--dimensional ``rotation--reduced phase space'' (that will be more precisely described in the next Section~\ref{Set up}) we fix coordinates \begin{eqnarray}\label{coord1}({\rm R}, {\rm G}, \Lambda, {\rm r}, {\rm g}, \ell)\end{eqnarray}
which equip such space with the two--form
$$\omega=d{\rm R}\wedge d{\rm r}+d{\rm G}\wedge d{\rm g}+d{\Lambda}\wedge d{\ell}\, . $$
To define the coordinates~\eqref{coord1}, we note that, as long as the Hamiltonian~\eqref{Kep} keeps to be negative, it generates motions on ellipses. We denote as $\mathbb E$ the ellipse generated by  Hamiltonian~\eqref{Kep} for a given initial datum $(\by, \bx)$. 
 Assuming $\mathbb E$ is not a circle, we let
\begin{itemize}
 \item[--] ${\rm R}
 $ is the ``radial velocity'' of $\bx'$; i.e., the projection of the velocity $\by'$ along the direction of $\bx'$;
 \item[--]${\rm G}$ is the Euclidean length of the angular momentum  $\bG={\mathbf x}\times {\mathbf y}$ 
 of $\bx$;
  \item[--]$\Lambda=\sqrt{a}$, where $a$ is the semi--major axis of ${\mathbb E}$;
 \item[--] ${\rm r}
 $ is the Euclidean length of $\bx'$;
 \item[--] ${\rm g}$ the angle detecting the perihelion of $\mathbb E$;
 \item[--] $\ell$  the ``mean anomaly of $\bx$''.  \end{itemize}
Precise formulae will be given below: compare Equation~\eqref{coord}. Here we only mention that, in a sense, the coordinates above are referred to a frame ``moving with $\bx'$'', in order to obtain reduction of rotations.
Using the coordinates above, and splitting the term inside parentheses in~\eqref{ham_4dof} as the sum of its $\ell$--average (denoted  as $\overline{\rm H}_{\CC}$
) and the zero--average part (denoted  as $\widetilde{\rm H}_{\CC}$), we arrive at

\begin{align}\label{ovlH}{\rm H}_{\rm 3b, \CC}({\rm R}, {\rm G}, \Lambda, {\rm r}, {\rm g}, \ell)=-\frac{1}{2\Lambda^2}+\delta\Big(\overline{\rm H}_{\CC}({\rm R}, {\rm G}, \Lambda, {\rm r}, {\rm g})+\widetilde{\rm H}_{\CC}({\rm R}, {\rm G}, \Lambda, {\rm r}, {\rm g}, \ell)\Big)\,. \end{align} 
In the notation, we have remarked  that, as an effect of the reduction, the system depends parametrically on the total angular momentum $\CC=\|\bC\|$.\\
The Hamiltonian~\eqref{ovlH} has an involved aspect. By no means it appears as -- or can be conjugated to -- a standard {\it close--to--be--integrable system}. These are systems of the form
$$\HH(\II, \varphi)=\hh(\II)+\mu f(\II, \varphi)$$
where $\mu$ is a very small parameter;  $(\II, \varphi)=(\II_1, \ldots, \II_n, \varphi_1, \ldots, \varphi_n)$ are coordinates -- usually named {\it action--angle} -- taking values in $V\times {\mathbb T}^n$, with $V\subset {\mathbb R}^n$ open and connected and ${\mathbb T}={\mathbb R}/(2\pi {\mathbb Z})$. We then recover such lack of structure by assigning to each term in~\eqref{ovlH} a ``relative weight''. We make two main assumptions. The former is that the Keplerian term~\eqref{Kep} is much greater than the zero--average terms
\begin{equation}\label{weight1}\left\|-\frac{1}{2\Lambda^2}\right\|\gg \delta\,\|\widetilde{\rm H}_{\CC}\|\end{equation}
where $\|\cdot\|$ is some norm on  functions.   Under condition~\eqref{weight1}, and provided that all the functions have a holomorphic extension on some small complex domain, perturbation theory (see \cite{arnold63}) allows us to conjugate the Hamiltonian~\eqref{ovlH} to 

\begin{equation}\label{remainder}{\rm H}_{\rm 3b, \CC}({\rm R}, {\rm G}, \Lambda, {\rm r}, {\rm g}, \ell)=-\frac{1}{2\Lambda^2}+\delta\overline{\rm H}_{\CC}({\rm R}, {\rm G}, \Lambda, {\rm r}, {\rm g})+{\rm O}_{2}({\rm R}, {\rm G}, \Lambda, {\rm r}, {\rm g}, \ell)\end{equation}
where  ${\rm O}_{2}$ denotes a remainder term, depending on all coordinates.
Let us look at the system which is obtained when the remainder is neglected: for such a system, 
 the first term in~\eqref{ovlH} becomes an inessential additive term for the averaged Hamiltonian which,  without loss of generality (see next Section~\ref{Set up} for a discussion), we fix at
\begin{equation}\label{Lambda}\Lambda=1\,.\end{equation}
Reabsorbing the parameter $\delta$ through a change of time, we are reduced to study the 2--degrees of freedom Hamiltonian $\overline{\rm H}_{\rm C}$, which is given by
\begin{align}\label{ham}\overline{\rm H}_{\rm C}({\rm R}, {\rm G}, {\rm r}, {\rm g})=\frac{{\rm R}^2}{2}+\frac{({\rm C}-{\rm G})^2}{2{\rm r}^2}-\alpha{\rm U}({\rm r}, {\rm G}, {\rm g})
 -\frac{\beta}{{\rm r}}
 \end{align}
 where we have assumed  \begin{eqnarray}\label{angular momenta}\bG\parallel \bC\parallel(\bC-\bG)\parallel \bk\end{eqnarray} so that $\|{\mathbf x}'\times{\mathbf y}'\|=\|\bC-\bG\|={\rm C}-{\rm G}$, and we have denoted as

 \begin{align}\label{Usb}{\rm U}({\rm r}, {\rm G}, {\rm g}):=\frac{1}{2\pi}\int_0^{2\pi}\frac{d\ell}{\|{\mathbf x}'-{\mathbf x}\|
 }
 \end{align}
 the simply $\ell$--averaged\footnote{Here, ``simply'' is used as opposed to the more familiar ``doubly'' averaged Newtonian potential, most often encountered in the literature; e.g.~\cite{arnold63, laskarR95, fejoz04, pinzari-th09, chierchiaPi11b, chierchiaPi11c}.} of the Newtonian\footnote{We call ``Newtonian potential'' the function $\displaystyle \frac{1}{\|\bx-\bx'\|}$. Note that the term $\gamma\by'\cdot \by$ has zero--average (being $\by$ proportional to the $\ell$--derivative of $\bx$ and $\by'$ $\ell$--independent), so it is merged in $\widetilde{\rm H}_{\CC}$, together with the zero--average part $\widetilde \UU$ of the Newtonian potential.} potential (obviously, written using the above coordinates), which turns to be $\CC$--independent (see formulae~\eqref{U} below).

 \noindent
 In order to describe the motions we are looking for, we rewrite $\overline\HH_{\CC}$ as
\begin{equation}\label{hamNEW}\overline\HH_{\CC}=\KK_{\CC}({\rm R, r})-\alpha \UU({\rm r, G, g})+\widetilde\KK_{\CC}({\rm G, r})\end{equation}
where
\begin{eqnarray*}
	\KK_{\CC}({\rm R, r}) =   \frac{\rm R^2}{2}  + \frac{\rm C^2}{2 {\rm r}^2} -\frac{\beta}{\rm r}\,,\quad
	\widetilde\KK_{\CC}({\rm G, r}) = \frac{{\rm G^2-2 C G}}{2 {\rm r}^2} \, .\\
\end{eqnarray*}
We look at regions of phase space where
\begin{equation}\label{weights}\|\KK_{\CC}\|\gg \alpha\| \UU\|\gg \|\widetilde\KK_{\CC}\| \, .
\end{equation}
which is our second assumption.

\begin{figure}[H]
\centering
\begin{subfigure}[b]{0.3\textwidth}
\includegraphics[width=\textwidth]{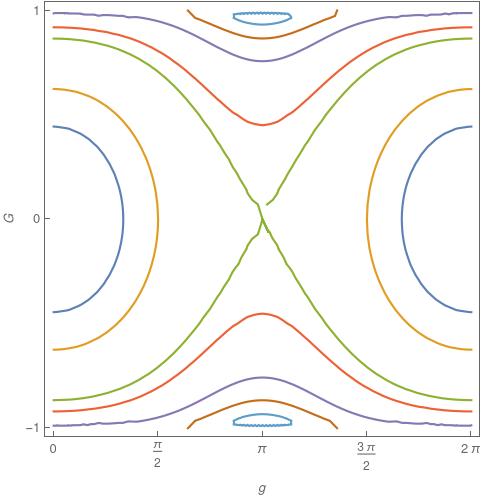}
\caption{$0<{\rm r}< 1$}
\end{subfigure}
\hfill
\begin{subfigure}[b]{0.3\textwidth}
\includegraphics[width=\textwidth]{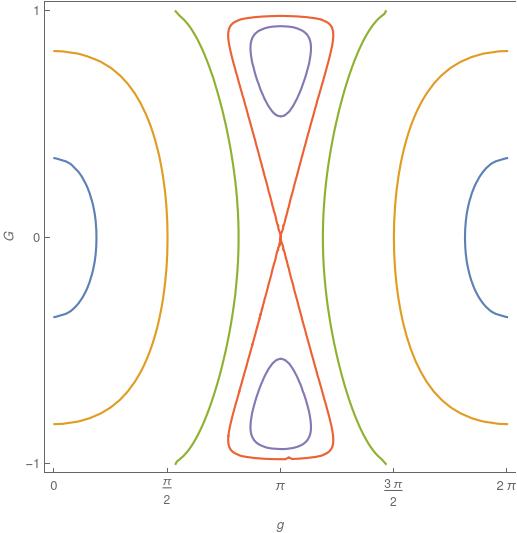}
\caption{$1<{\rm r}<2$}
\end{subfigure}
\hfill
\begin{subfigure}[b]{0.3\textwidth}
\includegraphics[width=\textwidth]{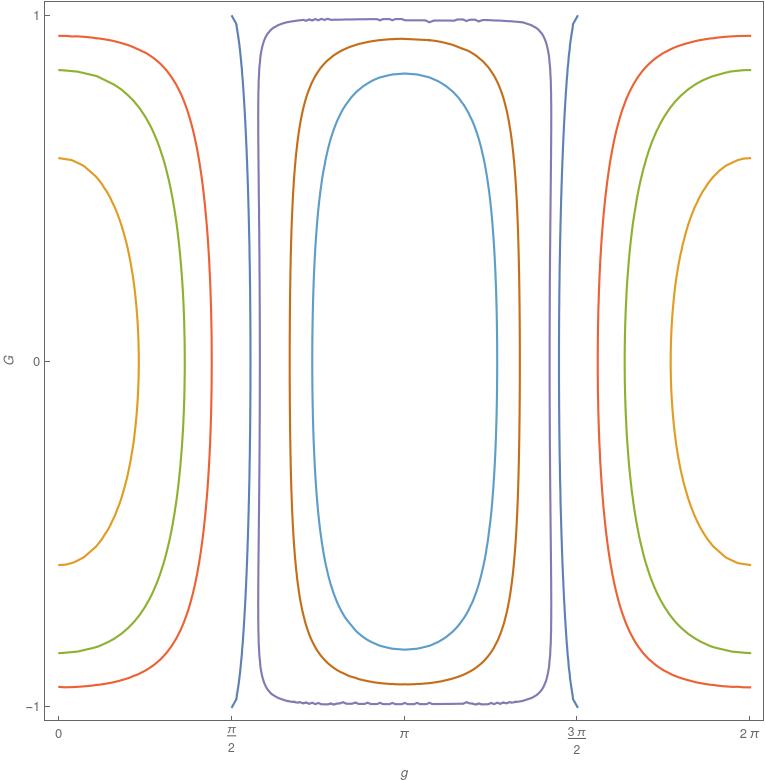}
\caption{${\rm r}>2$}
\end{subfigure}
\caption{Sections, at ${\rm r}$ fixed,  of the level surfaces of ${\rm E}$.  
}  
\label{figure1}
\end{figure}

\noindent 
Let us briefly comment on
inequalities~\eqref{weight1} and~\eqref{weights}. These inequalities aim 
to shape the Hamiltonian~\eqref{remainder} as a {\it three--scales} system, namely, composed, at a first order of approximation,  of three simpler terms of very different sizes. In particular, under such inequalities, one may argue that, at a first order of approximation, the motions $(\Lambda(t), {\rm R}(t), {\rm G}(t), \ell(t), {\rm r}(t), {\rm g}(t))$ of $\HH_{\rm 3b,\CC}$ are as follows:

\begin{conjecture}\label{picture of motion}\rm   
\item[--] $\Lambda(t)\sim1$ remains almost constant and $\ell(t)\sim t$ moves fast;
\item[--] the motion $({\rm R}(t), {\rm r}(t))$  is ruled by $\KK_{\CC}$;
\item[--] the motion $({\rm G}(t), {\rm g}(t))$  is ruled by the non--autonomous Hamiltonian $\UU({\rm r}(t), \cdot, \cdot)$.
$\quad\square$\end{conjecture}
Now, the term  ${\rm K}_{\rm C}$ is well--known. It consists of the one degree of freedom reduction of  a Hamiltonian as in~\eqref{Kep}, with  a fictitious angular momentum equal to $\CC$. The coordinate $\rm r$ moves as the length of a vector along a conic section (which can be an ellipse, parabola or hyperbola, according to the sign of the energy ${\rm K}_{\rm C}$) according to the Law of Equal Areas.
To understand  the dynamics generated by $\UU(\rm r, \cdot, \cdot)$, we need to  recall a property of such a function, pointed out in~\cite{pinzari19}. First of all, we remark that ${\rm U}$  is integrable. But the main point is that there exists a function ${\rm F}$ of two arguments such that
\begin{align}\label{relation***}{\rm U}({\rm r}, {\rm G}, {\rm g})={\rm F}(\rm r, {\rm E}({\rm r}, {\rm G}, {\rm g}))\end{align}
where
\begin{align}\label{euler}{\rm E}({\rm r}, {\rm G}, {\rm g})={\rm G}^2-{\rm r}\sqrt{1-{\rm G}^2}\cos{\rm g}\,.\end{align}
The function ${\rm E}({\rm r}, {\rm G}, {\rm g})$ above will be referred to as {\it Euler integral}, as it appears in the integration of the two--fixed centers Hamiltonian (also known as {\it Euler problem}). By~\eqref{relation***}, the level sets of ${\rm E}$, namely the curves \begin{align}\label{level curves}{\cal S}({\rm r}, {\cal E}):=\Big\{({\rm G}, {\rm g}):\quad {\rm G}^2-{\rm r}\sqrt{1-{\rm G}^2}\cos{\rm g}={\cal E}\Big\}\end{align} are also level sets of ${\rm U}$.  On the other hand, the phase portrait of ${\rm E}$  can be studied exactly, and this has been done in~\cite{pinzari20b}. We report the main results here. We fix a reference frame with ${\rm g}$ on the first axis, ${\rm G}$ on the second one. For the coordinates $({\rm g}, {\rm G})$, by the periodicity of ${\rm g}$, we consider a domain given by the rectangle $[0, 2\pi)\times (-1, 1)$. Then we have three cases.
  \begin{itemize}
\item[(a)]   $0<{\rm r}< 1$.  The point $(0, 0)$ is a minimum, while there are two symmetric maxima at $\big(\pi, \pm \sqrt{1-\frac{\rm {\rm r}^2}{4}}	\big)$ and one saddle at $(\pi, 0)$. The phase portrait includes two separatrices
\begin{eqnarray}\label{S0S1}\left\{\begin{array}{l}\displaystyle{\cal S}_0({\rm r})=\{({\rm G}, {\rm g}):\ {\rm E}({\rm r}, {\rm G}, {\rm g})={\rm r}\}
\\\\
\displaystyle{\cal S}_1({\rm r})=\{({\rm G}, {\rm g}):\ {\rm E}({\rm r}, {\rm G}, {\rm g})=1\}
\end{array}\right.
\end{eqnarray}
with ${\cal S}_0({\rm r})$ going through the saddle $(\pi, 0)$ and ${\cal S}_1({\rm r})$ through $(\frac{\pi}{2}, \pm 1 )$.
Rotational motions in between ${\cal S}_0({\rm r})$ and ${\cal S}_1({\rm r})$ do exist. ${\cal S}_0({\rm r})$ delimits librations about the minimum and rotations. ${\cal S}_1({\rm r})$ delimits different librations surrounding the maxima and the saddle and librational motions about the minimum.
 
 \item[(b)]   $1<{\rm r}< 2$. 
The minimum $(0, 0)$ persists, as well as the two symmetric maxima $\big(\pi, \pm \sqrt{1-\frac{\rm {\rm r}^2}{4}}	\big)$, the saddle at $(\pi, 0)$ and the separatrices~\eqref{S0S1}, with the difference, now, that ${\cal S}_1({\rm r})$ is inner with respect to ${\cal S}_0({\rm r})$, when looking from the minima. Rotational motions disappear, as in fact ${\cal S}_0({\rm r})$ delimits 
 librations about the  maxima and librations surrounding the maxima and the saddle, while
 ${\cal S}_1({\rm r})$ delimits different librations surrounding the maxima and the saddle and librational motions about the minimum.

\item[(c)]  ${\rm r}>2$. The saddle point and the separatrix ${\cal S}_0({\rm r})$ disappear, as and $(\pi, 0)$ turns to be a maximum, while $(\pi, 0)$ is still a minimum.
The phase portrait includes only the separatrix ${\cal S}_1(\rm r)$ in~\eqref{S0S1}, which delimits different librational motions about the minimum or the maximum.      
 \end{itemize}  
 The situation is represented in Figure~\ref{figure1}.

\begin{figure}[H]
\centering
\begin{subfigure}[b]{0.3\textwidth}
\includegraphics[width=\textwidth]{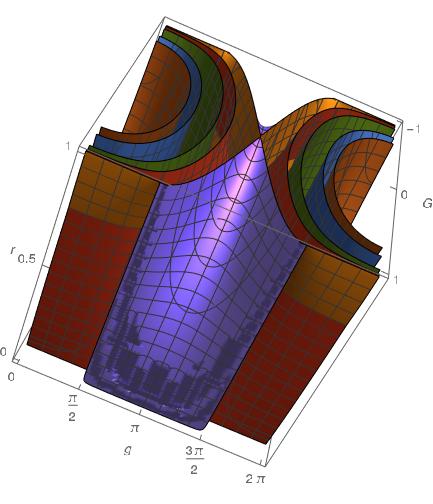}
\caption{$0<{\rm r}< 1$}
\end{subfigure}
\hfill
\begin{subfigure}[b]{0.3\textwidth}
\includegraphics[width=\textwidth]{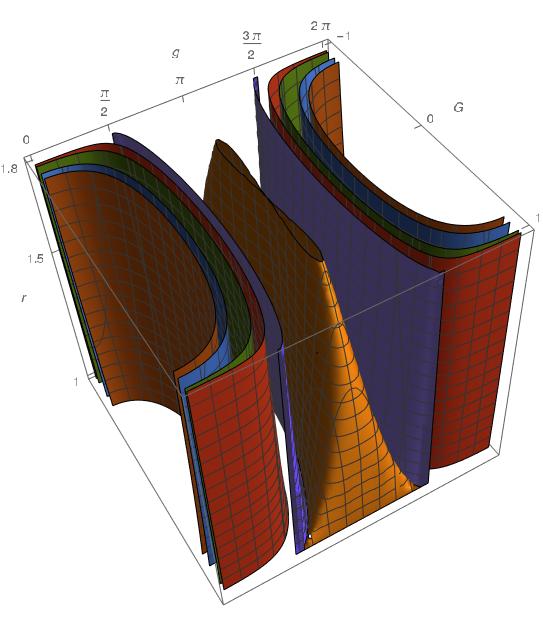}
\caption{$1<{\rm r}<2$}
\end{subfigure}
\hfill
\begin{subfigure}[b]{0.3\textwidth}
\includegraphics[width=\textwidth]{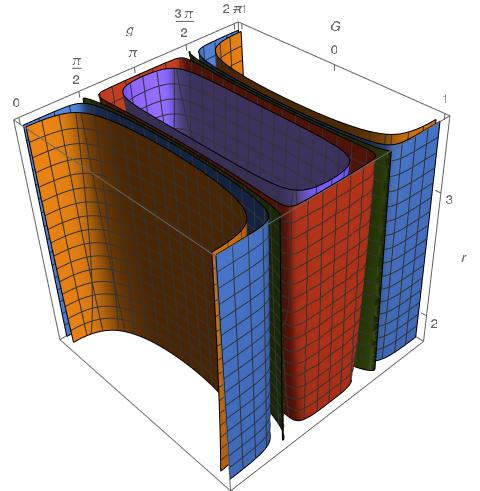}
\caption{${\rm r}>2$}
\end{subfigure}
\caption{Logs of the level surfaces of ${\rm E}$ in the space $({\rm r}, {\rm g}, {\rm G})$.  
}  
\label{figure2}
\end{figure}

\noindent
It is to be remarked, however, that the coordinate $\rm r$ stays fixed under $\rm E$, while it moves under $\overline\HH_\CC$. Therefore,  three--dimensional plots representing the manifolds corresponding to the ``lifted level sets''
\begin{equation}\label{eul_ell_surf}{\cal M}({\cal E})=\{({\rm r}, {\rm G}, {\rm g}):\ {\rm E}({\rm r}, {\rm G}, {\rm g})={\cal E}\}
\end{equation}
Such manifolds are represented in Figure~\ref{figure2}.

\begin{figure}[H]
\centering
\begin{subfigure}[b]{0.3\textwidth}
\includegraphics[width=\textwidth]{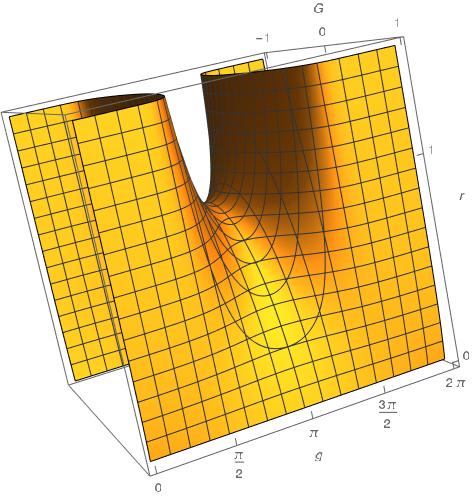}
\caption{$0<{\cal E}< 1$}
\end{subfigure}
\hfill
\begin{subfigure}[b]{0.3\textwidth}
\includegraphics[width=\textwidth]{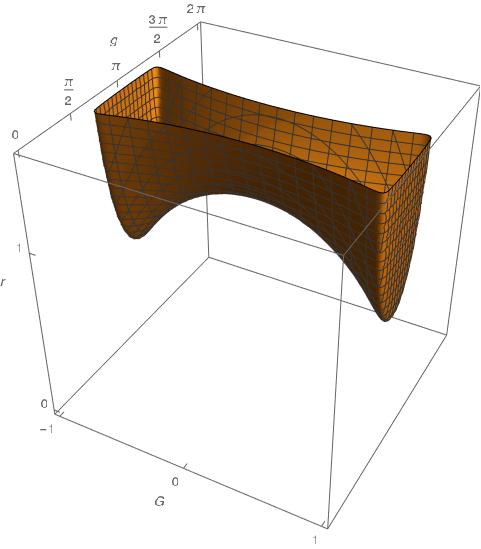}
\caption{$1<{\cal E}< 2$}
\end{subfigure}
\hfill
\begin{subfigure}[b]{0.3\textwidth}
\includegraphics[width=\textwidth]{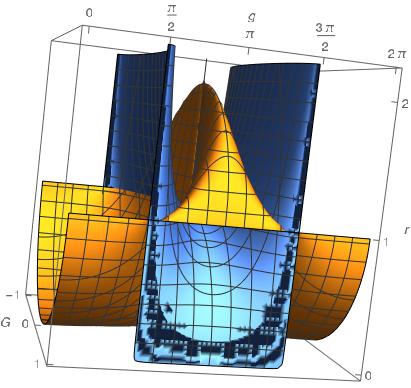}
\caption{}
\end{subfigure}
\caption{(a) and (b): The saddle point of ${\cal M}(\cal E)$; (c): ${\cal M}_0$ (yellow)  and ${\cal M}_1$ (blue).
}  
\label{figure3}
\end{figure}

\noindent
Each manifold ${\cal M}({\cal E})$ with $0<{\cal E}<2$ has a saddle at 
\begin{eqnarray}\label{saddles}({\rm r}_{\rm sad}, {\rm G}_{\rm sad}, {\rm g}_{\rm sad})=({\cal E}, 0, \pi)\qquad \forall\ 0<{\cal E}<2\,.\end{eqnarray}
The manifolds obtained
 ``lifting'' along the $\rm r$--direction the  curves ${\cal S}_0(\rm r)$, ${\cal S}_1(\rm r)$ in~\eqref{S0S1} will be denoted as
 \begin{eqnarray}\label{M0M1}\left\{\begin{array}{l}\displaystyle{\cal M}_0=\{({\rm r}, {\rm G}, {\rm g}):\ {\rm E}({\rm r}, {\rm G}, {\rm g})={\rm r}\}
\\\\
\displaystyle{\cal M}_1=\{({\rm r}, {\rm G}, {\rm g}):\ {\rm E}({\rm r}, {\rm G}, {\rm g})=1\}
\end{array}\right. \, .
\end{eqnarray}

\noindent
See Figure~\ref{figure3}.

\noindent
Combining the phase portraits above with Conjecture~\ref{picture of motion}, we pose the following
\begin{conjecture}\label{conjOLD}\rm
For a set of parameters and in a region of phase space where conditions~\eqref{weight1},~\eqref{weights} are verified,
\begin{itemize}
\item[i)]  the manifolds ${\cal M}({\cal E})$ are ``approximate invariant manifolds'' for  the Hamiltonian ${\rm H}_{\rm 3b, \CC}$ in~\eqref{ovlH}, at least for the time that $\rm r$ does not leave a fixed region (a), (b), (c) above; 
\item[ii)] the coupling between $\KK_\CC$ and $\UU$, the presence of the   ``disturbing term'' $\widetilde\KK_\CC$ and the remainder ${\rm O}_2$ in~\eqref{remainder}  are source of chaotic dynamics for ${\rm H}_{\rm 3b, \CC}$ in~\eqref{ovlH}, closely to ${\cal M}_0$. $\quad\square$
\end{itemize}
\end{conjecture}

\noindent
However, our numerical explorations will only support the following assertion.
\begin{conjecture}\label{conj}\rm
For a set of parameters and in a region of phase space where condition~\eqref{weights} is verified,
\begin{itemize}
\item[i)$'$]  the manifolds ${\cal M}({\cal E})$ are ``approximate invariant manifolds'' for  the Hamiltonian $\overline{\rm H}_{\CC}$ in~\eqref{ham}, at least for the time that $\rm r$ does not leave a fixed region (a), (b), (c) above; 
\item[ii)$'$] the coupling between $\KK_\CC$ and $\UU$ and the presence of the   ``disturbing term'' $\widetilde\KK_\CC$ are source of chaotic dynamics for $\overline{\rm H}_{\CC}$ in~\eqref{ham}, closely to ${\cal M}_0$.
$\quad\square$
\end{itemize}
\end{conjecture}

\noindent
Note that Conjecture~\ref{conj} is based on~\eqref{weights}, but does not need~\eqref{weight1}. This is precisely the reason that led us to relax Conjecture~\ref{conjOLD} to the form~\ref{conj}. Let us briefly comment on this.
\\
A typical difficulty in Celestial Mechanics is represented by the {\it lack of parameters}. A famous example goes back to V. I. Arnold, who, in the  paper~\cite{arnold63},  wanted to regard the $(1+n)$--body problem (in the planetary\footnote{The planetary $(1+n)$--body problem consists of the Newtonian attraction of $1+n$ masses $m_0$, $\ldots$, $m_n$, where $m_1$, $\ldots$, $m_n$ (``planets'') have comparable sizes, but much smaller than $m_0$ (``star'', or ``sun'').} version) as close to $n$ {\it independent}
Kepler Hamiltonians~\eqref{Kep}.
He had at his disposal only {\it one} parameter, given by the maximum ratio $\mu$ of the planets' masses to the sun's. In a very similar situation 
as for the Hamiltonian~\eqref{ovlH}, where the Keplerian approximation provides  motions for only the $(\Lambda, \ell)'$s coordinates,
using a {\it two--scale} approximation (a scale ``of order 1'' for the motions of the $\ell's$; a scale ``of order $\mu$''  for the motions of the ellipses), he ingeniously found a good approximation for the motions of {\it all} coordinates. To fulfil\footnote{The statement in~\cite{arnold63} has been completely proved in~\cite{fejoz04}. The study has been reconsidered in~\cite{chierchiaPi11b} for open problems after~\cite{arnold63, fejoz04}.} this, he required, besides the smallness of the parameter $\mu$, an additive condition (i.e., the smallness of eccentricities and inclinations of the planets' instantaneous ellipses of the planets) having the r\^ole of pushing away, in the Hamiltonian, remainder terms from the two leading scales terms. Now, inequalities~\eqref{weight1} and~\eqref{weights} have the precise scope of emulating Arnold's strategy, with the difference that, in our case, they provide a {\it three--scale} system. In particular,
\eqref{weight1} stresses that the velocity of $\ell$ is much larger than  the velocities of
$\rm(R, r)$ and $\rm(G, g)$, in turn separated by~\eqref{weights}.
Of course, the Hamiltonian $\overline\HH_\CC$ has a physical meaning only whenever~\eqref{weight1} is satisfied. 
However, what is, if existing, a ``natural'' choice of parameter masses and/or of additive conditions that make~\eqref{weight1}--\eqref{weights} true does not seem immediate to us. For this reason, we choose to investigate the motions of $\overline\HH_\CC$ independently whether condition~\eqref{weight1} is verified or not. Our interest in $\overline\HH_\CC$ is indeed precisely related to the Euler integral \eqref{euler}: we aim to find zones in the phase space of $\overline\HH_\CC$ where $\EE$ affords slow variations and, simultaneously, chaos is detected. A similar point of view has been taken up, on the other hand,  in the published papers~\cite{pinzari20a, chenPi2021, diruzzaDP20}.

\noindent
Before switching to technical parts, we recall related works, with no aim of completeness.
Chaos in many--body systems has been widely studied in the literature~\cite{delshamsKDS2019, guardiaMS2016, bolotin2006c,  fejozG2016, fejozGKR14, guardiaKZ19}. 
 For general information on chaotic phenomena, the reader may consult
~\cite{delshamsDS08, legaGuF2016, delshamsDLS00}. 
Closely related papers to the current one are the aforementioned~\cite{pinzari20a, chenPi2021, diruzzaDP20}. Specifically, in~\cite{pinzari20a} Conjecture~\ref{conj} has been proved in the case (c), while in~\cite{chenPi2021} it has been proved in the case (a), with $\rm r\ll 1$ and for motions very close to ${\cal S}_0(\rm r)$. Both such papers are rigorous proofs and are based well adapted normal form theory, so they unavoidably deal with ideal situations, where ``ideal'' means that the estimates on parameters are far from being optimal. In~\cite{diruzzaDP20} the case (c) has been reconsidered from the numerical point of view and the existence of chaotic motions among librations has been pointed out.

\noindent This paper is organised as follows.

\begin{itemize}
\item[--] In Section~\ref{Set up} we collect useful properties for the Hamiltonian~\eqref{ham}.
\item[--] In Section~\ref{Some remarks} we discuss conditions~\eqref{weight1} and~\eqref{weights} on two concrete examples.
\item[--] In Sections~\ref{eul_integral} and~\ref{Neighbourhoods} we focus on one of the examples and study the phase space around two orbits. In particular, we study the variations of the function~\eqref{euler} around one orbit which spends much time closely to  the saddle point~\eqref{saddles} of one of the manifolds~\eqref{eul_ell_surf}.
\item[--] In Section~\ref{Symbolic dynamics} we show the existence of chaos (in fact, of symbolic dynamics) in the region of the mentioned saddle point. This is the main result of the paper.
\item[--] In Section~\ref{Control of errors and conclusions} we discuss how we control numerical errors, draw some conclusions and foresee possible future works.
	\end{itemize}

\section{Facts to be known}\label{Set up}
Let us consider the Hamiltonian $\HH_{\rm 3b}$ in~\eqref{ham_4dof}, with $x_3=x_3'=y_3=y_3'=0$.  We define a canonical change of coordinates which reduces the invariance of $\HH_{\rm 3b}$ by rotations, via a canonical transformation \begin{eqnarray}\label{change}(\by', \by, \bx', \bx)\in {\mathbb R}^2\times {\mathbb R}^2\times{\mathbb R}^2\times{\mathbb R}^2\setminus\Delta\to (\CC, {\rm G}, \Lambda, {\rm R}, {\rm c}, {\rm g}, \ell, {\rm r})\in {\mathbb R}_+^3\times {\mathbb R}\times{\mathbb T}^3\times {\mathbb R}_+\end{eqnarray}
where
$$\Delta=\Big\{\bx={\mathbf 0}\Big\}\bigcup\Big\{ \bx'={\mathbf 0}\Big\}\bigcup\Big\{\bx-\bx'={\mathbf 0}\Big\}$$
is the ``collision set''. To define the new coordinates at right hand side of~\eqref{change}, we  denote as

 \begin{itemize}
   \item[{\tiny\textbullet}]  ${\mathbf i}=\left(
\begin{array}{lll}
1\\
0\\
0
\end{array}
\right)\, $, $\ {\mathbf j}=\left(
\begin{array}{lll}
0\\
1\\
0
\end{array}
\right)$ the directions of  a hortonormal frame in ${\mathbb R}^2\times\{\mathbf 0\}$ and ${\mathbf k}={\mathbf i}\times {\mathbf j}$ (``$\times$'' denoting, as usual, the ``skew--product''). We assume~\eqref{angular momenta}. 
  \item[{\tiny\textbullet}] after 
 fixing  a set of values of $({\mathbf y},  {\mathbf x})$ where the Kepler Hamiltonian~\eqref{Kep}
 takes negative values, let ${\mathbb E}$ denote the elliptic orbit with initial values $({\mathbf y}_0, {\mathbf x}_0)$ in such set;
  \item[{\tiny\textbullet}]  ${\mathbf P}$, with $\|{\mathbf P}\|=1$,  the direction of the perihelion  of ${\mathbb E}$, assuming  ${\mathbb E}$ is not a circle;
     \item[{\tiny\textbullet}]    $\alpha_{\mathbf w}({\mathbf u}, {\mathbf v})$  the oriented angle  from ${\mathbf u}$ to ${\mathbf v}$ relatively to the positive orientation established by ${\mathbf w}$, if ${\mathbf u}$, ${\mathbf v}$ and ${\mathbf w}\in {\mathbb R}^3\setminus\{\mathbf 0\}$, with ${\mathbf u}$, ${\mathbf v}\perp{\mathbf w}$.
    \end{itemize}
Then the coordinates at the right hand side of~\eqref{change} are defined via

 \begin{align}\label{coord}
\left\{\begin{array}{l}\displaystyle {\rm C}=\|{\mathbf x}\times {\mathbf y}+{\mathbf x}'\times {\mathbf y}'\|\\\\
\displaystyle {\rm G}=\|{\mathbf x}\times {\mathbf y}\|\\\\
\displaystyle {\rm R}=\frac{\mathbf y'\cdot \mathbf x'}{\|\mathbf x'\|}\\\\
\displaystyle \Lambda= \sqrt{ a}
\end{array}\right.\qquad\qquad \left\{\begin{array}{l}\displaystyle {\rm c} =\alpha_{\mathbf k}(\mathbf i, \mathbf x')
\\\\
\displaystyle  {\rm g}=\alpha_{\mathbf k}({\mathbf x'},\mathbf P)
\\\\
\displaystyle  {\rm r}=\|\mathbf x'\|\\\\
\displaystyle \ell={\rm mean\ anomaly\ of\ {\mathbf x}\ in\ \mathbb E}
\end{array}\right.
\end{align}
We recall that the mean anomaly of $\bx$ is defined as  as the area of the elliptic sector spanned by $\bx$ relatively to the perihelion of $\mathbb E$, normalised to $2\pi$ (refer to Figure~\ref{model}).   \\
With a language which goes back to Liouville--Arnold theorem, the coordinates $\CC$, ${\rm G}$ and $\Lambda$ will be called ``actions'', for  being conjugated to ${\rm c}$, $\rm g$ and $\ell$, which take values in $\mathbb T$, hence, are called ``angles''. The coordinates~\eqref{coord} are singular when ${\rm G}=\Lambda$ (corresponding to vanishing eccentricity of $\mathbb E$. In that case, $\mathbf P$ is not defined) or $\rm r=0$ (as $\rm c$ is not defined), so we should safely exclude such values from our domain.  Observe however that the Hamiltonian~\eqref{ham_4dof} is $\rm c$--independent by its discussed SO(2) invariance, and the singularity at $\rm G=\Lambda$ could be -- if needed --  easily eliminated switching to the ``Poincar\'e'' transformation $(\Lambda, \rm G, \ell, \rm g)\to (\Lambda, p, \lambda, q)=(\Lambda, \sqrt{2(\Lambda-\rm G)}\cos{\rm g}, \ell+\rm g, -\sqrt{2(\Lambda-\rm G)}\sin{\rm g})$.
The canonical character of the coordinates~\eqref{coord} has been discussed, in a more general setting, in~\cite{pinzari19}. 
Using the coordinates~\eqref{coord}, the Hamiltonian $\HH_{\rm 3b}$ turns to be ${\rm c}$--independent, as  the action $\CC$ is a first integral for it. 
Then, we regard it as a ``fixed parameter'', skipping it from actions.
Another first integral, namely the action $\Lambda$, appears when
taking the $\ell$--average of~\eqref{ham_4dof}, as discussed in the previous section.   In order to further simplify the discussion, it turns to be useful to  remark the following scaling property. Switching to the a--dimensional and canonical coordinates 
 $$\widehat{\rm R}:={\rm R}{\Lambda}\,,\quad\widehat{\rm G}:=\frac{\rm G}{\Lambda}\,,\quad \widehat{\rm r}:=\frac{\rm r}{\Lambda^2}\,,\quad \widehat{\rm g}:={\rm g}$$
(possible because $\Lambda$ is a ``parameter'') one has the following identities
\begin{eqnarray}\label{scaling}
&&{\rm H}_{{\rm C}, \Lambda}({\rm R}, {\rm G}, {\rm r}, {\rm g})={\Lambda^{-2}}{\rm H}_{\widehat{\rm C}, 1}(\widehat{\rm R}, \widehat{\rm G}, \widehat{\rm r}, \widehat {\rm g})\nonumber\\\nonumber\\
&& {\rm U}_\Lambda({\rm r},{\rm G},  {\rm g})={\Lambda^{-2}}{\rm U}_{ 1}(\widehat{\rm r}, \widehat{\rm G}, \widehat {\rm g})\nonumber\\\nonumber\\
&& {\rm E}_\Lambda({\rm r},{\rm G},  {\rm g})={\Lambda^2}{\rm E}_{ 1}(\widehat{\rm r}, \widehat{\rm G}, \widehat {\rm g})
\end{eqnarray}
with $\widehat{\rm C}$ being the ratio $\frac{\rm C}{\Lambda}$.
The equalities in~\eqref{scaling} allow us to assume~\eqref{Lambda} once forever and eliminate the ``hats'' and subfixes ${}_{1}$.  
As a result, $\overline\HH_{\rm C}$ depends on   $3$ parameters only, namely
$\alpha$, $\beta$ and ${\rm C}$, and  is  reduced to
2 degrees of freedom, ruled by the coordinates $(\rm R, \rm G, \rm r, \rm g)$.
We  provide the explicit expression of $\UU$, under the choice~\eqref{Lambda}.
 Using, alternatively, the true anomaly $\nu$ and the eccentric anomaly $\xi$, we have

\begin{eqnarray}\label{U}
	\UU(\rm r, G, g)& = & \frac{{\rm G}^3}{2 \pi} \int_0^{2 \pi}  \frac{d \nu}{(1 + e \cos \nu)\sqrt{{\rm r}^2(1 + e \cos \nu)^2  - 2{\rm G}^2{\rm r}(1 + e \cos \nu)  \cos ({\rm g}+\nu)+{\rm G}^4}}\nonumber\\
	&=&\frac{1}{2 \pi} \int_0^{2 \pi}  \frac{(1-e\cos\xi)d \xi}{\sqrt{(1 - e \cos \xi)^2  - 2{\rm r}(\cos\xi-e)\cos{\rm g}+2{\rm r}{\rm G} \sin\xi\sin{\rm g}+{\rm r}^2}}
\end{eqnarray}
 with
$$e=\sqrt{1-{\rm G}^2}$$
being the eccentricity. 

\noindent
As a consequence of relation~\eqref{relation***} and as $\EE$ depends on $\rm g$ only via its cosinus while the other terms in~\eqref{ham} do not depend on $\rm g$, we remark the following symmetry:
\begin{proposition}
The Hamiltonian~\eqref{ham} does not change replacing $\rm g$ with $2k\pi-\rm g$, $k\in \mathbb Z$.
\end{proposition}

\noindent
In fact, this symmetry reflects in all orbits of $\overline\HH_\CC$; see, e.g., the orbits $\Gamma_{\rm s}$, $\Gamma_{\rm u}$ mentioned in Section~\ref{eul_integral}.

\noindent
As, in our experiments, we are   going to consider a global region of phase space, we need to establish the singularities of $\UU$. Below, we shall briefly show that
\begin{proposition}[\cite{pinzari20a, pinzari20b,chenPi2021}]	\label{prop: sing}
 The function $\UU$ is singular if and only if \,$0<\rm r<2$ and  $({\rm G}, {\rm g})\in{\cal S}_0(\rm r)$.
 \end{proposition}
Namely, the manifold ${\cal S}_0({\rm r})$  looses its meaning of saddle separatrix in the Hamiltonian~\eqref{ham} (discussed in the previous section) to gain the title of ``singular manifold''. 
 In~\cite{chenPi2021} the rate of divergence of $\UU$  has been established to be logarithmic, with respect to the distance from ${\cal S}_0(\rm r)$.
 
 \noindent
In this paper, we focus on a region of phase space where $0<\rm r<1$, so as to deal with the respective cases (a) in Figures~\ref{figure1}--\ref{figure2}. Ideally, we would be tempted to perform computations by replacing $\UU$ with a polynomial 
\begin{eqnarray}\label{UN}\UU_N({\rm r, G, g})=  \sum_{n=0}^{N} P_{n}({\rm G, g}) \cdot  {\rm r}^{n}\end{eqnarray}
 with sufficiently high degree $N$, provided to keep at a finite distance from ${\cal S}_0(\rm r)$. However, in this expansion  the coefficients $P_{n}({\rm G, g})$ are proportional to   negative powers of $\rm G$, as one immediately recognises from~\eqref{U}. This means that regions in phase space with very small values of $\rm G$ would not be covered by such an approximation, while we precisely aim to look at such regions. On the other hand,  by Proposition~\ref{prop: sing}, $\rm G=0$ is not a singularity, if ${\rm g}\neq \pi$ (as $(0, \pi)$ is the only point of ${\cal S}_0(\rm r)$ with ${\rm G}=0$). Therefore, instead of~\eqref{UN}, we consider a ``renormalised'' expansion of the fom   
 \begin{equation}
\label{funcion_U}
\UU_N({\rm r, G, g})=  \sum_{n=0}^{N} Q_{n}({\rm r, G, g}) \cdot  {\rm r}^{n}
\end{equation}
which differs from~\eqref{UN} by orders of ${\rm r}^{-N-1}$. The expansion~\eqref{funcion_U} is possible because of the relation~\eqref{relation***}. Indeed, by such relation, $\UU$ depends on $({\rm G}, {\rm g})$ only via ${\rm E}(\rm r, G, g)$. Therefore, picking up, for any fixed level set~\eqref{level curves} with ${\cal E}\geq 0$, the point of ${\rm E}$ with coordinates $(\rm G, g)=(\sqrt{\cal E}, \frac{\pi}{2})$,  we have the identity
$$
{\rm F}({\rm r}, {\cal E})=\UU\left(\rm r, {\sqrt{\cal E}}, \frac{\pi}{2}\right)
$$
 This identity reflects in the  expansion~\eqref{UN}, providing the  expansion~\eqref{funcion_U}, with
 $$Q_n({\rm r, G, g})=P_n\left(\sqrt{{\rm E}({\rm r, G, g})}, \frac{\pi}{2}\right)\,,\qquad \forall\ \  ({\rm r, G, g}):\ {\rm E}({\rm r, G, g})\geq 0$$ From the procedural point of view, we remark that in the expansion~\eqref{funcion_U} only the terms with $n=2k$ even survive, as (as one readily sees using, e.g., a Legendre polynomials expansion) the function in~\eqref{U} is even in $\rm r$ when ${\rm g}=\frac{\pi}{2}$.

 \noindent
 We conclude this section with the

\noindent
{\bf Proof of Proposition~\ref{prop: sing}} The first expression in~\eqref{U} shows, for the function under the integral,  a pole of order $1$ 
(corresponding to the zero of the expression under the square root, and understood a collision between $\bx$ and $\bx'$) when the following equalities are satisfied
\begin{equation}\label{singularities}\nu=-{\rm g}\quad{\rm mod}\ 2\pi\,,\qquad {\rm r}(1+e\cos{\rm g})={\rm G}^2\end{equation}
and a pole of order $2$ at $\nu=\pi$ when ${\rm G}=0$ (corresponding to the zero of $(1 + e \cos \nu)$ and understood a collision between $\bx$ and the unit mass).
Observe that the second equation in~\eqref{singularities} is nothing else than the equation of ${\cal S}_0({\rm r})$. The singularity at ${\rm G}=0$ leaves instead $\UU$ perfectly regular, as the second expression in 
\eqref{U} gives
$$\UU({\rm r}, 0, {\rm g})=\frac{1}{2 \pi} \int_0^{2 \pi}  \frac{(1-\cos\xi)d \xi}{\sqrt{(1 -  \cos \xi)^2  + 2{\rm r}(1-\cos\xi)\cos{\rm g}+{\rm r}^2}}$$
This shows that the only possibility of singularity for $\UU({\rm r}, 0, {\rm g})$ is 
when $0<{\rm r}<2$ and ${\rm g}=\pi$. But this is already counted in ${\cal S}_0(\rm r)$. $\quad \square$

 \section{Discussion of~\eqref{weight1} and~\eqref{weights} on two examples}\label{Some remarks}
 
In the introduction,  we mentioned that the simultaneous fulfilment of  inequalities~\eqref{weight1} and~\eqref{weights} does depend {\it only} by choice of the parameters of the system -- in our case,  $\alpha$, $\beta$ and $\CC$, but also needs a careful choice of the  phase space. In agreement with the numeric nature of the paper,  in this section we investigate  the question on two concrete examples. We pick two triples of values for $\alpha$, $\beta$ and $\CC$, and, for each triple, we consider motions of different kind. We check that, while the inequality~\eqref{weights} is met along all the  orbits under examination, unfortunately,~\eqref{weight1} is not.

\subsection*{Example 1}

In the first example, we take\footnote{We recall that $\alpha, \beta$ are the two independent mass parameters and they are uniquely linked to $\kappa, \mu$ which turn out to be  $\kappa=0.01787503$, $\mu=0.02153618$.}
 \begin{equation}
 \label{parameter_2}
 	  \alpha = 100 \, ,\quad \beta = 120 \, , \quad \CC = 10 \,.
 \end{equation}
In order to check~\eqref{weights} we consider two orbits, $\Gamma_1$ and $\Gamma_2$, of the Hamiltonian $\overline\HH_\CC$ in~\eqref{ham} on two {\it different} energy levels, ${\cal H}_1$, ${\cal H}_2$,  of  $\overline\HH_\CC$, but with initial data
 $({\rm R}_i, {\rm G}_i,{\rm r}_i,   {\rm g}_i)$ chosen so that the triplets  $({\rm r}_i, {\rm G}_i, {\rm g}_i)$ 
coincide with the saddle points~\eqref{saddles} of the manifolds ${\cal M}({\cal E}_i)$, with 
 $${\cal E}_i=\EE({\rm r}_i, {\rm G}_i, {\rm g}_i)=\left\{
\begin{array}{lll}
0.2\quad&{\rm if}\quad &i=1\\\\
0.8\bar 3&{\rm if} &i=2
\end{array}
\right. \, .$$ 
 In fact, we take 
  \begin{eqnarray}
  \label{in_cond_2}
 &&{\cal H}_1:=613.75\, ,\qquad	\Gamma_1:\ \left\{
 	\begin{array}{l}
 		\rm R_{1} = 0 \\
 		\rm G_1 = 0 \\
 		\rm r_1 =0.2 \\
 		\rm g_1 = \pi
	\end{array}
 \right. \, ,\\\nonumber\\\nonumber\\
 &&{\cal H}_2:=-155.025\, ,\qquad	\Gamma_2:\ \left\{
 	\begin{array}{l}
 		\rm R_{2} = 0 \\
 		\rm G_2 = 0 \\
 		\rm r_2 =0.8\bar3 \\
 		\rm g_2 = \pi
	\end{array}
 \right. \, .
 \end{eqnarray}
 We  remark that 
${\cal E}_1$ has been chosen so that $\KK_\CC$ is initially positive, while
${\cal E}_2$ has been chosen so that $\KK_\CC$ is initially negative (in fact, at its\footnote{As well known, $\KK_\CC$ attains its minimum, given by $-\frac{\beta^2}{2\CC^2}$, when $\rm R=0$ and
$${\rm r}=\frac{\CC^2}{\beta}=\frac{100}{120}=0.8\bar 3\,.$$} minimum). \\
In order to check~\eqref{weight1}, we consider the orbits $\widetilde\Gamma_1$ and $\widetilde\Gamma_2$, of the whole Hamiltonian $\HH_{\rm 3b, \CC}$ in~\eqref{ovlH}, departing from the  initial data obtained completing the respective initial values of $\Gamma_1$ and $\Gamma_2$ with $\Lambda_{i}=1$ (as prescribed in~\eqref{Lambda}) and $\ell_{i}=\pi$.
 \begin{figure}[H]
\centering
\includegraphics[width=7.5cm,height=5cm,draft=false]{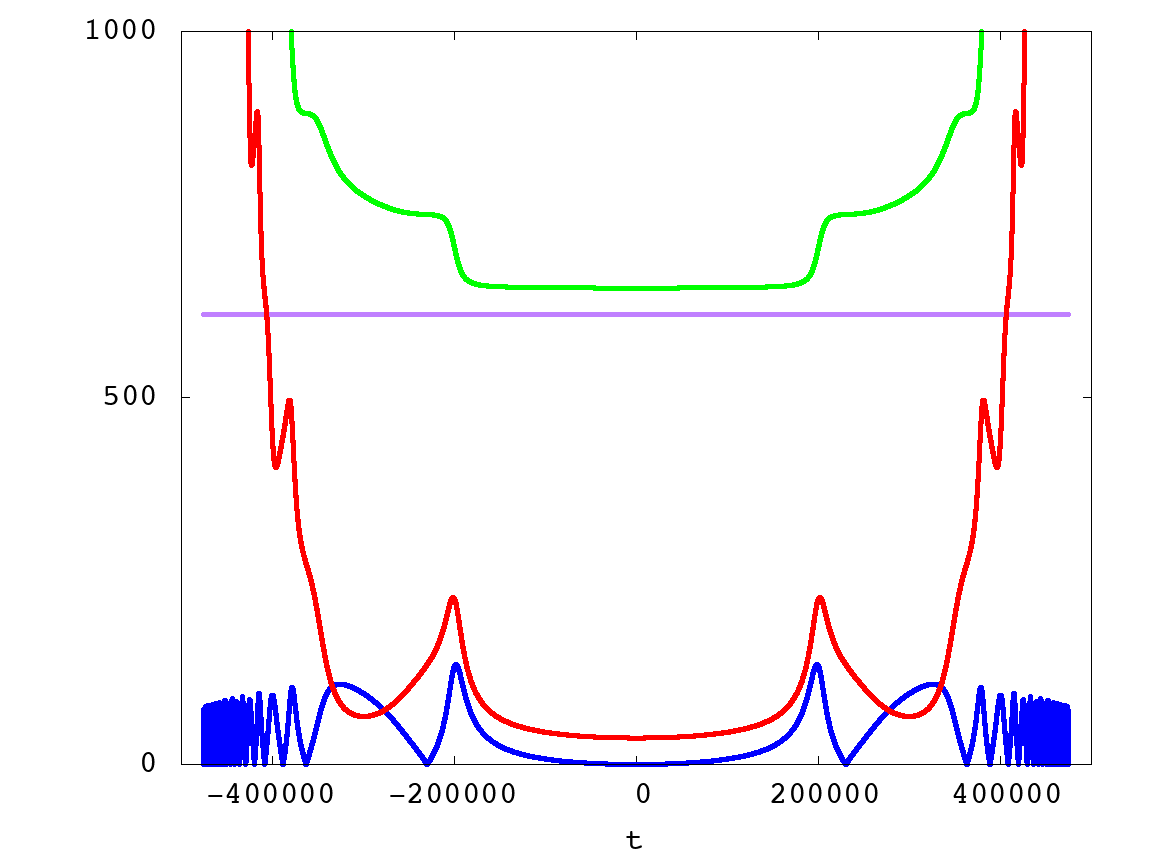}
\includegraphics[width=7.5cm,height=5cm,draft=false]{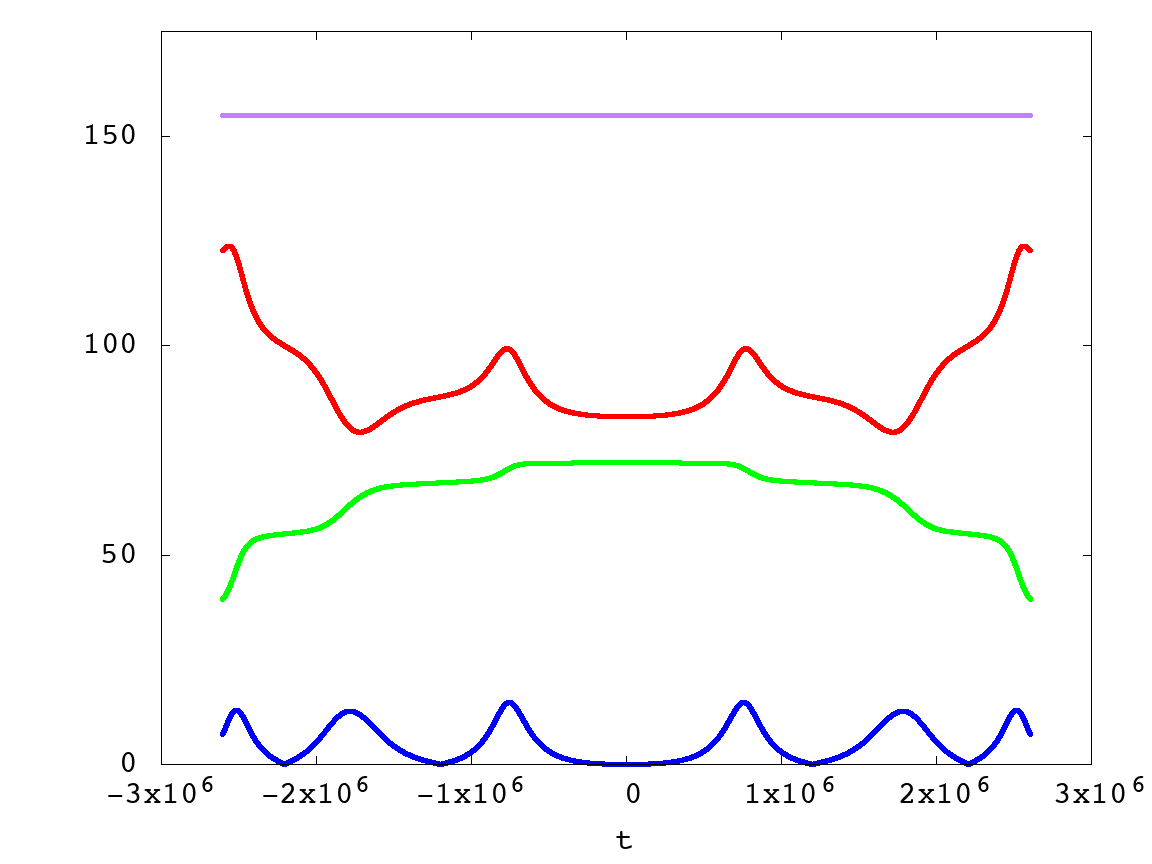}
	\caption{Graphs of the absolute values of  $\KK_\CC$ (green), $\alpha\UU$ (red) and $\widetilde \KK_\CC$ (blue) along $\Gamma_1$ (left) and $\Gamma_2$ (right). The purple line, representing the total energy~\eqref{ham}, is reported for comparison. Incidentally, plotting the total energy is a well known useful tool to check the correctness of numerical integrations}. \label{energia2} 
\end{figure}

 \begin{figure}[H]
\centering
\includegraphics[width=7.5cm,height=5cm,draft=false]{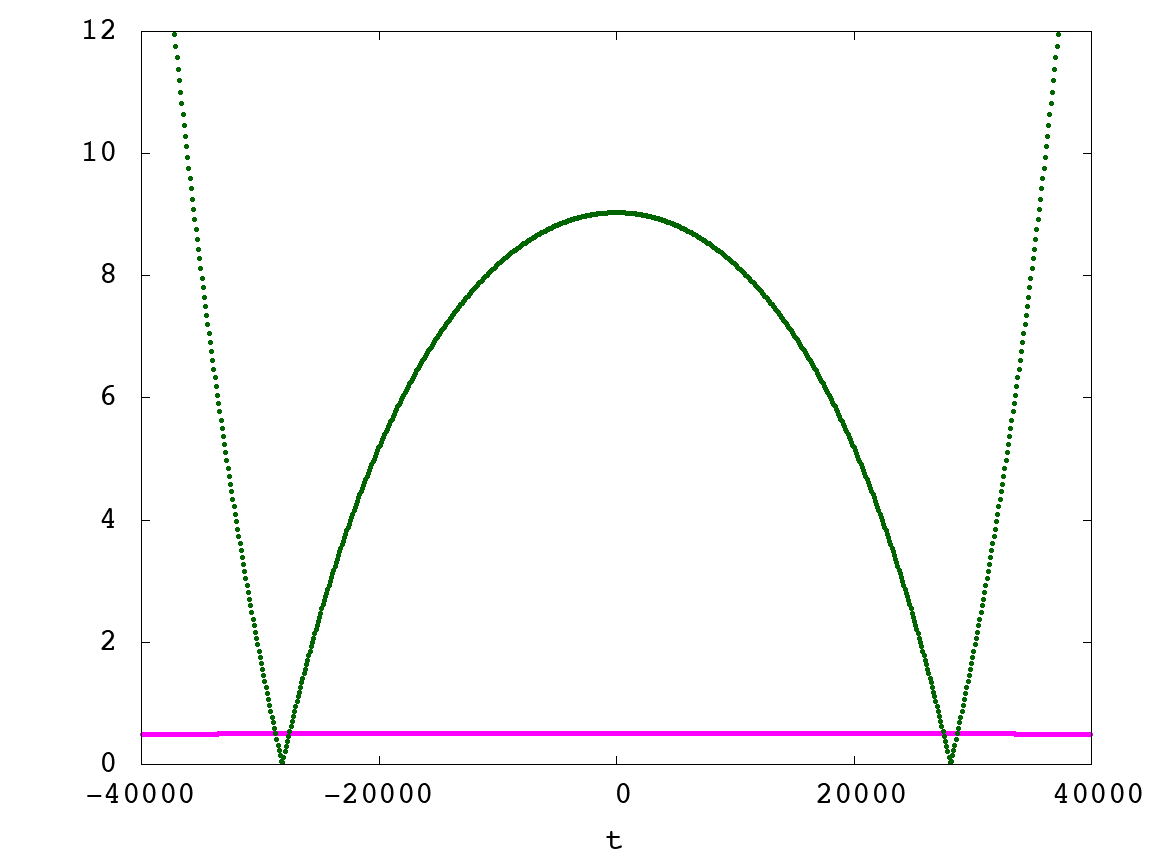}
\includegraphics[width=7.5cm,height=5cm,draft=false]{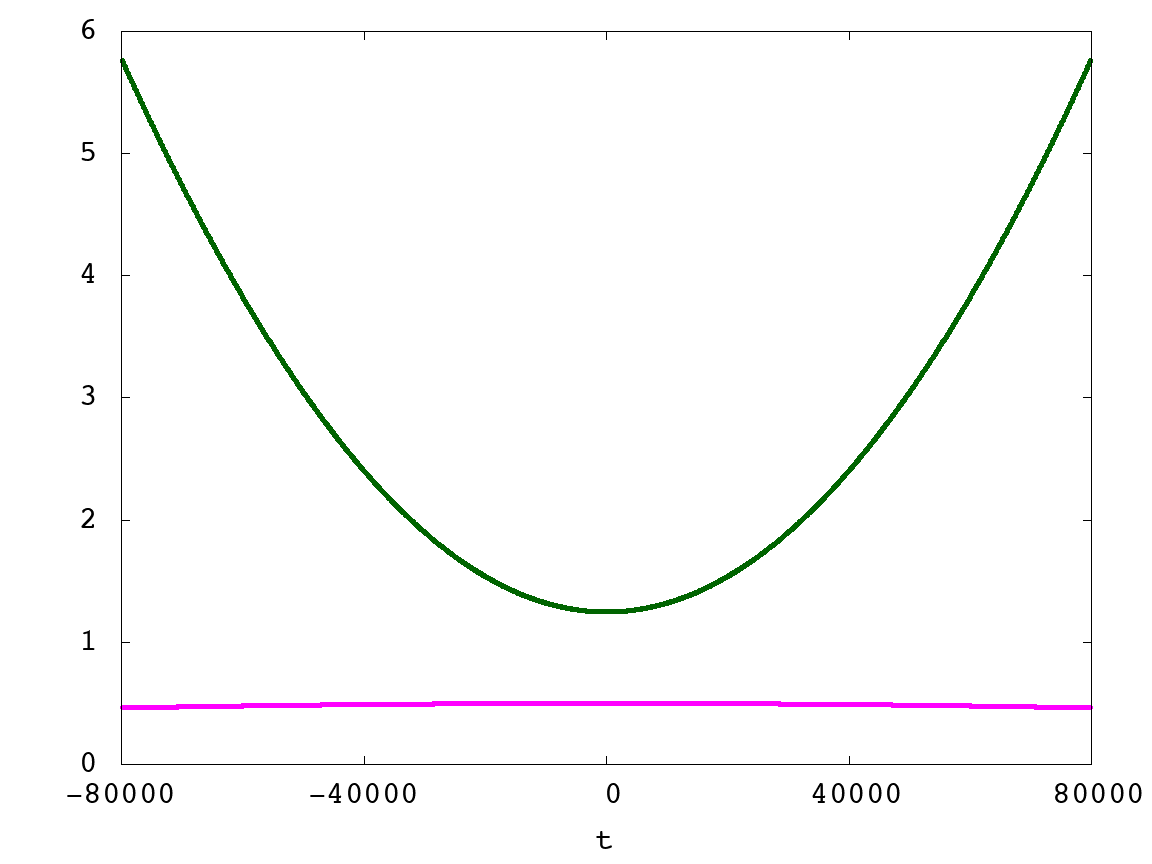}
	\caption{
Graphs of the absolute values of $\widetilde\HH_\CC$ (green) and $\left\|-\frac{1}{2\Lambda^2}\right\|$ (pink)
	along $\widetilde\Gamma_1$ (left) and $\widetilde\Gamma_2$ (right).} \label{htildegamma12} 
\end{figure}

\noindent
The results are plotted in Figures~\ref{energia2} and~\ref{htildegamma12}. Figure~\ref{energia2} shows that, along $\Gamma_2$, $\KK_\CC$ remains ``incapsulated'' at its initial value for much longer a time than along $\Gamma_1$, a somewhat expected fact. However, for the part of the graph represented in such figures,  relations~\eqref{weights} are well maintained along   $\Gamma_1$ and $\Gamma_2$ as well. 
Figure~\ref{htildegamma12} clearly says that, unfortunately, the inequality in~\eqref{weight1} does not hold nor along $\widetilde\Gamma_1$  or along $\widetilde\Gamma_2$.

 \begin{figure}[H]
\centering
\includegraphics[width=7.5cm,height=5cm,draft=false]{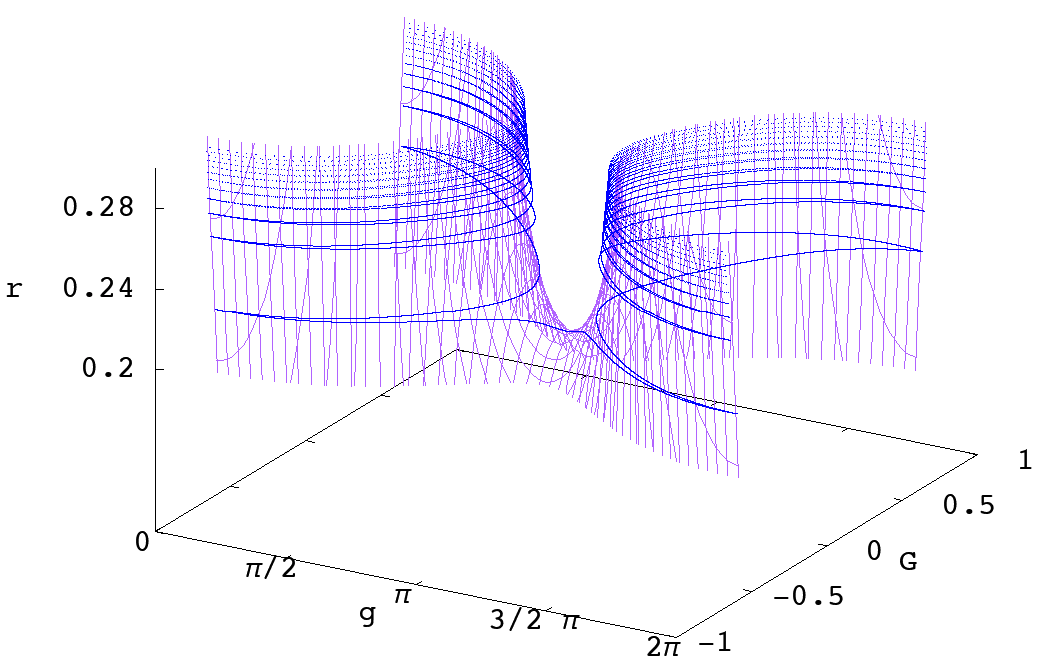}
\includegraphics[width=7.5cm,height=5cm,draft=false]{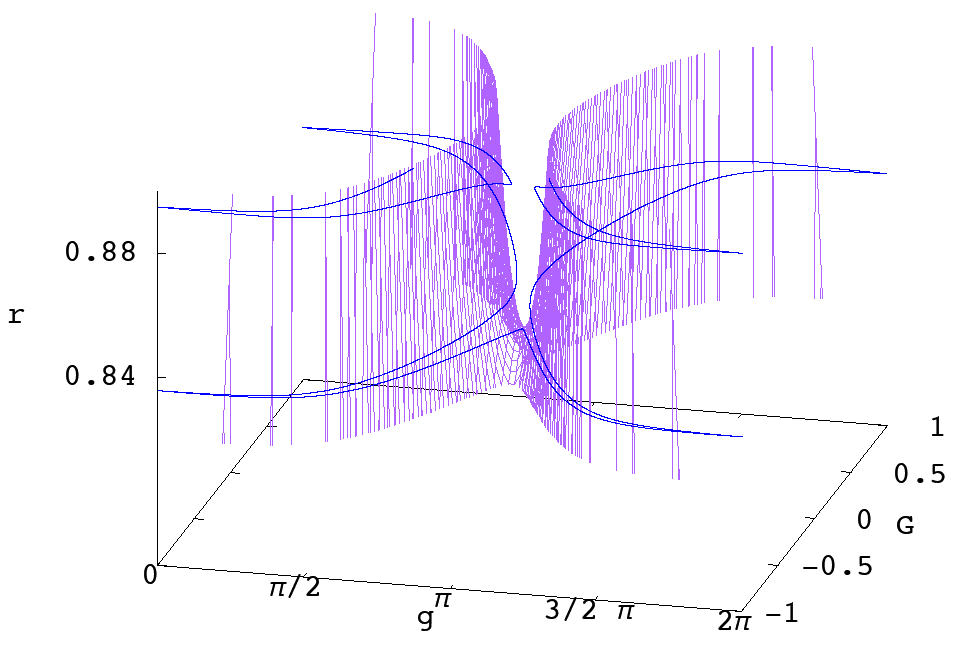}
	\caption{
	Orbits $\Gamma_1$ (left) and $\Gamma_2$ (right) in blue and the respective manifolds ${\cal M}({\cal E}_1)$, ${\cal M}({\cal E}_2)$ in purple.} \label{330_3D} 
\end{figure}

 \subsection*{Example 2}
 We choose\footnote{In this case $\kappa, \mu$ turn out to be  $\kappa=0.06814254$, $\mu=0.02725702$.}
 \begin{equation}
 \label{parameter_1}
\alpha = 50 \, ,\quad \beta = 20 \, , \quad \CC =1.6\, 
 \end{equation}
 and we fix the energy level 
  for the Hamiltonian~\eqref{ham} 
 with the value
  \begin{equation}\label{energy}{\cal H}=-76.887 \, .
 \end{equation}
 On such energy level, we choose three orbits, which we denote as
 $\Gamma_{\rm s}$, $\Gamma_{\rm 0}$ and $\Gamma_{\rm u}$,
 respectively determined by the following initial data
 \begin{eqnarray}\label{stabledatum}
&&\Gamma_{\rm s}:\qquad \left\{
 \begin{array}{lll}\rm R_{\rm s}   =-11.367\\
\rm G_{\rm s} =0.993\\
\rm r_{\rm s} =0.132\\
\rm g_{\rm s} =2.759
\end{array}
\right. \, , \\\nonumber\\\nonumber\\
\label{genericdatum} 
&&\Gamma_{\rm 0}:\qquad\left\{
 \begin{array}{lll}
 \rm R_{\rm 0}   = -9.075\\
\rm G_{\rm 0} = 0.5\\
\rm r_{\rm 0} = 0.132\\
\rm g_{\rm 0} =   \pi
\end{array}
\right.\, ,\\\nonumber\\\nonumber\\
\label{unstabledatum} 
&&\Gamma_{\rm u}:\qquad\left\{
 \begin{array}{lll}\rm R_{\rm u}   =10.331\\
\rm G_{\rm u} =0.718\\
\rm r_{\rm u} =0.132\\
\rm g_{\rm u} =5.699
\end{array}
\right. \, .
\end{eqnarray}
  \begin{figure}[H]
\centering
\includegraphics[width=5cm,height=4cm,draft=false]{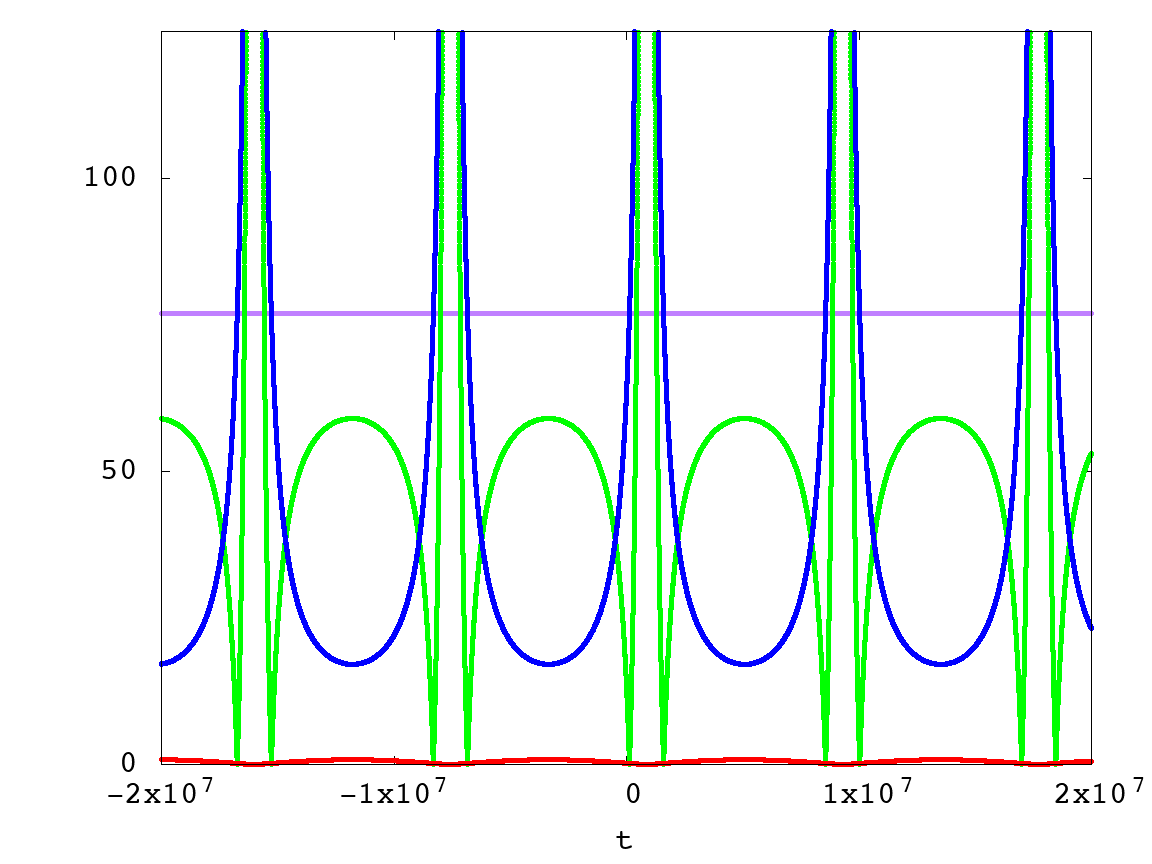}
\includegraphics[width=5cm,height=4cm,draft=false]{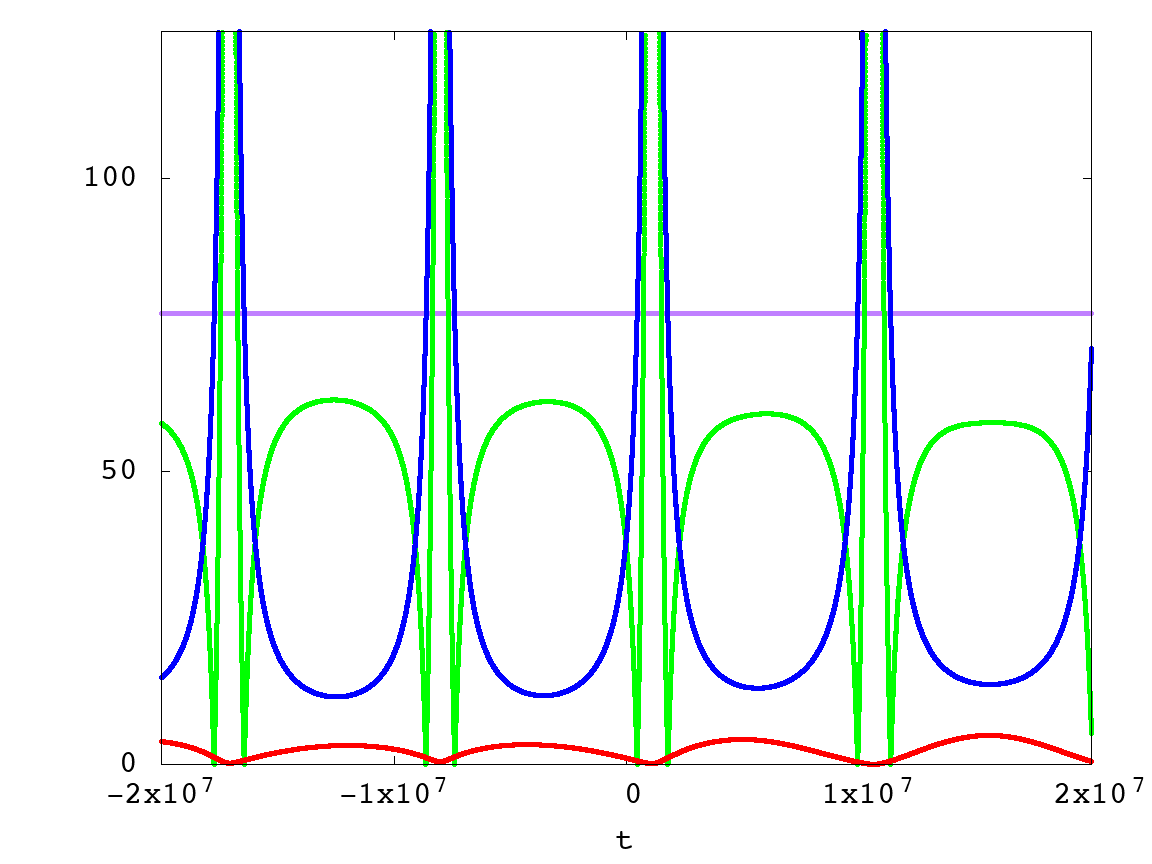}
\includegraphics[width=5cm,height=4cm,draft=false]{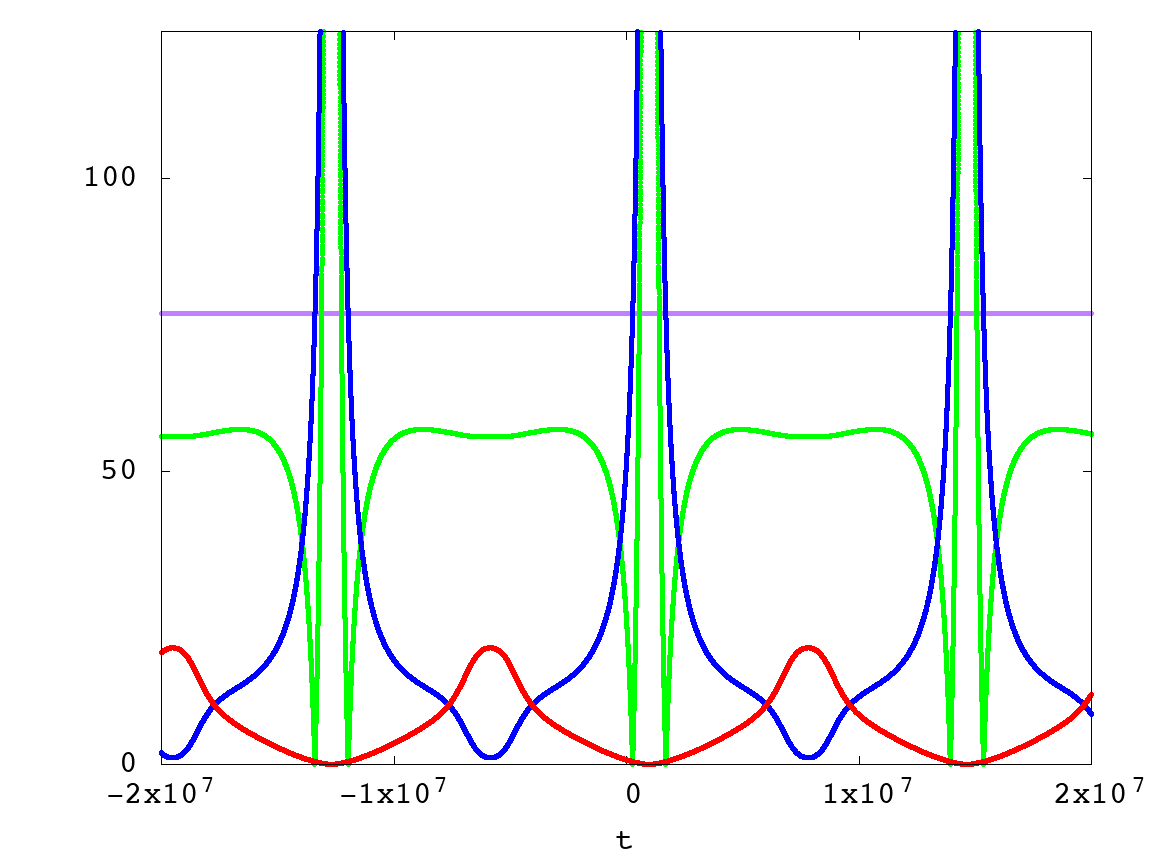}
	\caption{ Graphs of the absolute values of   $\KK_\CC$ (green), $\alpha\UU$ (red) and $\widetilde \KK_\CC$ (blue) along $\Gamma_{\rm s}$ (left), $\Gamma_0$ (centre), $\Gamma_{\rm u}$ (right). The purple line represents the total energy~\eqref{ham}. } \label{energia1} 
\end{figure}

  \begin{figure}[H]
\centering
\includegraphics[width=5cm,height=4cm,draft=false]{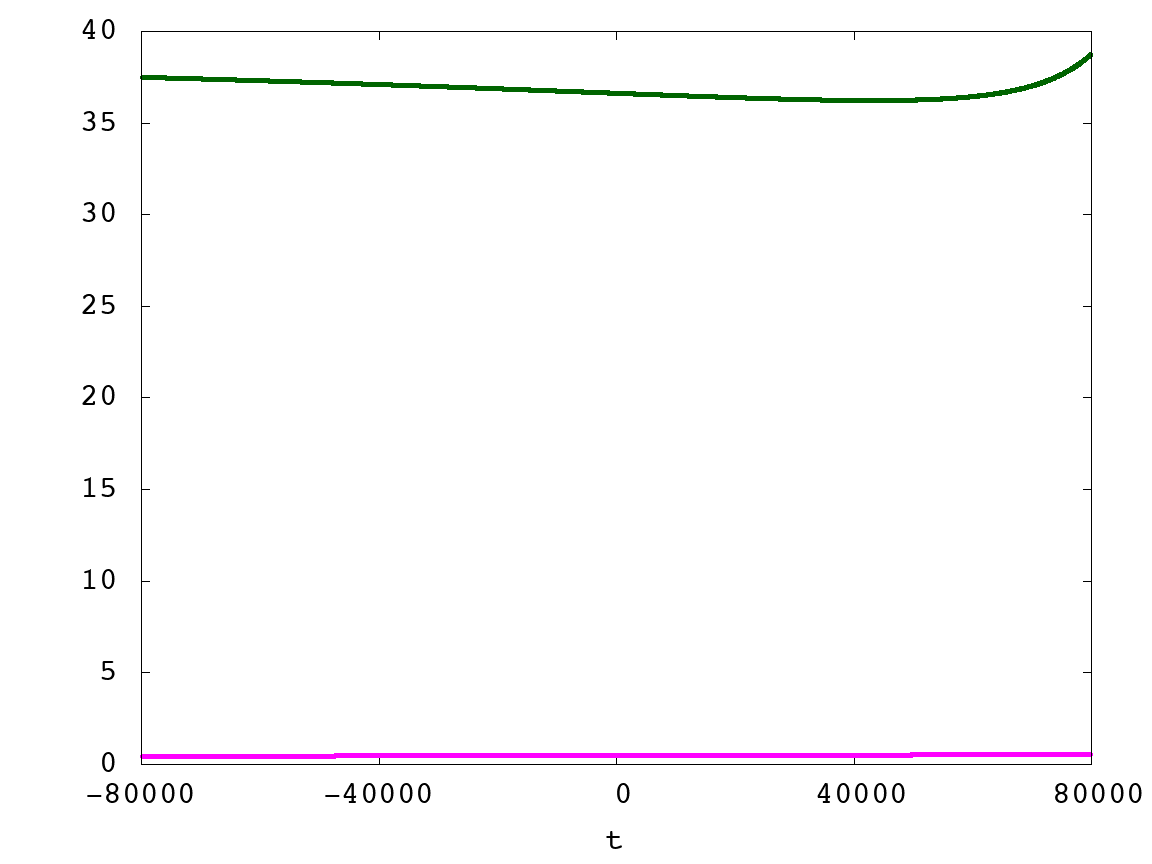}
\includegraphics[width=5cm,height=4cm,draft=false]{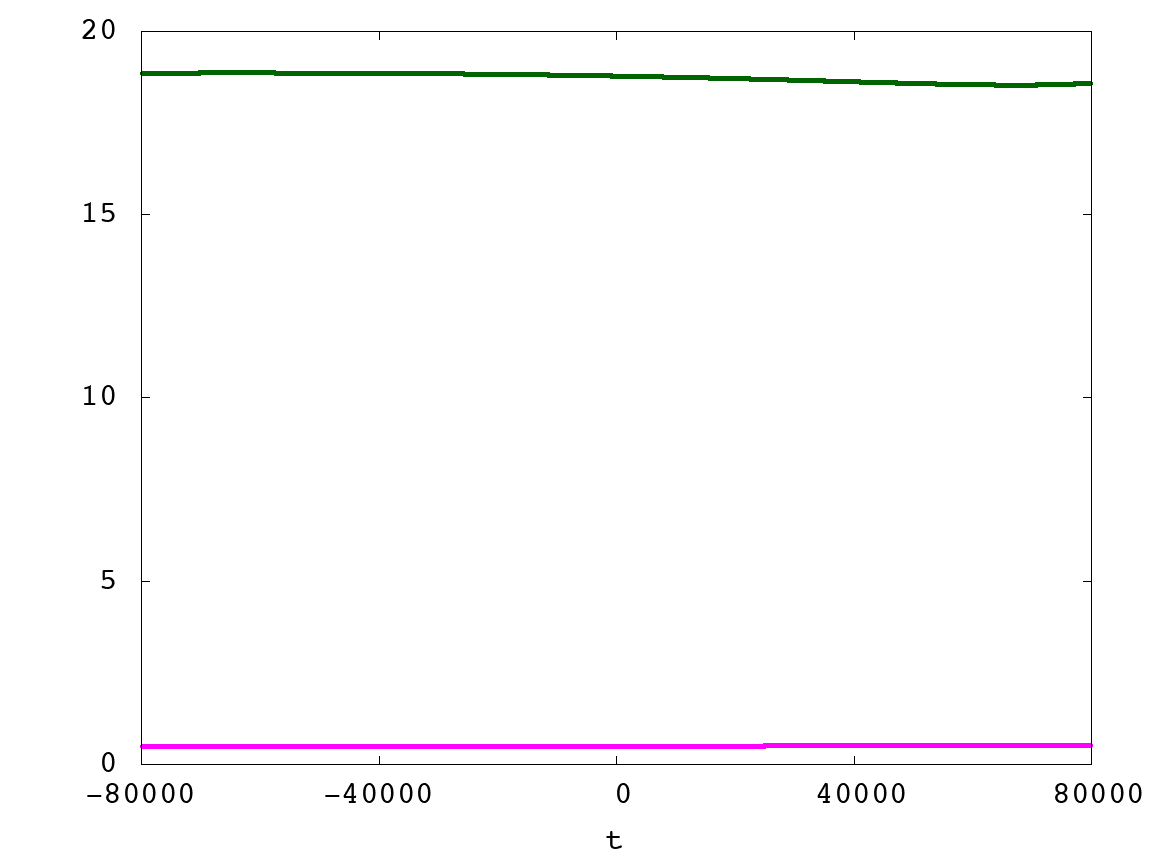}
\includegraphics[width=5cm,height=4cm,draft=false]{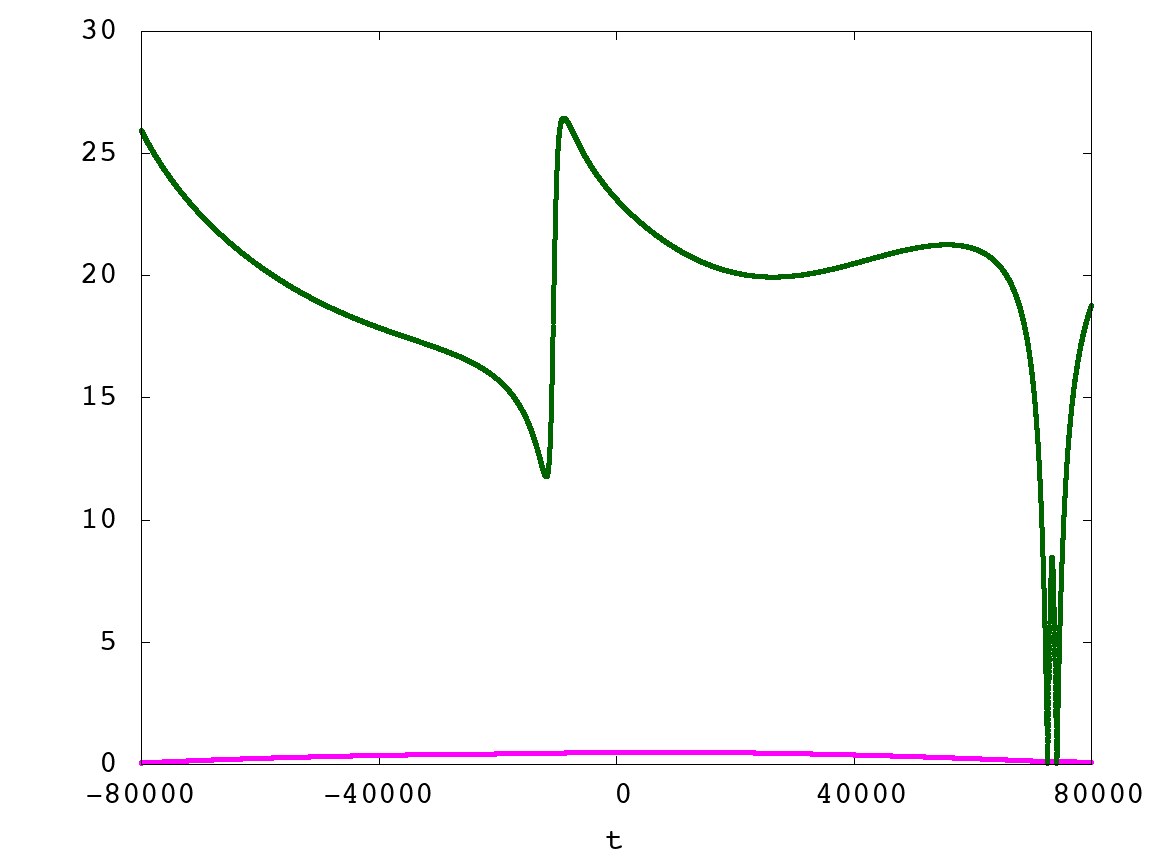}
	\caption{ Graphs of the absolute values of $\widetilde\HH_\CC$ (green) and $\left\|-\frac{1}{2\Lambda^2}\right\|$ (pink)
	along $\widetilde\Gamma_{\rm s}$ (left), $\widetilde\Gamma_0$ (centre), $\widetilde\Gamma_{\rm u}$ (right).
	} \label{htildegammas0u} 
\end{figure}

\noindent
As in the previous example,  inequality~\eqref{weight1} is illustrated on the orbits $\widetilde\Gamma_{\rm s}$, $\widetilde\Gamma_0$ and $\widetilde\Gamma_{\rm u}$ of the whole Hamiltonian $\HH_{\rm 3b, \CC}$ in~\eqref{ovlH}, departing from  initial data obtained completing the respective initial data of $\Gamma_{\rm s}$, $\Gamma_0$ and $\Gamma_{\rm u}$ with $\Lambda_i=1$ and $\ell_i=\pi$.

The results  are plotted in Figures~\ref{energia1} and~\ref{htildegammas0u}. In  Figure~\ref{energia1}, the zones where $\KK_\CC$ and $\widetilde\KK_\CC$ diverge correspond to the coordinate $\rm r$ approaching $0$. As a consequence, we have that~\eqref{weights} is not satisfied on the entire orbits, but only on the portion around the maximum of $|\KK_\CC|$, which corresponds with the zone around the minimum of $|\widetilde\KK_\CC|$. As in the previous example, Figure \ref{htildegammas0u}  shows that the inequality in \eqref{weight1} is not met along any of the orbits $\widetilde\Gamma_{\rm s}$, $\widetilde\Gamma_0$ and $\widetilde\Gamma_{\rm u}$.\\
\vskip.1in
\noindent
As mentioned above, notwithstanding the negative results of Figures \ref{htildegamma12} and \ref{htildegammas0u}, justified  by the considerations  in the introduction, from now on, we focus on the dynamical properties of the Hamiltonian $\overline\HH_\CC$ in \eqref{ham}. Our goal is to check slow variations of the Euler integral in some chosen region of phase space and co--existence of chaotic phenomena. At this respect, we remark that,  even though Figure~\ref{energia2} of Example 1 is encouraging, proving existence of chaos closely to ${\cal M}_0$ along this way seems really hard. The difficulty is that, even though the initial point of $\Gamma_i$ has been chosen {\it precisely} on the saddle of ${\cal M}({\cal E}_i)$, and, for a long time, the orbits maintains to be very close to ${\cal M}({\cal E}_i)$, however, the coordinate $\rm r$ increases such in a way to  leave the region $0<\rm r<1$ (hence, the region of the saddle) in a short time; see Figure~\ref{330_3D}. For this reason, in the rest of the paper we shall be focused on the orbits $\Gamma_{\rm s}$ and $\Gamma_{\rm u}$ in Example 2, where the motion of $\rm r$ is sufficiently slow.

\section{The orbits $\Gamma_{\rm s}$ and $\Gamma_{\rm u}$}\label{eul_integral}
Let us consider the Hamiltonian~\eqref{ham}, with $\alpha$, $\beta$ and $\CC$ as in~\eqref{parameter_1}. We fix the value ${\cal H}$ of the energy  as in~\eqref{energy} and we reduce the coordinate $\rm R$ via \begin{eqnarray}\label{energy reduction}\rm R=\pm \sqrt{2\left(\rm {\cal H}-\alpha\UU(r, G, g)-\frac{\rm(C-G)^2}{2{\rm r}^2}+\frac{\beta}{\rm r}\right)}\,,\end{eqnarray}
with the sign being chosen by continuity. We look at the motion of the triplet $\rm (r, G, g)$ in a 3--dimensional space.

\begin{figure}[H]
 \centering
 \includegraphics[width=7cm,height=5cm,draft=false]{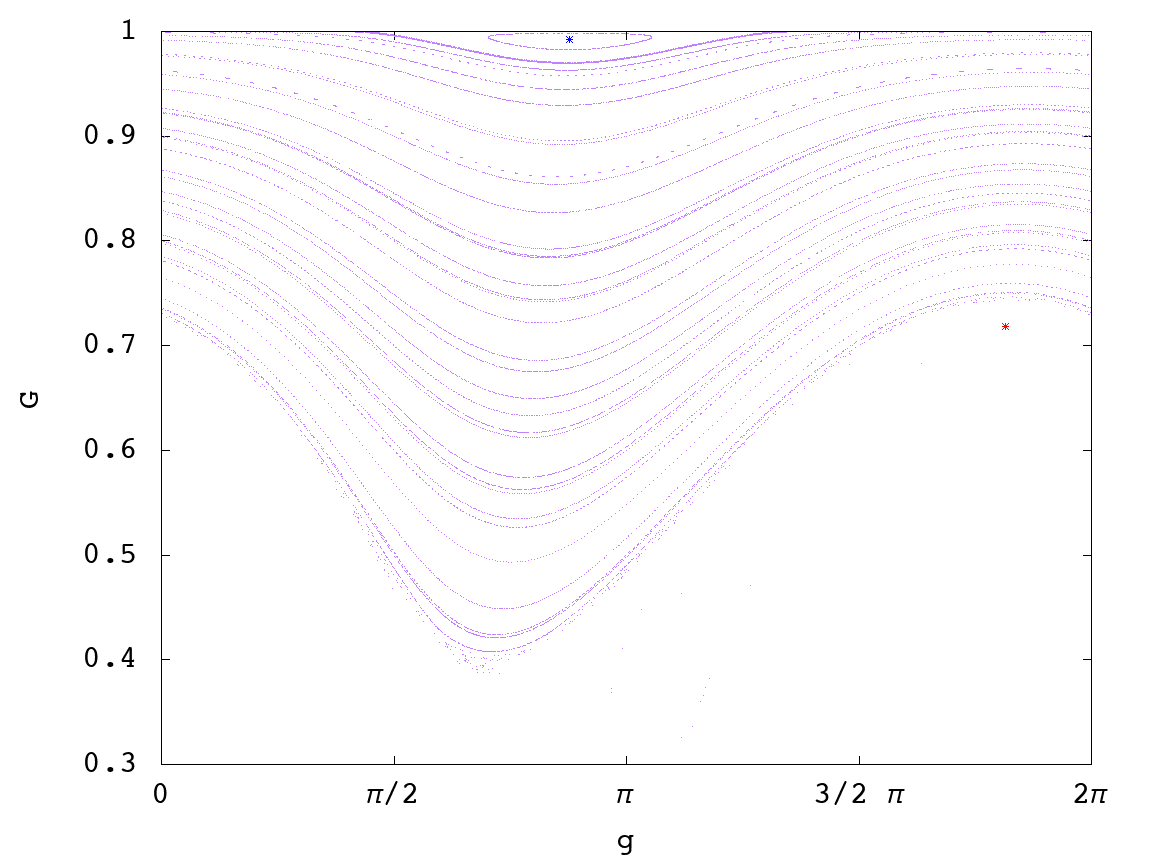}
  \includegraphics[width=8cm,height=5cm,draft=false]{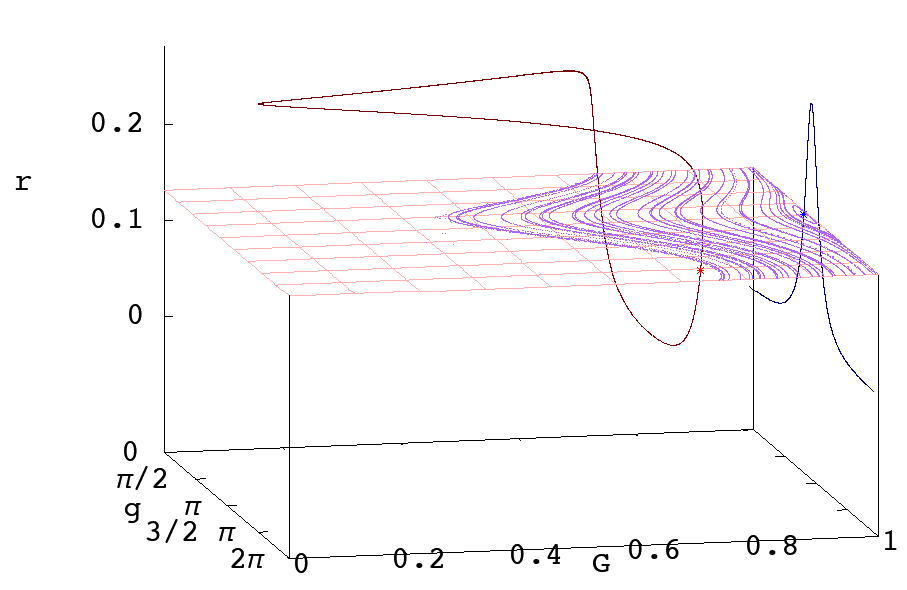}
 \caption{Left: the Poincar\'e map~\eqref{poinc} with ${\cal H}$ in~\eqref{energy} and $\Pi_{\rm s}$ orthogonal to $\Gamma_{\rm s}$ at $(\rm r_{\rm s}, G_{\rm s}, g_{\rm s})$. Right: spatial visualisation with $\Gamma_{\rm s}$, $\Gamma_{\rm u}$ in blue and in red, respectively.}  \label{Poi_section}
 \end{figure}

\noindent
 We empirically find a periodic orbit of $\HH$ in~\eqref{ham} in correspondence of the initial datum~\eqref{stabledatum}.
We denote as $\Gamma_{\rm s}$ the projection of such  orbit in the space $\rm (r, G, g)$.
We choose $\Pi_{\rm s}$ as the plane orthogonal to $\Gamma_{\rm s}$
at $\rm (r_{\rm s}, G_{\rm s}, g_{\rm s})$.  
We  construct a 2--dimensional map
\begin{equation}\label{poinc}{\cal P}_{\!\textrm{\tiny$\cal H$,$\Pi_{\rm s}$}}:\qquad ({\rm G}, {\rm g})\to ({\rm G}', {\rm g}')
\end{equation}
where $({\rm G}', {\rm g}')$ is the first return value on $\Pi_{\rm s}$.
By construction, $\rm (G_{\rm s}, g_{\rm s})$ is a fixed point of ${\cal P}_{\!\textrm{\tiny$\cal H$,$\Pi_{\rm s}$}}$.
The images of the map ${\cal P}_{\!\textrm{\tiny$\cal H$,$\Pi_{\rm s}$}}$ with $\Pi_{\rm s}$ as said are depicted in Figure~\ref{Poi_section}, left. A Newton Algorithm is used to find other fixed points, besides $\rm (G_{\rm s}, g_{\rm s})$.  Another point is actually found $\rm (G_{\rm u}, g_{\rm u})$, which (using the equation of $\Pi_{\rm s}$ and of the energy reduction~\eqref{energy reduction}) unfolds to the quadruplet~\eqref{unstabledatum}.
Amazingly, we did not find other fixed points of ${\cal P}_{\!\textrm{\tiny$\cal H$,$\Pi_{\rm s}$}}$:
\begin{numevid}\label{fixed points poinc} If ${\cal H}$ is as in~\eqref{energy} and $\Pi_{\rm s}$ is orthogonal to $\Gamma_{\rm s}$ at $\rm (r_s, G_s, g_s)$, then  $\rm (G_s, g_s)$, $\rm (G_u, g_u)$ are the only fixed points of ${\cal P}_{\!\textrm{\tiny$\cal H$,$\Pi_{\rm s}$}}$.
$\quad\square$\end{numevid}
In Table~\ref{fixed_point}  we report the value of the two fixed points. This is to be compared with the situation studied in~\cite{diruzzaDP20}, where several hyperbolic points in a chaotic region were numerically detected.
\begin{table}[H]
\centering
\begin{tabular}{|c|c|}
   \hline
$\rm (G_{\rm s}, g_{\rm s})$ & $(0.992515 , 0.878179 \,\pi )$ \\
\hline
$\rm (G_{\rm u}, g_{\rm u})$&  $( 0.717909, 1.81405 \, \pi )$  \\
\hline
\end{tabular}.    
\caption{Values of fixed points.}\label{fixed_point}. 
\end{table}
\noindent
The computation of the eigenvalues of the linear part of ${\cal P}_{\!\textrm{\tiny$\cal H$,$\Pi_{\rm s}$}}$ at $\rm (G_{\rm s}, g_{\rm s})$ and $\rm (G_{\rm u}, g_{\rm u})$
 assigns to $\rm (G_{\rm s}, g_{\rm s})$ the character of elliptic fixed point (for having complex eigenvalues), and to $\rm (G_{\rm u}, g_{\rm u})$ the 
character of hyperbolic fixed point (for having real eigenvalues, one inside, one outside the unit circle) for ${\cal P}_{\!\textrm{\tiny$\cal H$,$\Pi_{\rm s}$}}$ (see Table~\ref{eigenvalues}).
\begin{table}[H]
\centering
\begin{tabular}{|c|c|c|}
\hline
   & eigenvalue 1 & eigenvalue 2\\
   \hline
$\rm (G_{\rm s}, g_{\rm s})$ & $0.99568665  + i \, 0.99278032$ & $0.99568665  - i \, 0.99278032$ \\
\hline
$\rm (G_{\rm u}, g_{\rm u})$&  $- 0.051632609$ & $- 19.366447$ \\
\hline
\end{tabular}. 
\caption{Eigenvalues of fixed points $\rm (G_{\rm s}, g_{\rm s})$ and $\rm (G_{\rm u}, g_{\rm u}).$} \label{eigenvalues}
\end{table}
\noindent Therefore, we shall refer to the periodic orbits $\Gamma_{\rm s}$, $\Gamma_{\rm u}$ through  $\rm (r_{\rm s}, G_{\rm s}, g_{\rm s})$, $\rm (r_{\rm u}, G_{\rm u}, g_{\rm u})$ in the space $\rm (r, G, g)$ as ``elliptic'', ``hyperbolic'' periodic orbit, respectively. Such orbits are depicted in
 Figure~\ref{Poi_section}, right, where also the plane $\Pi_{\rm s}$ is visualised.

 \noindent
  In the next, in order to study the validity of Conjecture~\ref{conj}, i), we study how the Euler integral changes along such orbits. We shall see that $\Gamma_{\rm s}$ is immersed in a region of phase space close to ${\cal M}_1$ at all times, while $\Gamma_{\rm u}$ spends much time close to ${\cal M}_0$.

 \subsection*{Spread of E about $\bm\Gamma_{\mathbf s}$}
\begin{figure}[H]
 \centering
 \includegraphics[width=3.5cm,height=3cm,draft=false]{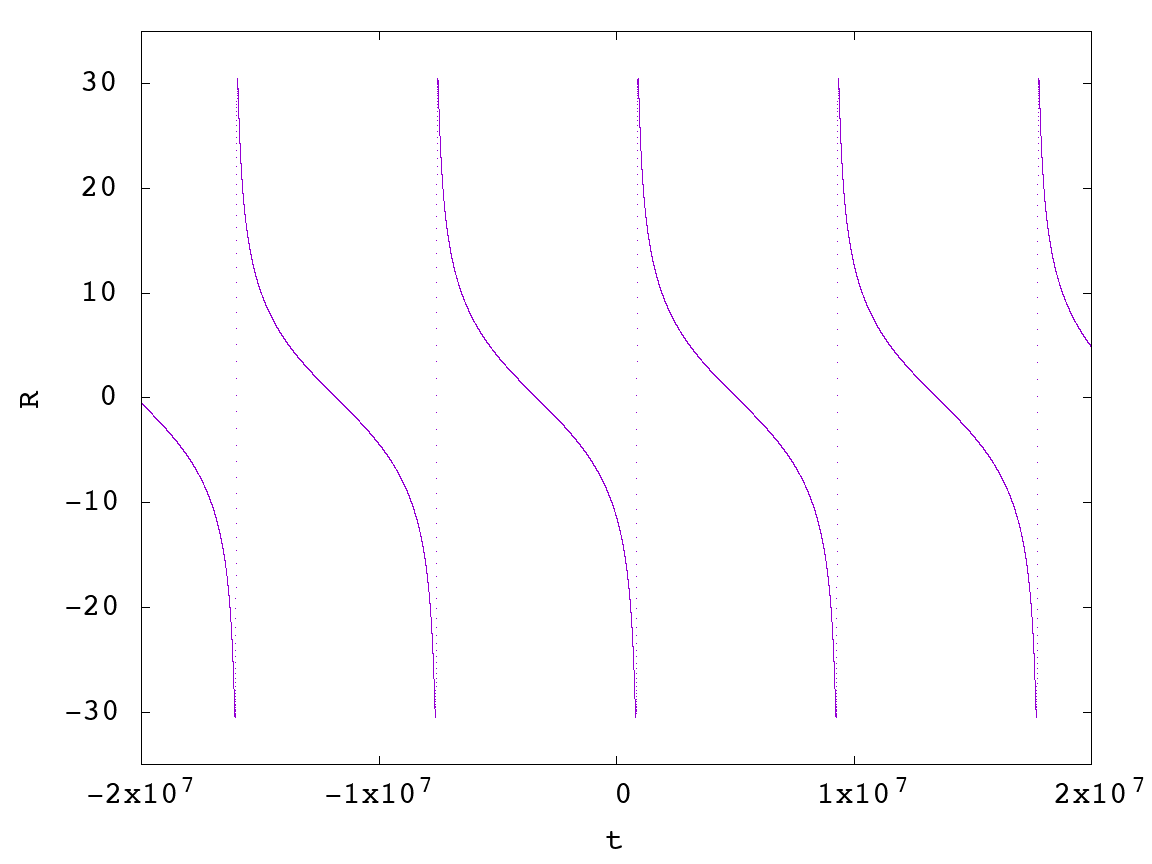}
 \includegraphics[width=3.5cm,height=3cm,draft=false]{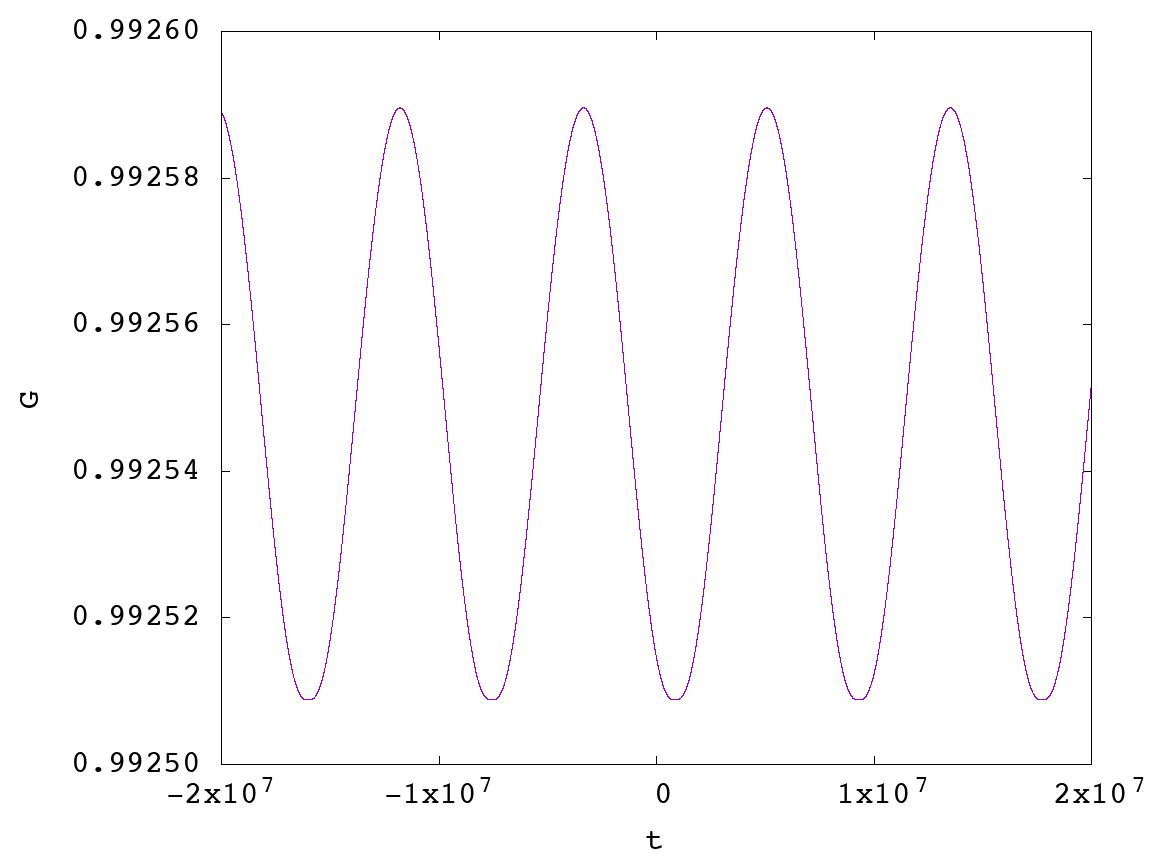}
 \includegraphics[width=3.5cm,height=3cm,draft=false]{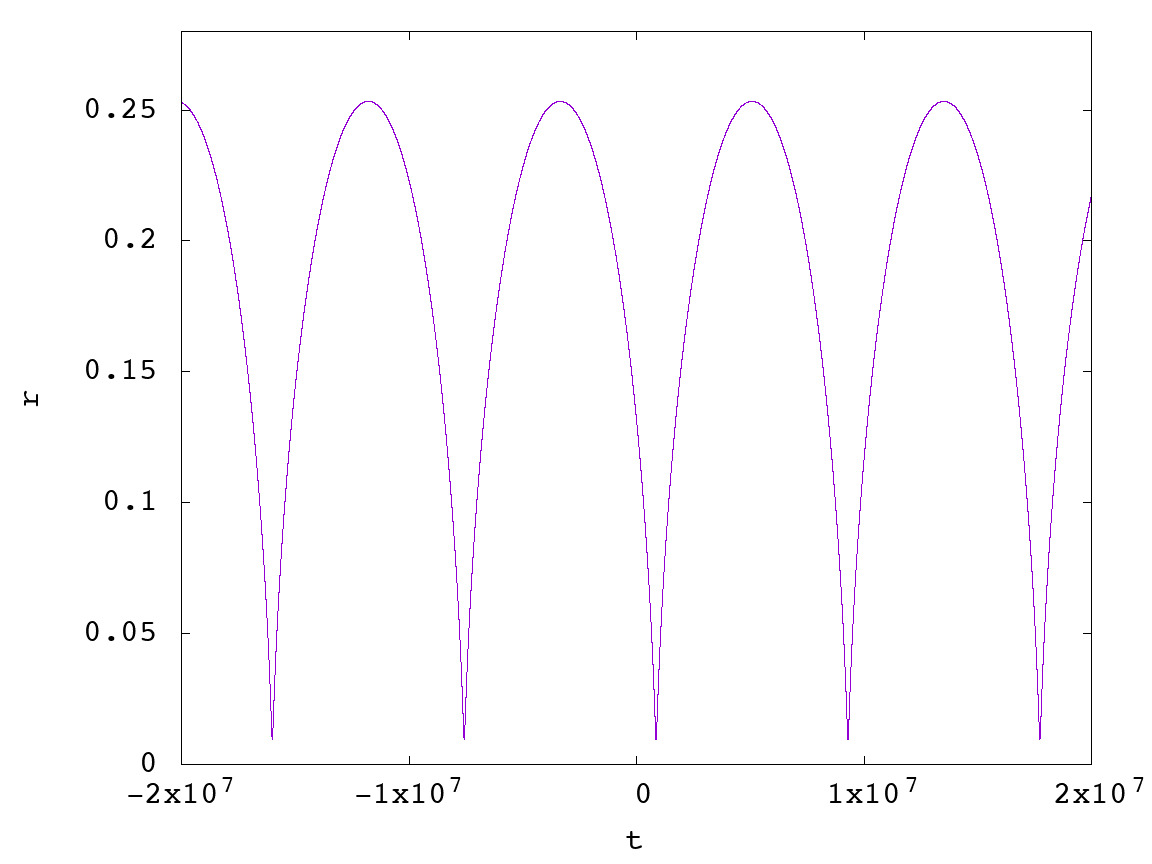}
 \includegraphics[width=3.5cm,height=3cm,draft=false]{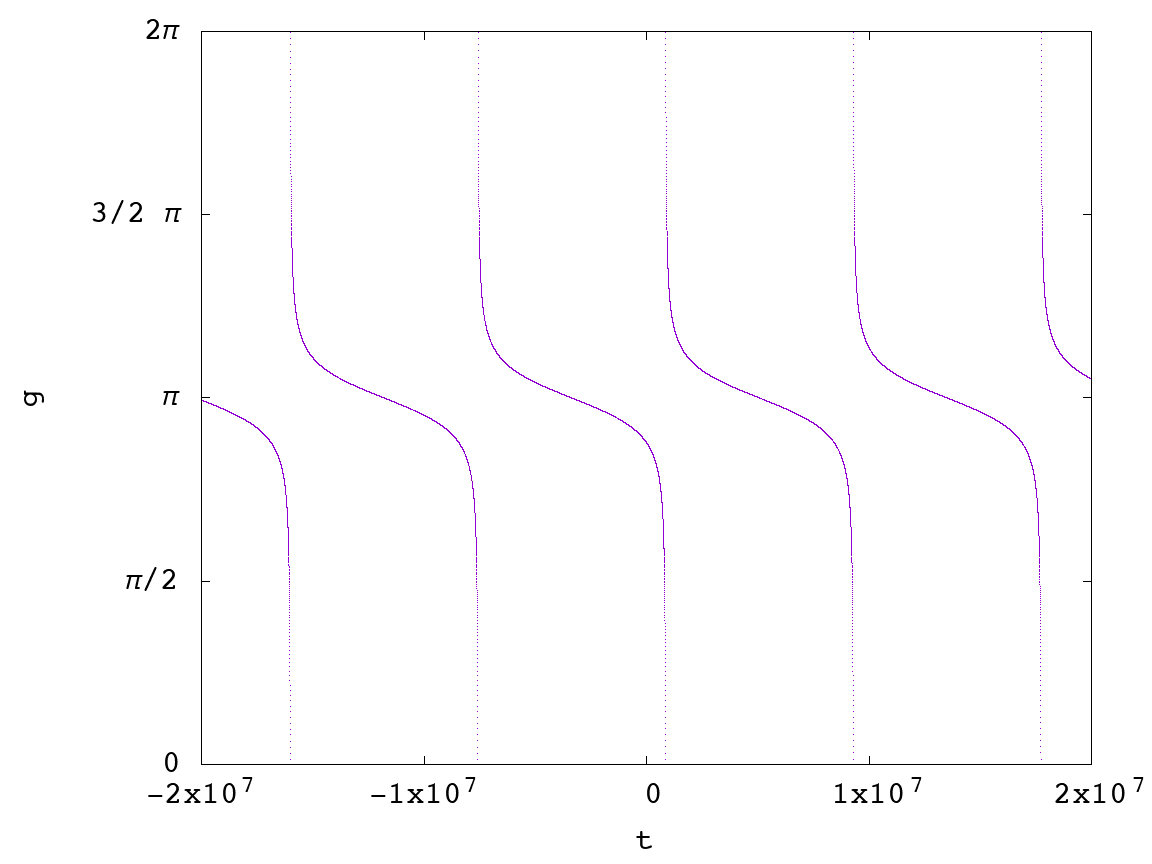} \\
 \includegraphics[width=3.5cm,height=3cm,draft=false]{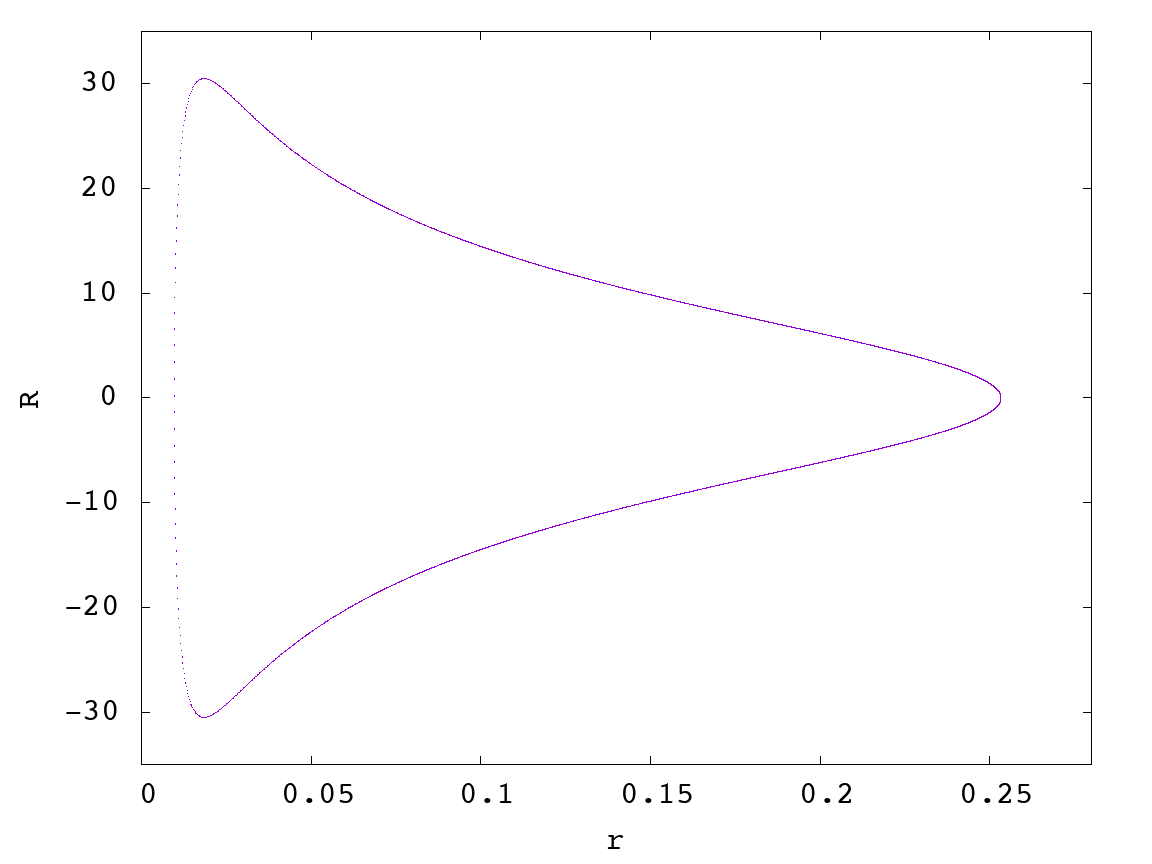}
 \includegraphics[width=3.5cm,height=3cm,draft=false]{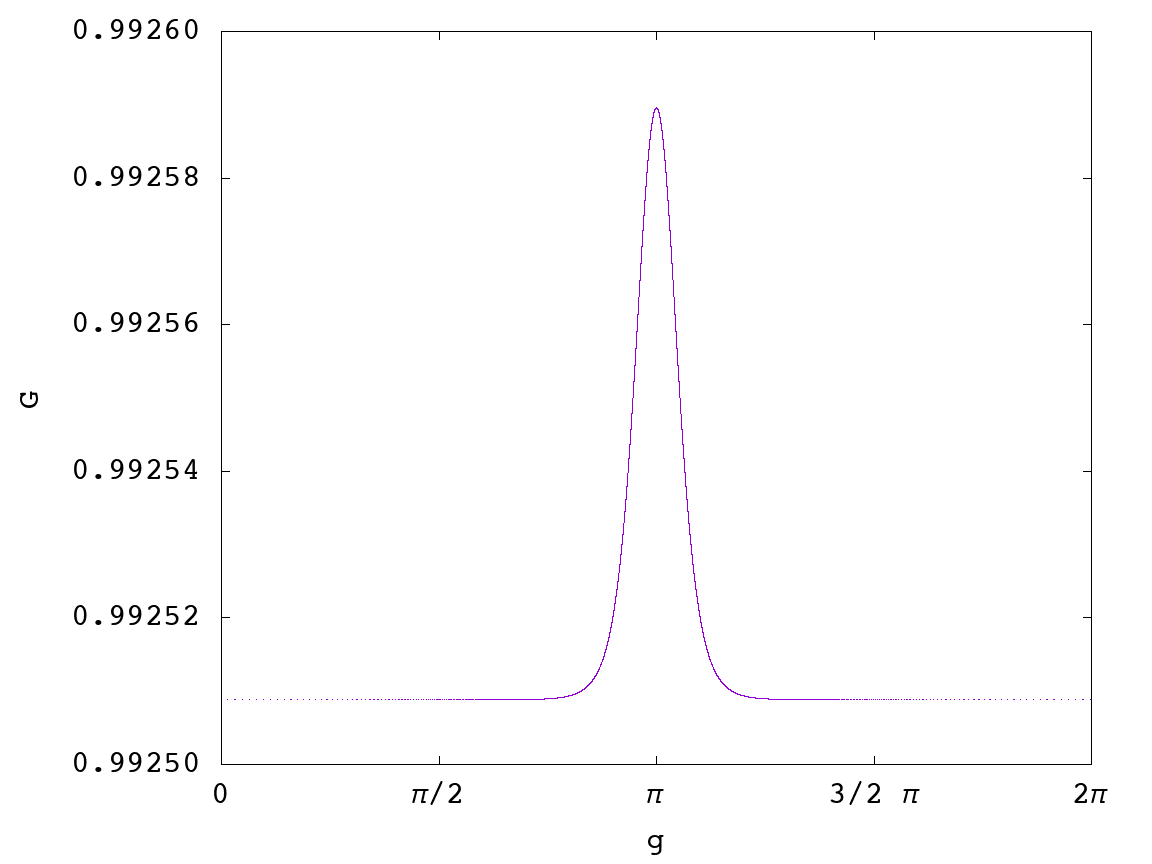}
 \includegraphics[width=3.5cm,height=3cm,draft=false]{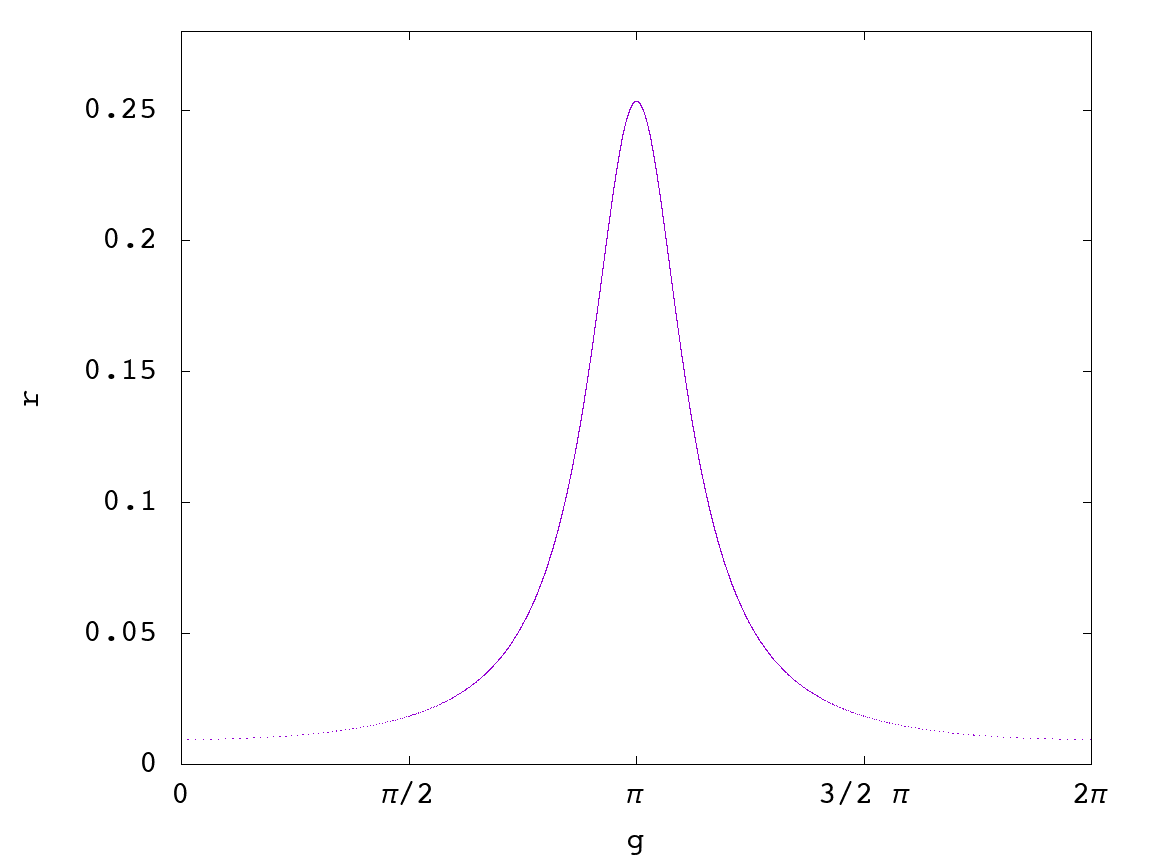}
 \includegraphics[width=3.5cm,height=3cm,draft=false]{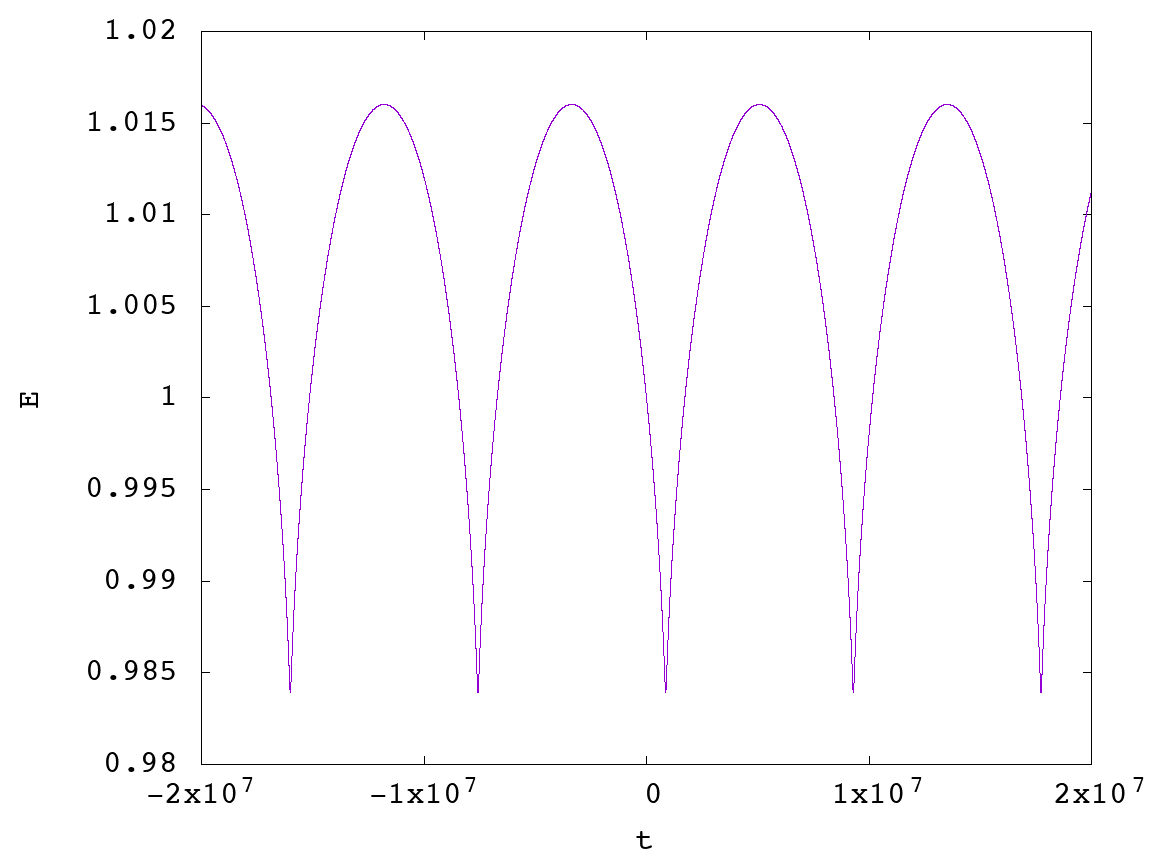}
 \caption{ Details on the orbit $\Gamma_{\rm s}$.}  \label{334_orbit}
 \end{figure}

  In the top panel of Figure~\ref{334_orbit},  the time variations of  $\rm R$, $\rm G$, $\rm r$, $\rm g$ along $\Gamma_{\rm s}$  are represented. In particular, 
 $\rm g$ spans $[0, 2 \pi]$ while 
 $\rm G$ 
 has a very short range of variation (the relative variation changes periodically by a factor of order of $10^{-4}$, so it turns to be quasi--constant). In the bottom panel, we have represented the motion in the planes $\rm(R,r)$, $\rm(G,g)$, $\rm(r,g)$ and the variation of $\rm E$ along the orbit. The latter plot shows that $\rm E$ varies a little, taking values very close to $1$. This means that the orbit is in a zone of phase space very close to ${\cal S}_1(\rm r)$ (see~\eqref{S0S1}) with $\rm r$ taking the mentioned values. 

\begin{figure}[H]
\centering
\includegraphics[width=8.5cm,height=6cm,draft=false]{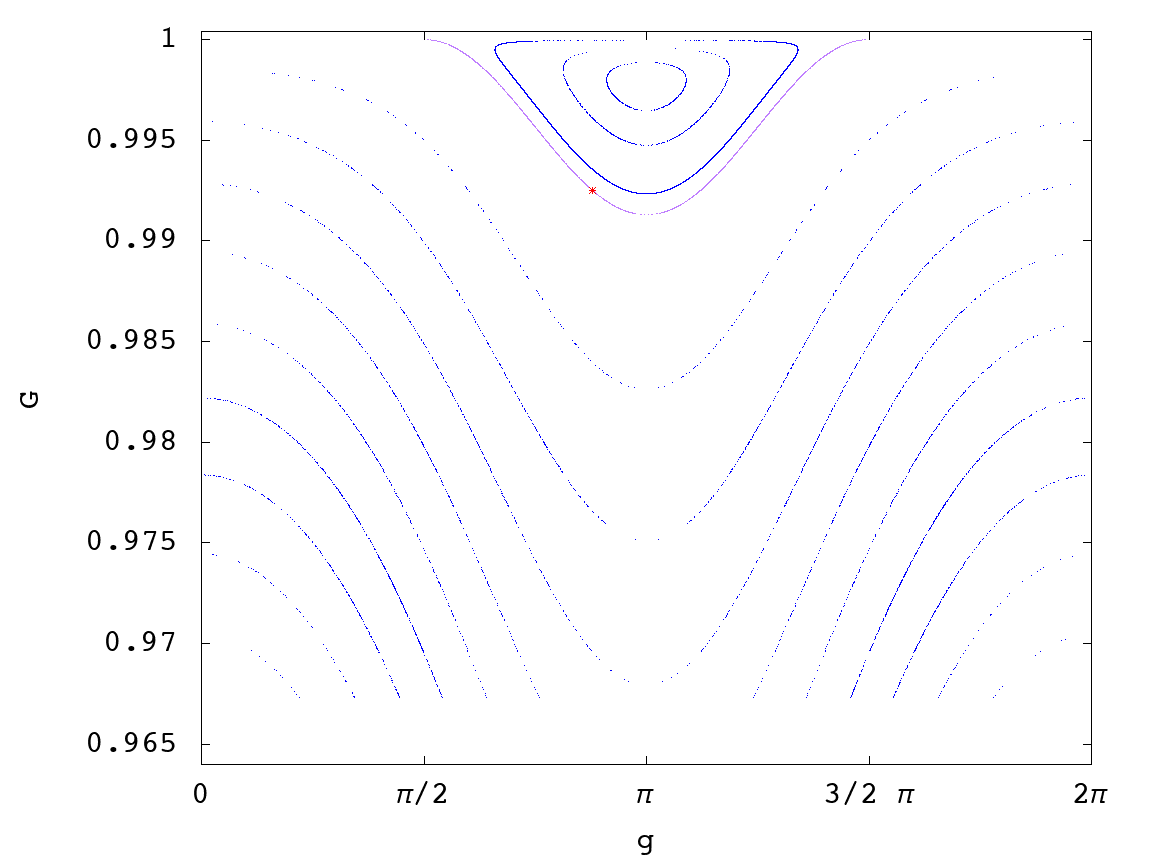} \\
\includegraphics[width=7.5cm,height=6cm,draft=false]{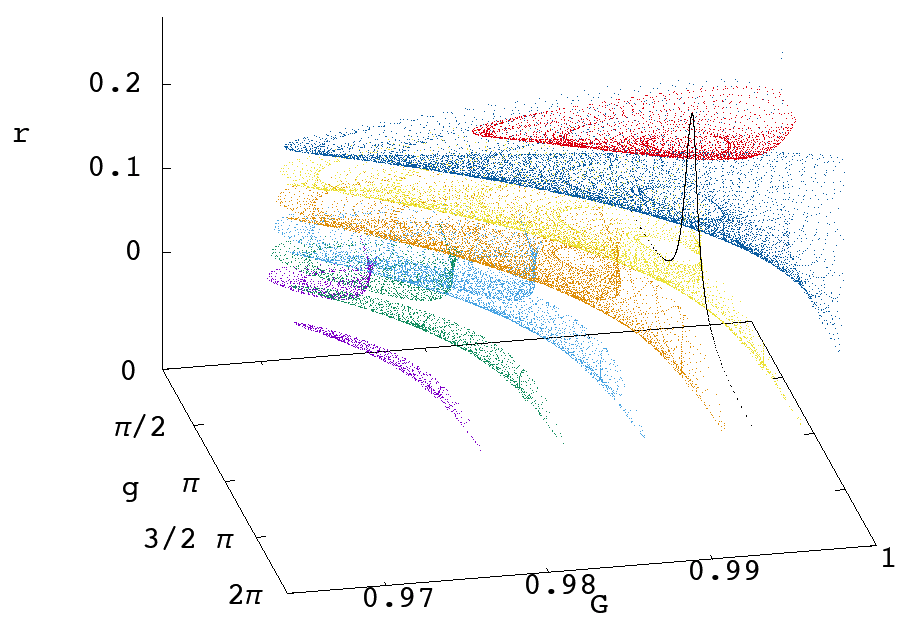}
\includegraphics[width=7.5cm,height=6cm,draft=false]{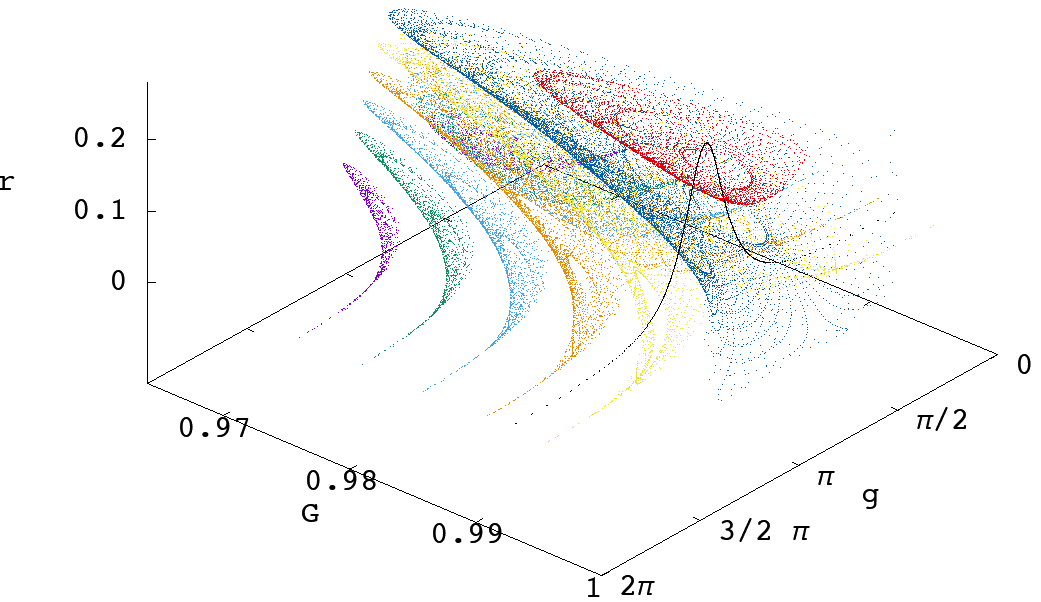}
\caption{
Curves (top) ${\cal S}({\rm r}_0, \theta)$ and manifolds (bottom) ${\cal M}(\theta)
$ intersected by the $\overline\HH_\CC$ evolution, 
	with initial data
	$\rm (r_0, G_0, g_0)\in{\cal M}_1$, and the curve $\Gamma_{\rm s}$.
 }
 \label{fase_rfix_cost} 
\end{figure}

\noindent 
In order to inspect the variations of $\rm E$ in a neighbourhood of $\Gamma_{\rm s}$, we proceed as follows. We choose a grid of initial conditions $\rm (R_0,  G_0, r_0, g_0)$ on the same energy level~\eqref{energy}, and verifying \begin{eqnarray}\label{incond1}
\rm (r_0, G_0, g_0)\in{\cal M}_1
	\end{eqnarray}
	(with ${\cal M}_1$ as in~\eqref{M0M1})
and let the system evolve under $\overline\HH_\CC$. Then in the plane $\rm (G, g)$ we mark a point whenever $\rm (G, g)\in{\cal S}({\rm r}_0, \theta)$, with ${\cal S}({\rm r}, {\cal E})$ as in~\eqref{level curves} for the initial values in the grid. We find that 
\begin{numevid}\label{numevid_gammas} The only non void level manifolds ${\cal M}(\theta)$ intersected by
orbits of $\overline\HH_\CC$ with initial data in ${\cal M}_1$
 are those with \begin{equation}\label{theta}\theta \in [0.91,1.01]\,.\end{equation}
$\quad\square$
\end{numevid}
We report the results in Figure~\ref{fase_rfix_cost} (top), with the purple curve corresponding to $\theta=1$ and the blue curves to different values of $\theta$ in~\eqref{theta}. The red point in the figure represents $\rm (G_s, g_s)$. 

\noindent
We also provide a spatial visualisation,  reporting in Figure~\ref{fase_rfix_cost}  (bottom)  manifolds ${\cal M}(\theta)$ in~\eqref{eul_ell_surf} which are intersected by the time evolution of  initial data in ${\cal M}_1$, with $\theta$ as in~\eqref{theta}.

\noindent
 \subsection*{Spread of E about $\bm\Gamma_{\mathbf u}$}
 \begin{figure}[H]
 \centering
 \includegraphics[width=3.5cm,height=3cm,draft=false]{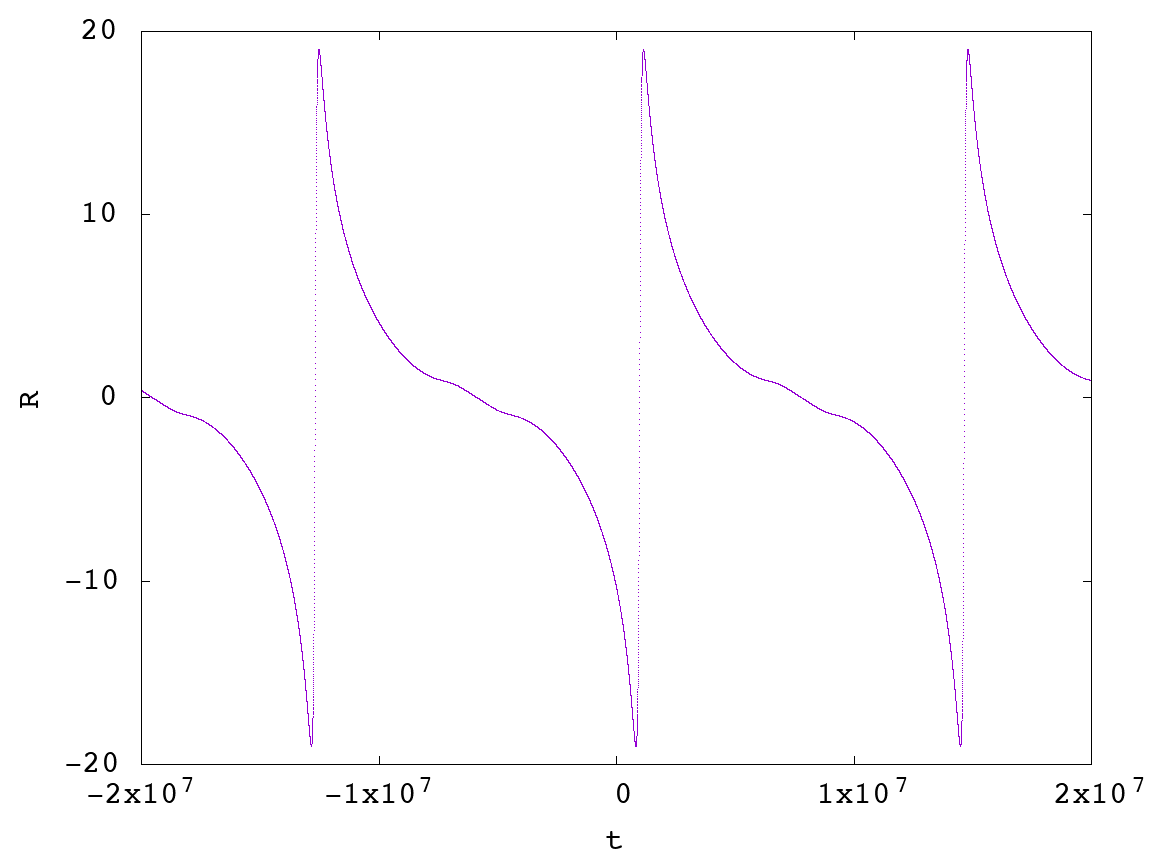}
 \includegraphics[width=3.5cm,height=3cm,draft=false]{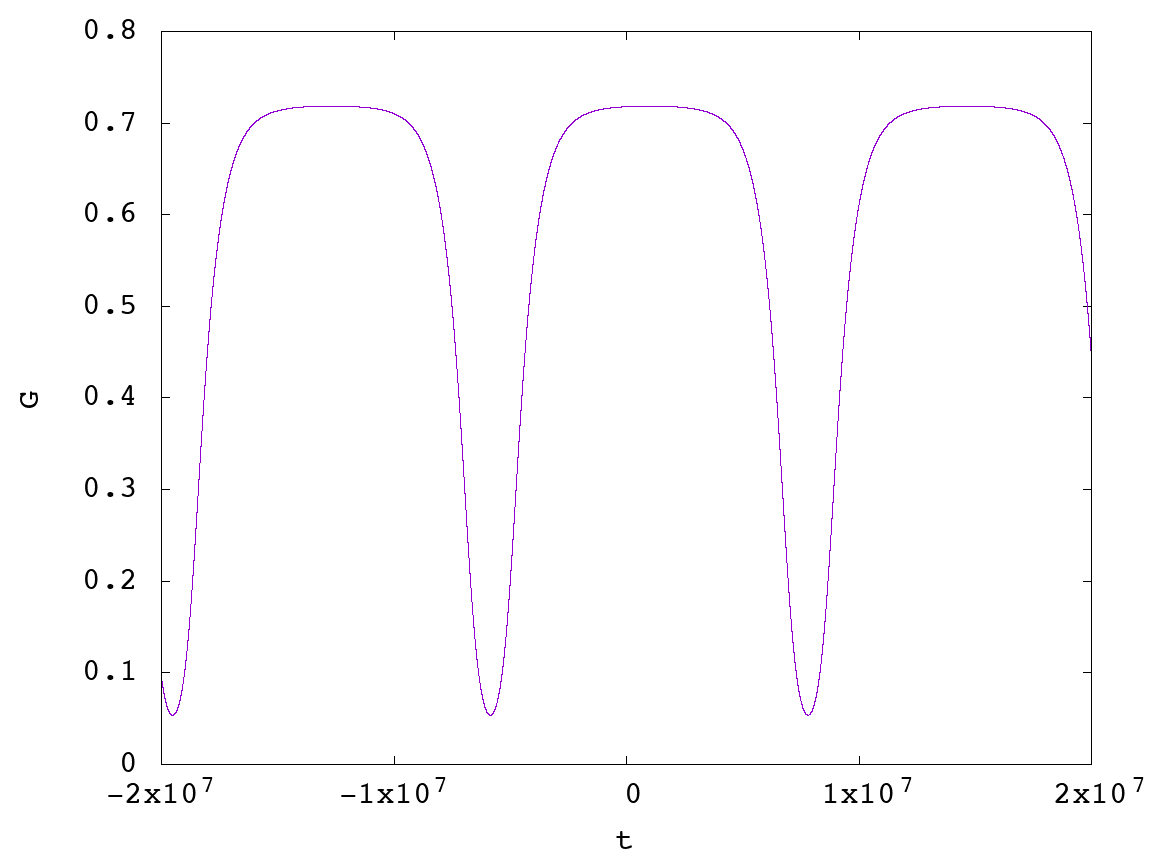}
 \includegraphics[width=3.5cm,height=3cm,draft=false]{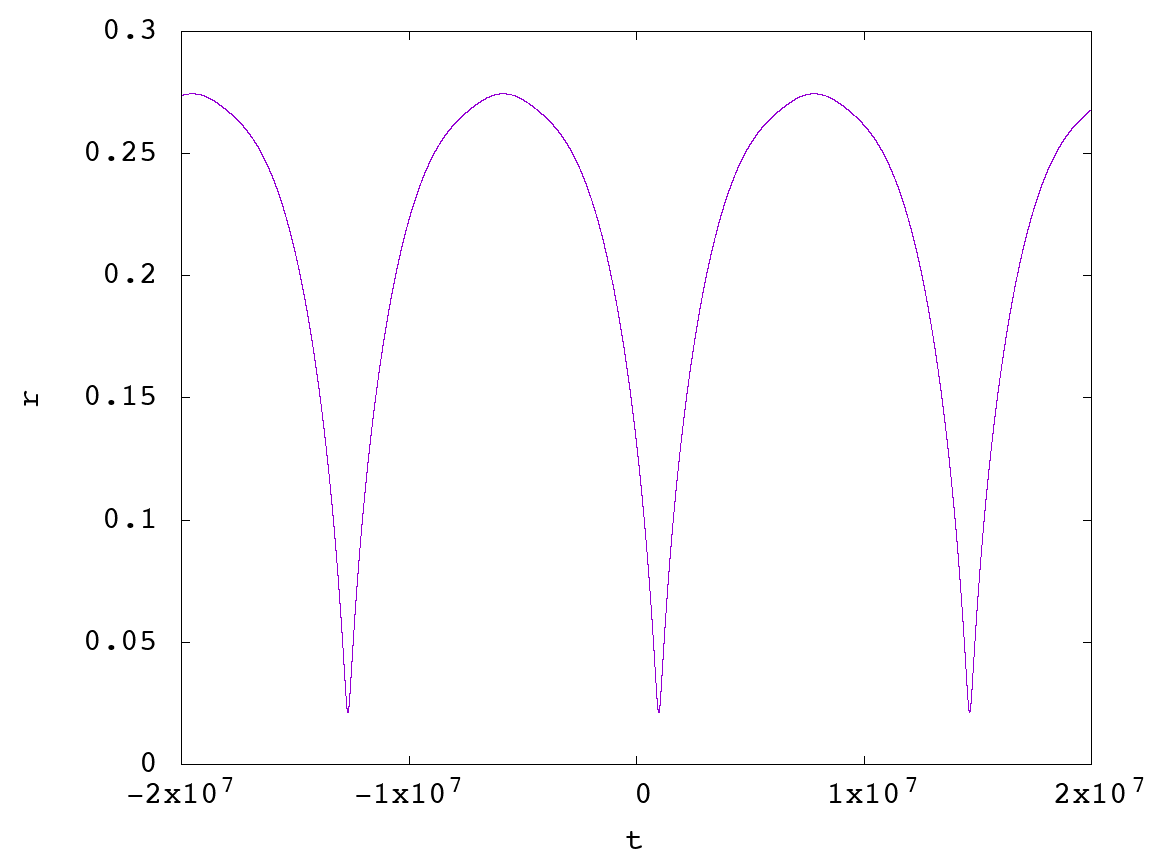}
 \includegraphics[width=3.5cm,height=3cm,draft=false]{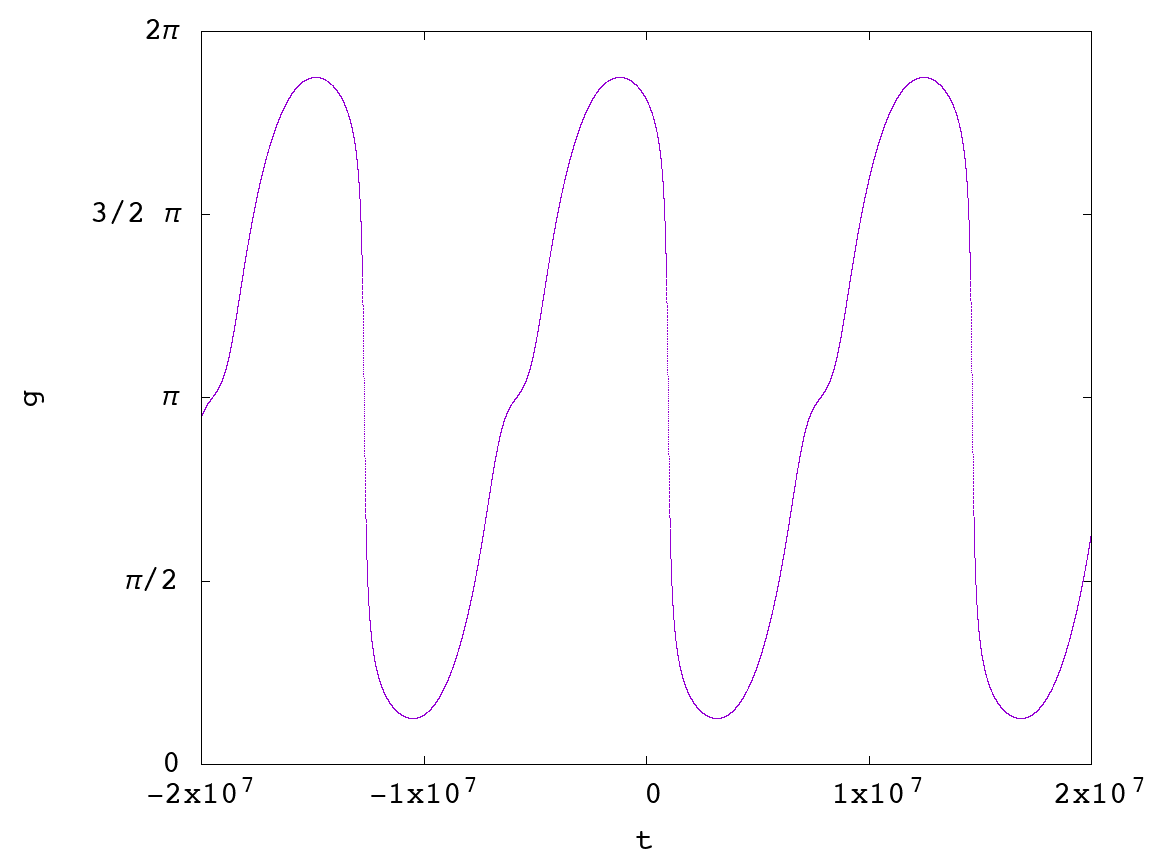} \\
 \includegraphics[width=3.5cm,height=3cm,draft=false]{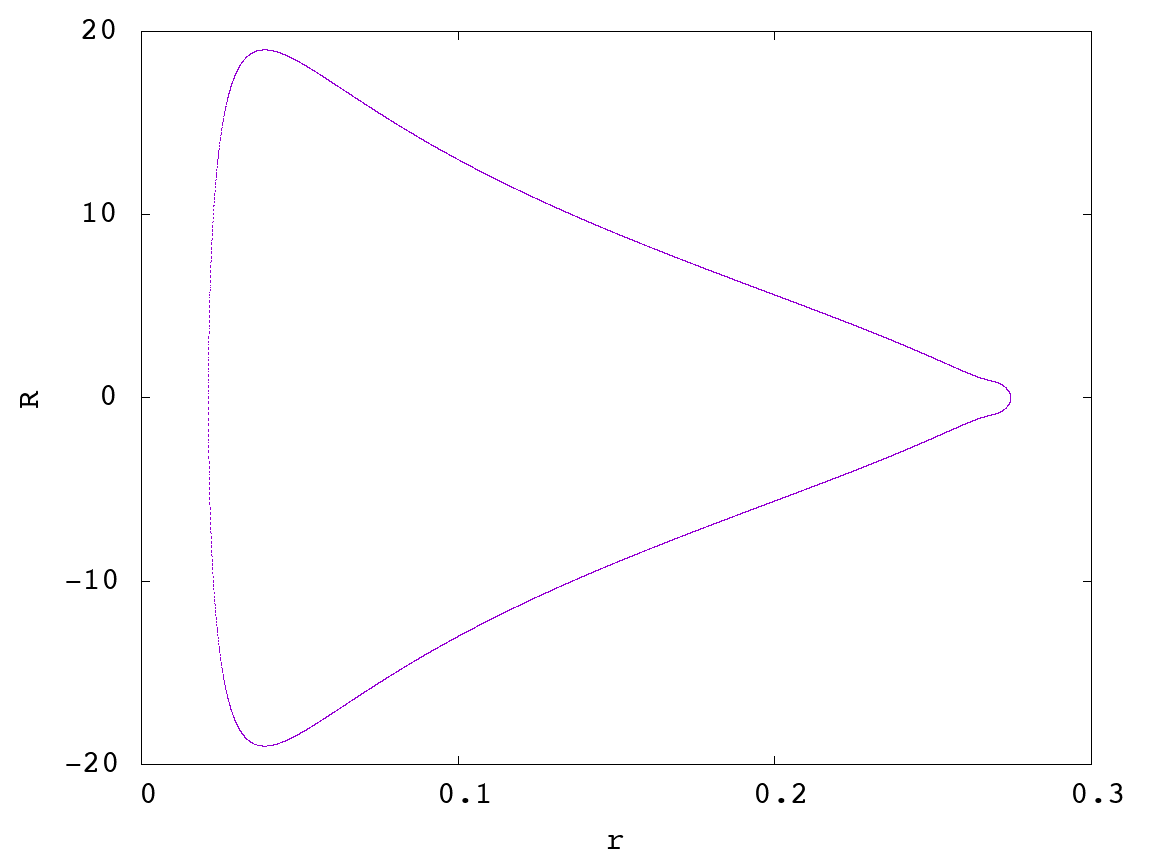}
 \includegraphics[width=3.5cm,height=3cm,draft=false]{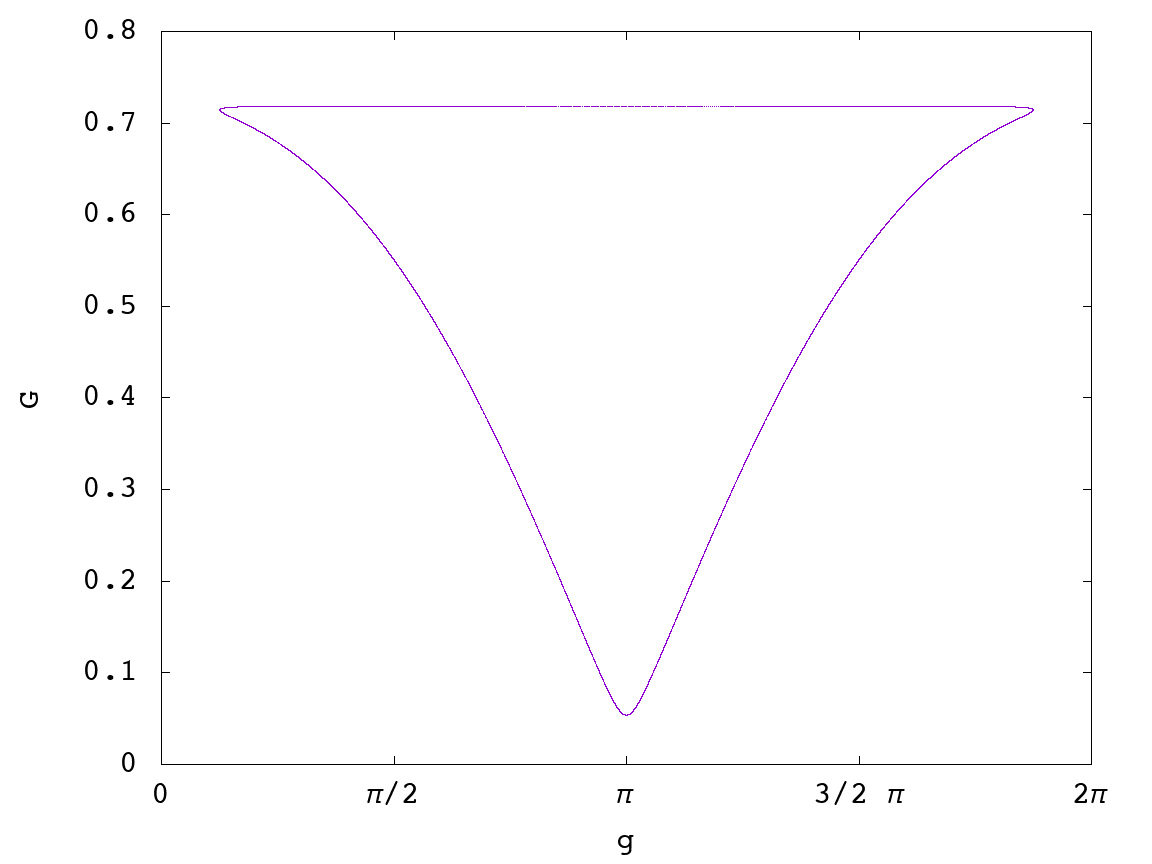}
 \includegraphics[width=3.5cm,height=3cm,draft=false]{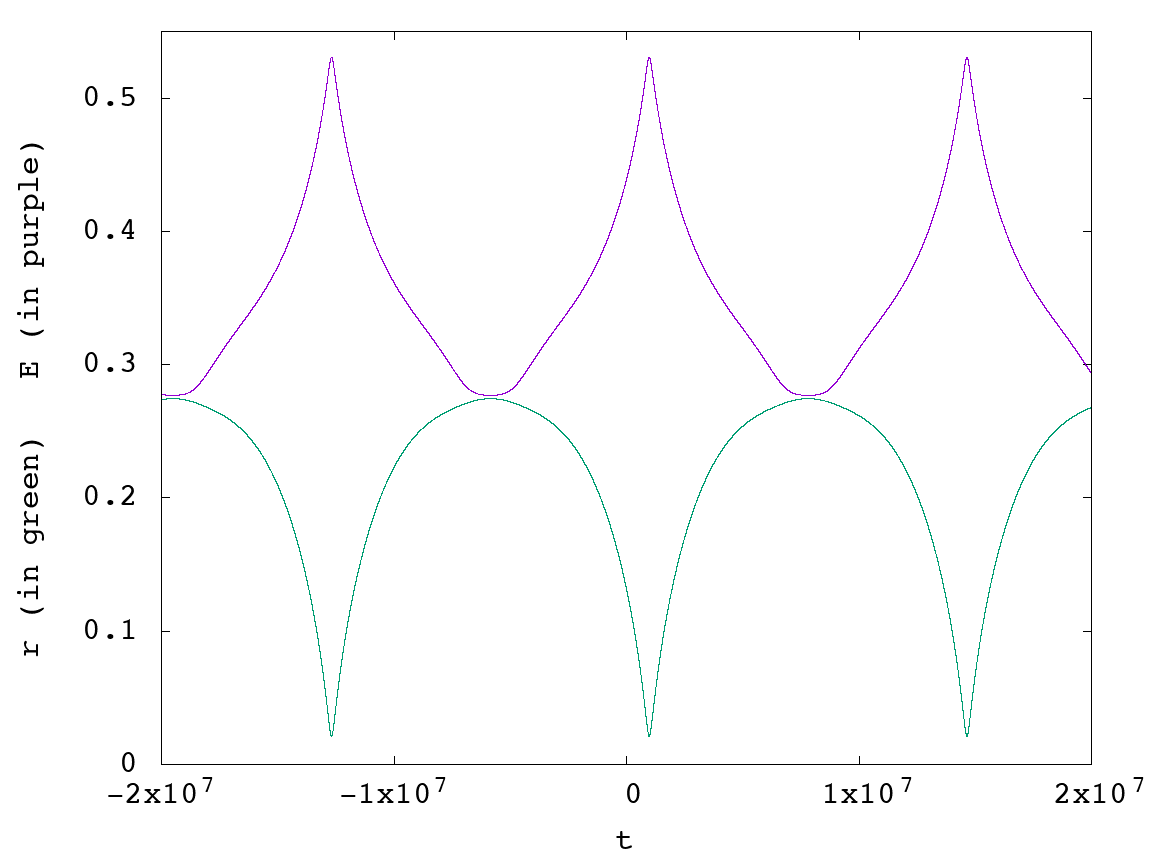}
 \includegraphics[width=3.5cm,height=3cm,draft=false]{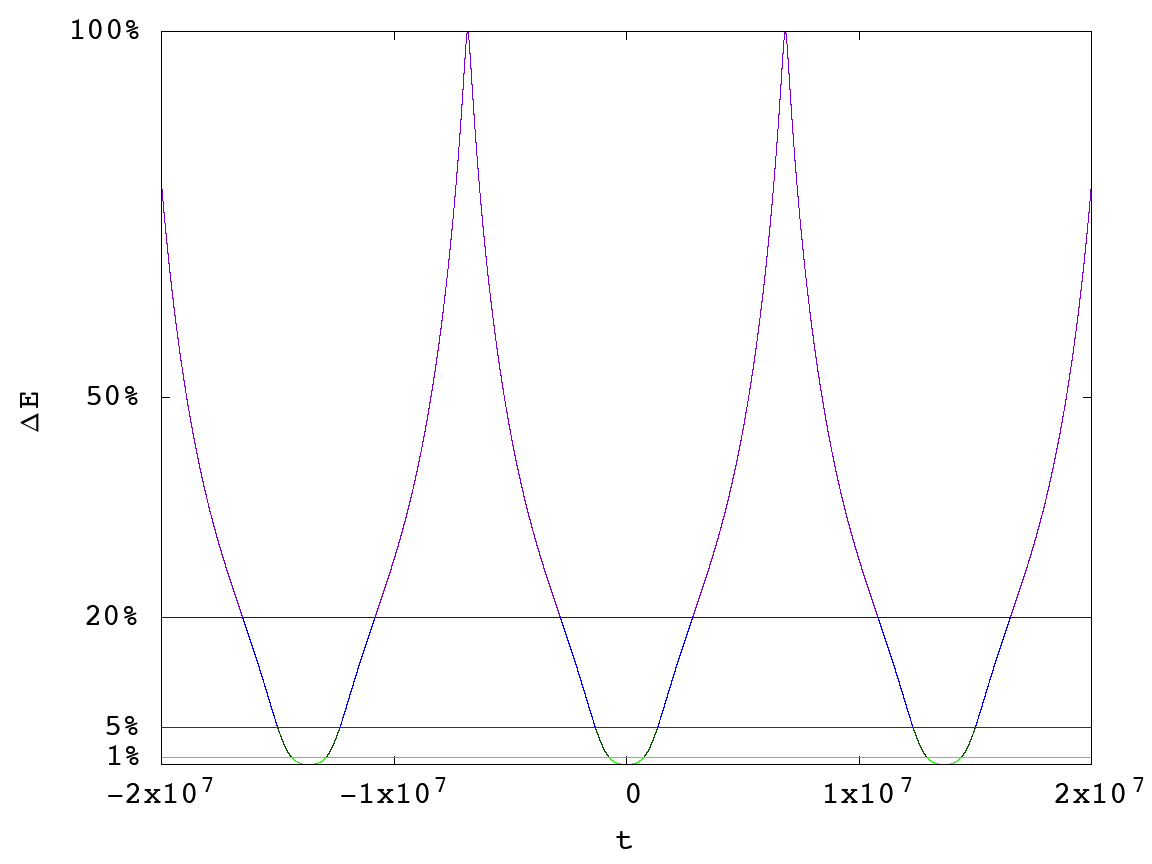}
 \caption{Evolution of the periodic hyperbolic orbit. Upper from left to right, respectively, variables \rm (R, G, r, g) versus time. Bottom from left to right: orbit in the planes $\rm (R, r), \rm (G, g)$, respectively; variation of the Euler integral (in violet) compared with variable $\rm r$ (in green); percentage variations of Euler integral (in blue less than 20\%, in dark--green less than 5\%, in green less than 1\%). }  \label{335_orbit}
 \end{figure}
\begin{figure}[H]
\centering
\includegraphics[width=7.5cm,height=5cm,draft=false]{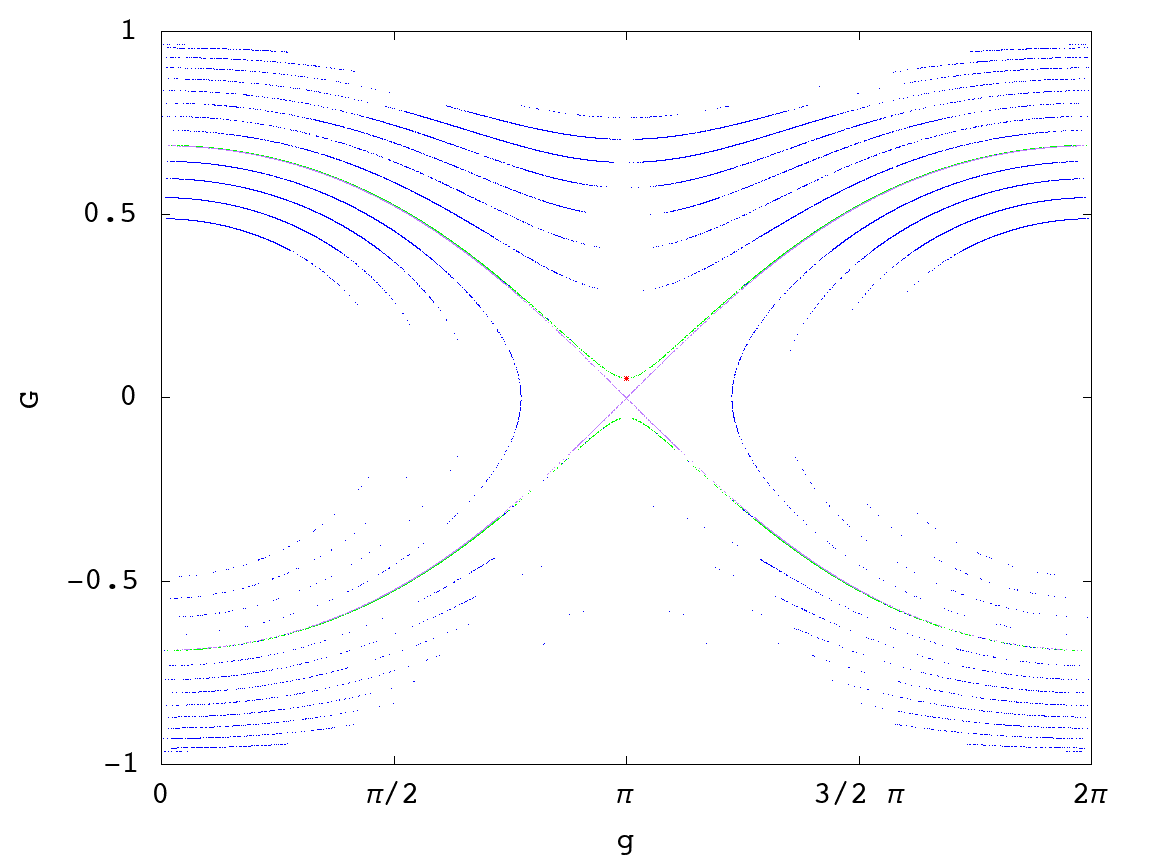}\\
\includegraphics[width=7.5cm,height=5cm,draft=false]{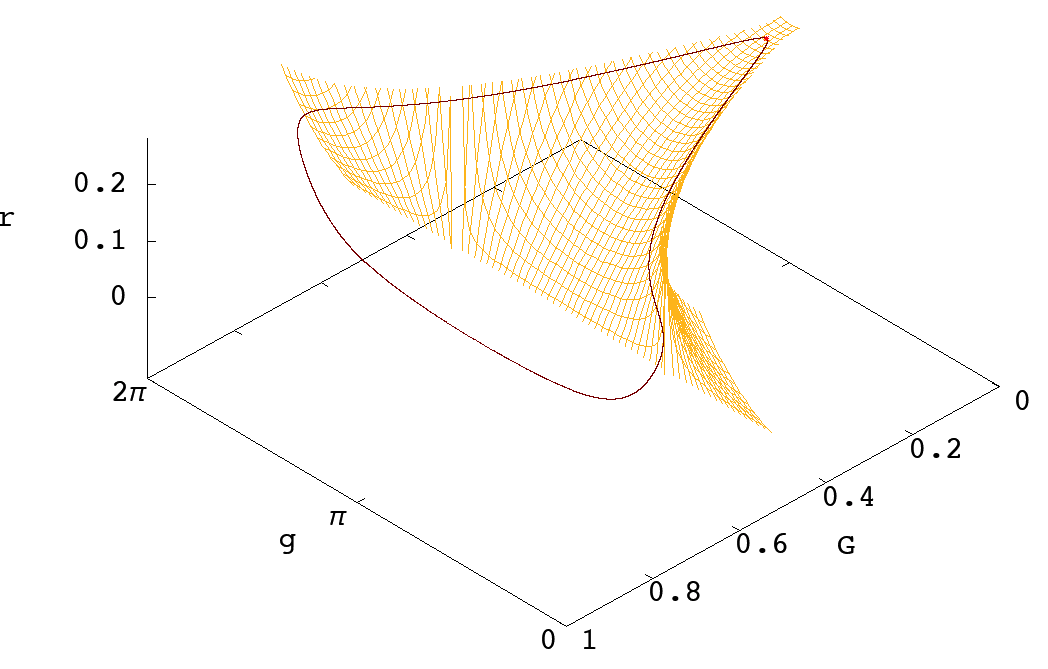}
\includegraphics[width=7.5cm,height=5cm,draft=false]{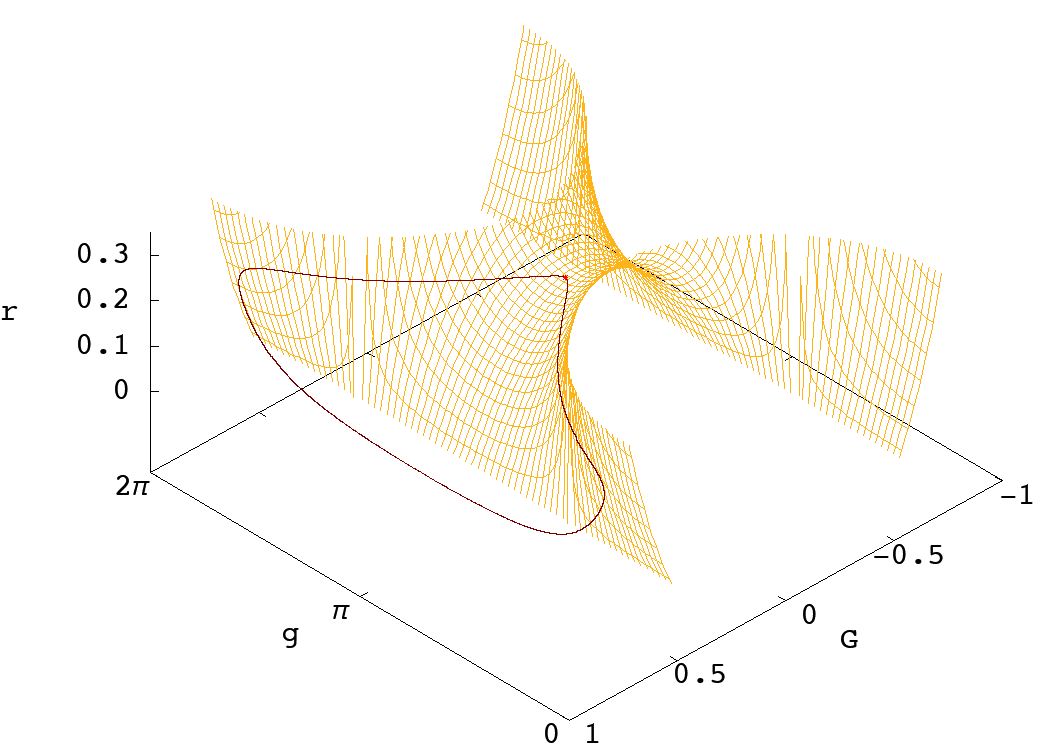}
\caption{Top: curves ${\cal S}({\rm r}_0, \varrho\EE_{\rm min})$ intersected by the $\overline\HH_\CC$ evolution, with initial data $\rm (r_0, G_0, g_0)\in{\cal M}(\EE_{\rm min})$.  Bottom: position of $\Gamma_{\rm u}$ relatively to the manifold ${\cal M}(\EE_{\rm min})$.}\label{contour_iperb}
\end{figure}

\begin{figure}[H]
\centering
\includegraphics[width=7.5cm,height=5cm,draft=false]{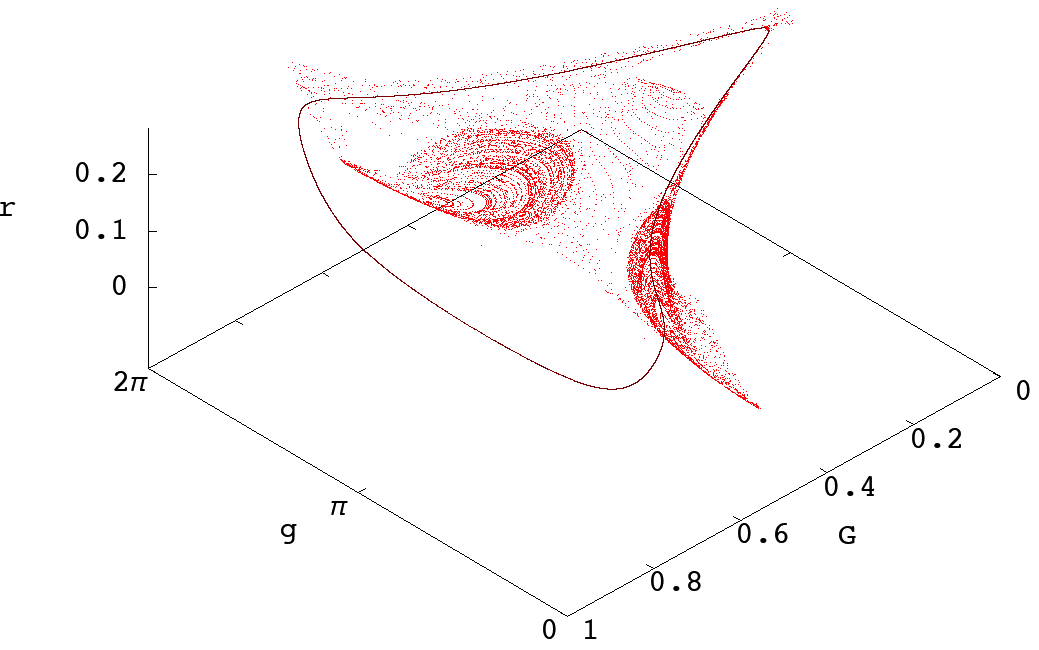}
\includegraphics[width=7.5cm,height=5cm,draft=false]{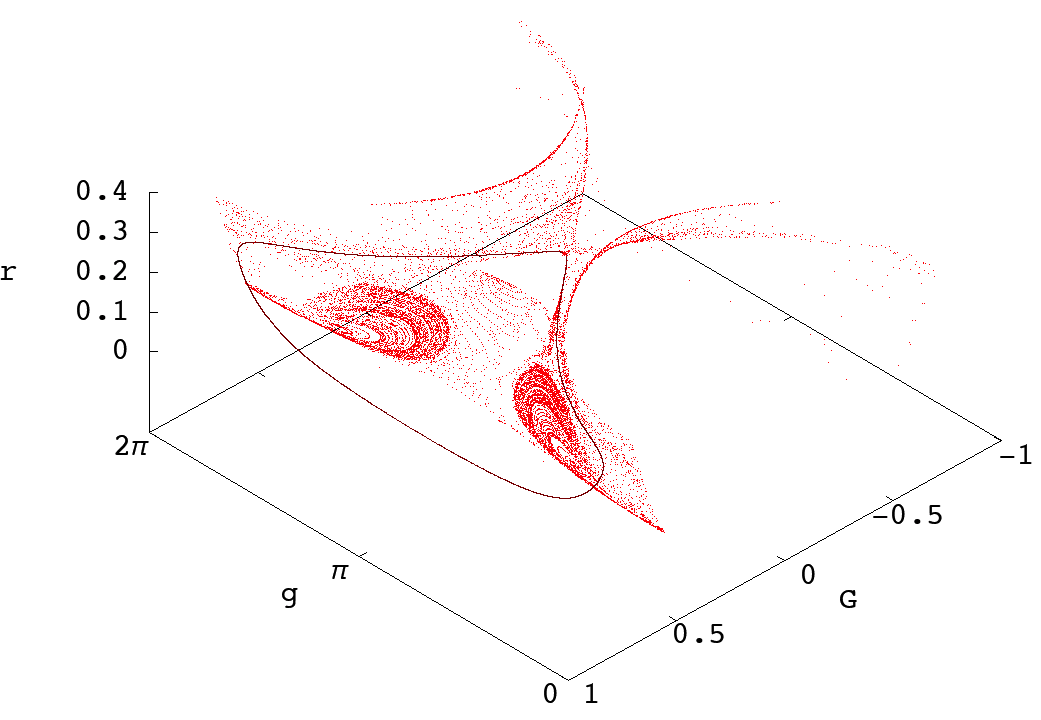}\\
\includegraphics[width=7.5cm,height=5cm,draft=false]{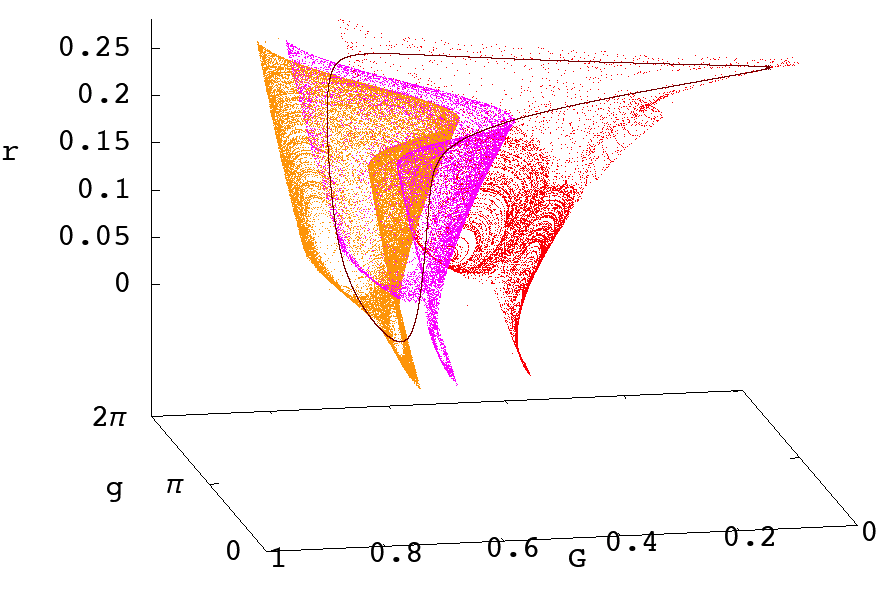}
\includegraphics[width=7.5cm,height=5cm,draft=false]{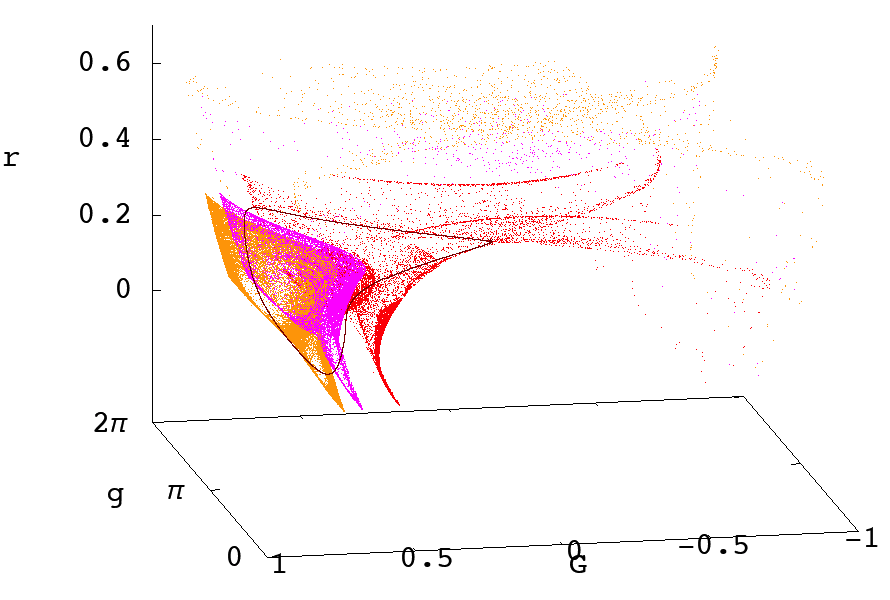}
	\caption{Top: details on the positions of $\Gamma_{\rm u}$ and the part of ${\cal M}(\EE_{\rm min})$ intersected by orbits evolving from ${\cal M}(\EE_{\rm min})$ under $\overline\HH_\CC$.  Bottom: manifolds ${\cal M}({\cal E})$ intersected by the $\overline\HH_\CC$ evolution, with initial data $\rm (r_0, G_0, g_0)\in\cal M{(\EE_{\rm min})}$, with ${\cal E}\in\big\{\EE_{\rm min}\,,\ \EE_{\rm mid}\,,\ \EE_{\rm max}\big\}$} \label{hyp_man} 
\end{figure}

\begin{figure}[H]
\centering
\includegraphics[width=7.5cm,height=5cm,draft=false]{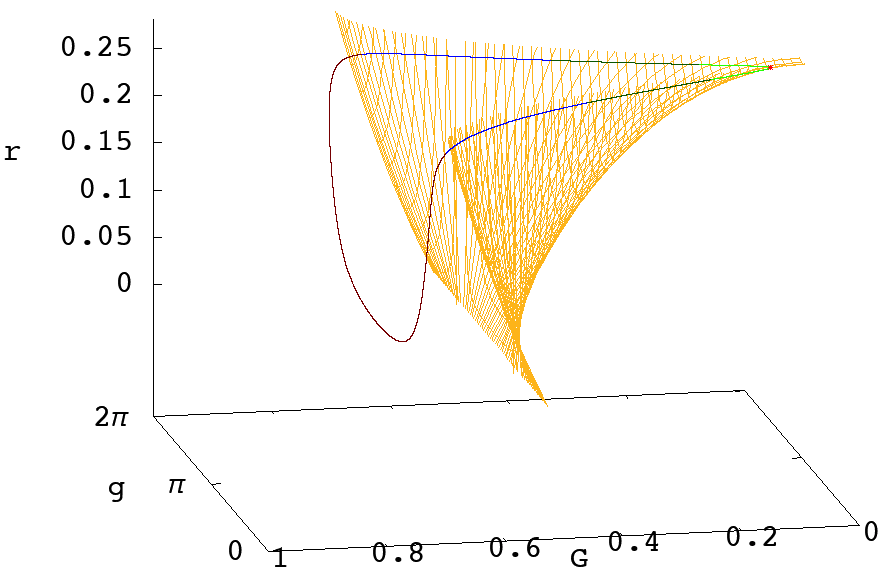}
	\caption{Variations of Euler integral (in blue less than 20\%, in dark--green less than 5\%, in green less than 1\%) along the orbits $\Gamma_{\rm u}$ related to ${\cal M}(\EE_{\rm min})$}. \label{eps1_eps2} 
\end{figure}

\begin{figure}[H]
\centering
\includegraphics[width=15cm,height=10cm,draft=false]{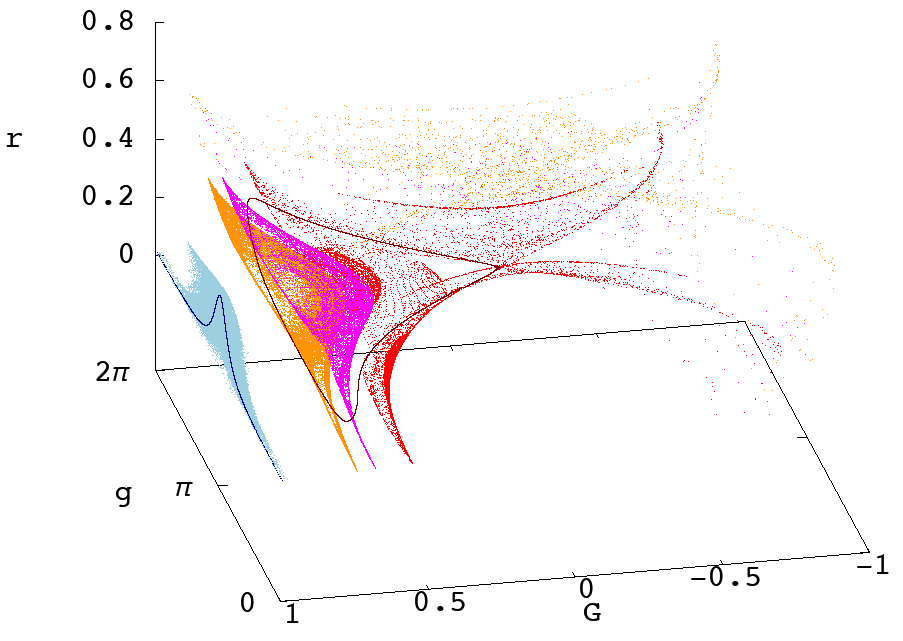}
	\caption{Comparison between the spread of $\EE$ associated to $\Gamma_{\rm s}$ (light blue) and $\Gamma_{\rm u}$ (orange to red).} \label{ell_hyp} 
\end{figure}

  
\noindent The time variations of the coordinates $\rm R$, $\rm G$, $\rm r$ and $\rm g$ along the hyperbolic periodic orbit are depicted in the upper panel of Figure~\ref{335_orbit}. The bottom panel, from left to right, shows the hyperbolic orbit in the planes $\rm (R, r)$ and $\rm (G, g)$; the time variation of $\rm E$ along the orbit compared with the variation of variable $\rm r$. Due to the variation of the velocity along the orbit, $\rm E$ spends most of time \lq\lq close\rq\rq \, to its initial value. With this we mean that, if $T$ the period of the orbit, and $T_{20\%}$,$T_{5\%}$, $T_{1\%}$ the time when, respectively, the variation of $\rm E$ is less than 20\% (straight line blue in right bottom plot in Figure~\ref{335_orbit}), less than 5\% (straight line dark-green in the same plot) and less than 1\% (straight line green), the follow relations hold:
\begin{equation*}
	T_{20\%} = 41.5 \% \cdot T \quad , \quad T_{5\%}= 19.4 \% \cdot T \quad , \quad T_{1\%}= 10.9 \% \cdot T \, \, .
\end{equation*}
Moreover, in Figure~\ref{335_orbit} bottom--third from left, we note that the Euler integral $\rm E$ of the hyperbolic orbit reaches its minimum $\EE_{\rm min}$ when $\rm r$ is maximum and its maximum $\EE_{\rm max}$ when $\rm r$ is minimum. Moreover, the maximum value $\rm r_{\rm max}$ of $\rm r$ is slightly  less than the minimum $\rm E_{\rm min}$ of $\rm E$, as in fact
$$\rm r_{\rm max}= \varrho_0 \EE_{\rm min}\,,\qquad  \varrho_0=0.991\,.$$
 In particular,  $\rm E$ reaches its minimum along the orbit in a region of phase spase which is very close to ${\cal S}_0(\rm r_{\rm max})$. Analogously to the case of the elliptic orbit,   we plot, in the plane $\rm (G, g)$,  level curves  ${\cal S}(\rm r_{\rm max}, \varrho \cdot \EE_{min})$ intersected by the $\overline\HH_\CC$--evolution under a grid of initial values $\rm (r_0, G_0, g_0)\in {\cal M}(\EE_{min})$. We find that 
 \begin{numevid}
 The only non void level manifolds ${\cal M}(\varrho\EE_{min})$ intersected by
orbits of $\overline\HH_\CC$ with initial data in ${\cal M}(\EE_{min})$
 are those with
 \begin{equation}\label{varrho}\varrho\in [0\,,\ 3.1] \,.\quad \square
 \end{equation} 
 \end{numevid}
 We report the results in Figure~\ref{contour_iperb}, top.  By construction, the value (marked in red) of $\rm (G, g)$ on $\Gamma_{\rm u}$ at the  time when $\rm E=E_{\rm min}$ belongs to the curve (plotted in green) with  $\varrho=1$, while ${\cal S}_0(\rm r_{\rm max})$ is obtained for $\varrho=0.991$ (plotted in purple). In blue, we plot curves for different values of $\varrho$ in~\eqref{varrho}. For comparison, in Figure~\ref{contour_iperb}, bottom, we report the position of $\Gamma_{\rm u}$ relatively to the manifolds ${\cal M}(\rm E_{\rm min})$ (yellow). We note that $\varrho\EE_{min}$ with $\varrho= 3.1$ corresponds at ${\cal E} = 0.858673$ which is less than the minimum admissible value provided in  Numerical Evidence~\ref{numevid_gammas}.
\noindent
Finally, in order to measure how the level manifolds ${\cal M}(\cal E)$ which are ``touched'' by $\Gamma_{\rm u}$ spread under $\overline\HH_{\rm C}$, we proceed as follows. We pick a box of initial values $\rm(R_0, G_0, r_0,  g_0)$ on the energy level~\eqref{energy} and  $\rm(r_0, G_0, g_0)\in \cal M{(\EE_{\rm min})}$ and let the system evolve under $\overline\HH_{\rm C}$. We mark a point on the three--dimensional space $\rm (r, G, g)$ whenever the orbits intersects ${\cal M}(\cal E)$ , with
\begin{equation}\label{varrho1}
\cal E\in\big\{\EE_{\rm min}\,,\ \EE_{\rm mid}\,,\ \EE_{\rm max}\big\}
\end{equation}
where $\EE_{\rm min}$, $\EE_{\rm mid}$, $\EE_{\rm max}$ are, respectively, the minimum (0.276991), middle (0.438944), maximum (0.530668) value of $\EE$ along $\Gamma_{\rm u}$ (see Figure~\ref{335_orbit}). We obtain the picture in Figure~\ref{hyp_man}, bottom, where the orange is for  ${\cal M}(\EE_{\rm min})$, magenta is for ${\cal M}(\EE_{\rm mid})$ and red for ${\cal M}(\EE_{\rm max})$.
\noindent In the Figure~\ref{eps1_eps2}, $\Gamma_{\rm u}$  is plotted and we report in blue the part of the orbit such that $\rm E$ --along the orbit-- varies  less than $20\%$ compared to its minimum value $\EE_{min}$, in dark--green the part with variation less than $5\%$ and in green the part with variation less than $1\%$ (compare with last plot in Figure~\ref{335_orbit}).
\noindent
We conclude this section with a visualisation of  the ``spread of $\EE$'' about $\Gamma_{\rm s}$ and $\Gamma_{\rm u}$. In Figure~\ref{ell_hyp} we plot manifolds ${\cal M}({\cal E})$ intersected under the $\overline\HH_\CC$--evolution, with initial data in ${\cal M}({\cal E})$, 
	with $\cal E$ as in Figure~\ref{fase_rfix_cost} 
 (light blue) and~\eqref{varrho1} (orange to red).

\section{Neighbourhoods of $\Gamma_{\rm u}$}\label{Neighbourhoods}

In this section, we show several 2--dimensional maps associated with $\Gamma_{\rm u}$, but constructed in different ways, and detect chaotic phenomena.

\begin{figure}[H]
\centering
\includegraphics[width=7.5cm,height=5cm,draft=false]{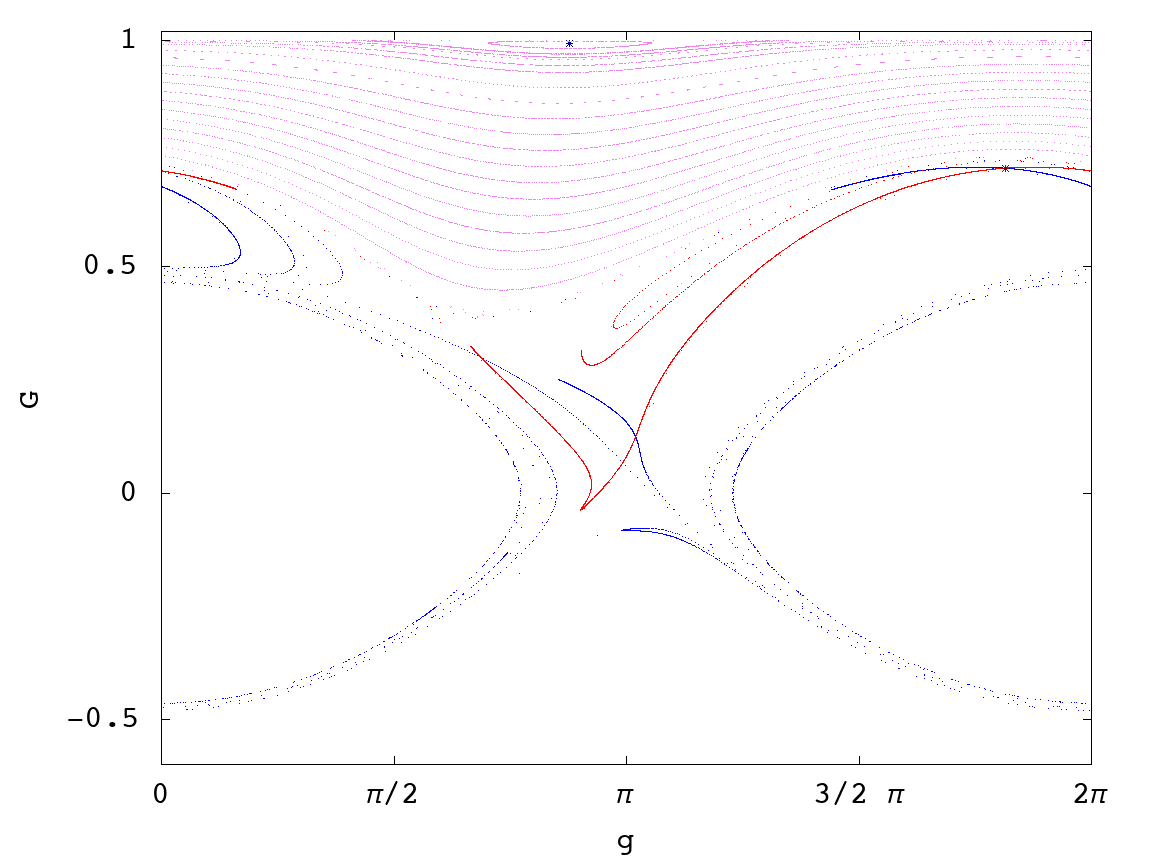}
\includegraphics[width=7.5cm,height=5cm,draft=false]{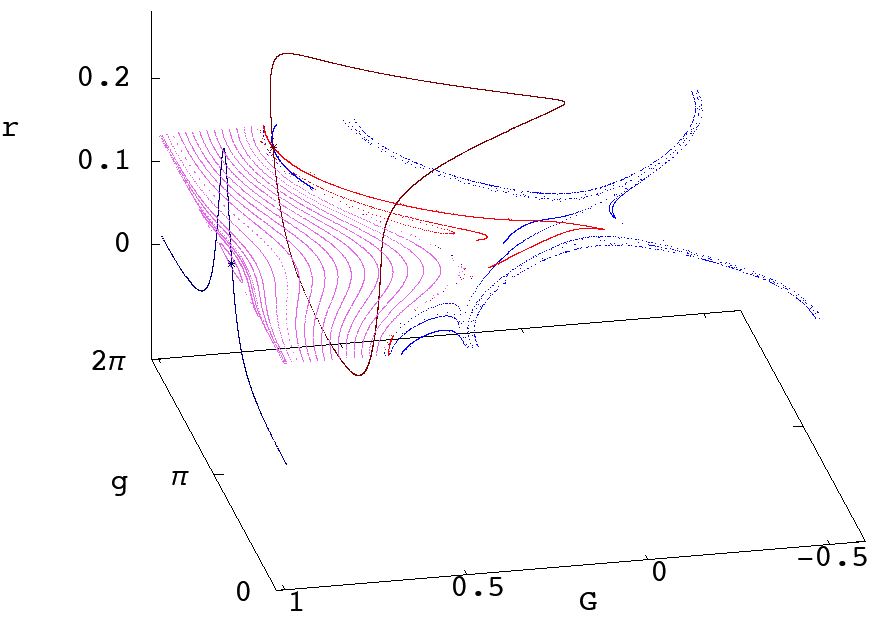}
	\caption{Plane (left) and spatial (right) visualisation of ${\cal P}_{\!\textrm{\tiny$\cal H$,$\Pi_{\rm s}$}}$. 
	 } \label{st_unst_1} 
\end{figure}
\vskip.2in
\noindent
1. The first map is ${\cal P}_{\!\textrm{\tiny$\cal H$,$\Pi_{\rm s}$}}$,  in~\eqref{poinc}, defined, we recall,  as the first return map on a plane $\Pi_{\rm s}$ orthogonal to the orbit $\Gamma_{\rm s}$ in Figure~\ref{fase_rfix_cost}  at the point $\rm (r_s, G_s, g_s)$, defined as in~\eqref{stabledatum}. As outlined in the Numerical Evidence~\ref{fixed points poinc}, the map ${\cal P}_{\!\textrm{\tiny$\cal H$,$\Pi_{\rm s}$}}$ does not show other fixed point than $\rm (G_s, g_s)$ and $\rm (G_u, g_u)$. We find the following
\begin{numevid}[transverse homoclinic intersection for ${\cal P}_{\!\textrm{\tiny$\cal H$,$\Pi_{\rm s}$}}$]
The stable\footnote{We recall that local stable and unstable manifolds associated to the hyperbolic fixed point $\bx^*$ of a map $\cal P$ are defined as 
 \begin{eqnarray*}
 	\WW^s_{loc} & = & \Big\{x \  \Big | \ \|
 	{\cal P}^n(x)-x^*\| \rightarrow 0,\, \ n\in \N_+ \ , \ n \rightarrow \infty \Big \} \, , \\
 	\WW^u_{loc} & = & \Big\{x \  \Big | \ \|{\cal P}^{-n}(x)-x^*\| \rightarrow 0,\, \ n\in \N_+ \ , \ n \rightarrow \infty \Big \} \, .
 \end{eqnarray*}} and unstable manifolds departing from  and arriving at $\rm (G_u, g_u)$ under ${\cal P}_{\!\textrm{\tiny$\cal H$,$\Pi_{\rm s}$}}$ 	have transverse intersection.
$\quad\square$\end{numevid}
At this respect, in Figure~\ref{st_unst_1} the following objects are visible: \begin{itemize}
	\item[--] the elliptic (dark--blue) and the hyperbolic (dark--red) fixed points; \item[--] rotational tori (purple);
	 \item[--] chaotic motions (dotted purple);
	 \item[--] the transverse homoclinic intersection between the stable (blue) and unstable (red) manifolds from $\rm (G_u, g_u)$.
	 \end{itemize}

\begin{figure}[H]
\centering
\includegraphics[width=5cm,height=4cm,draft=false]{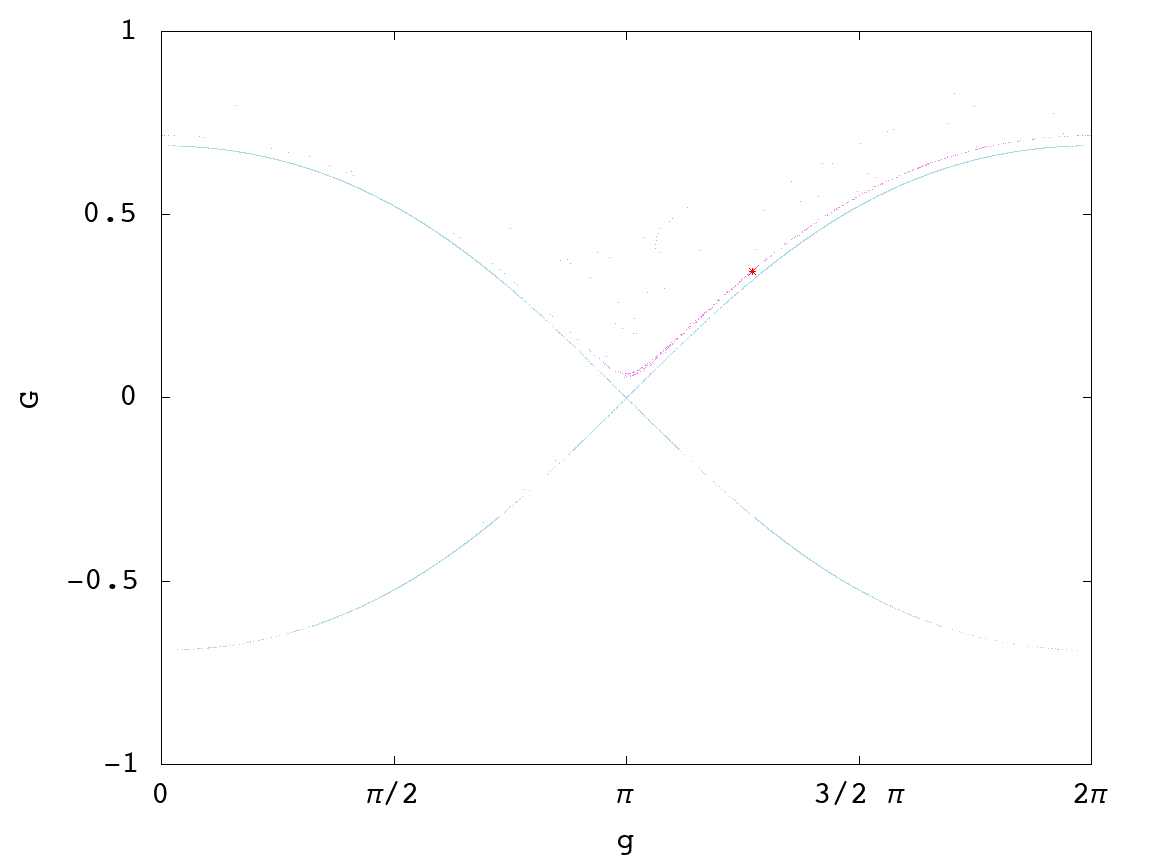}
\includegraphics[width=5cm,height=4cm,draft=false]{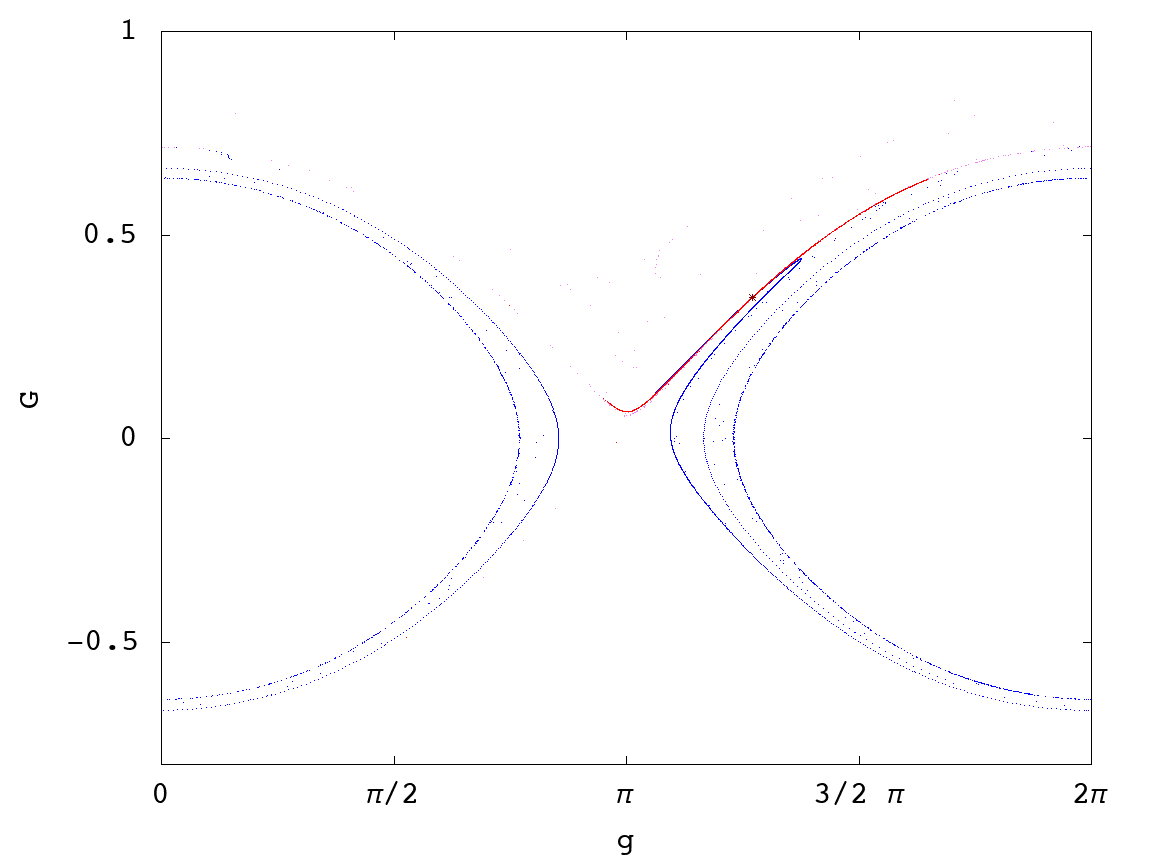}
\includegraphics[width=5cm,height=4cm,draft=false]{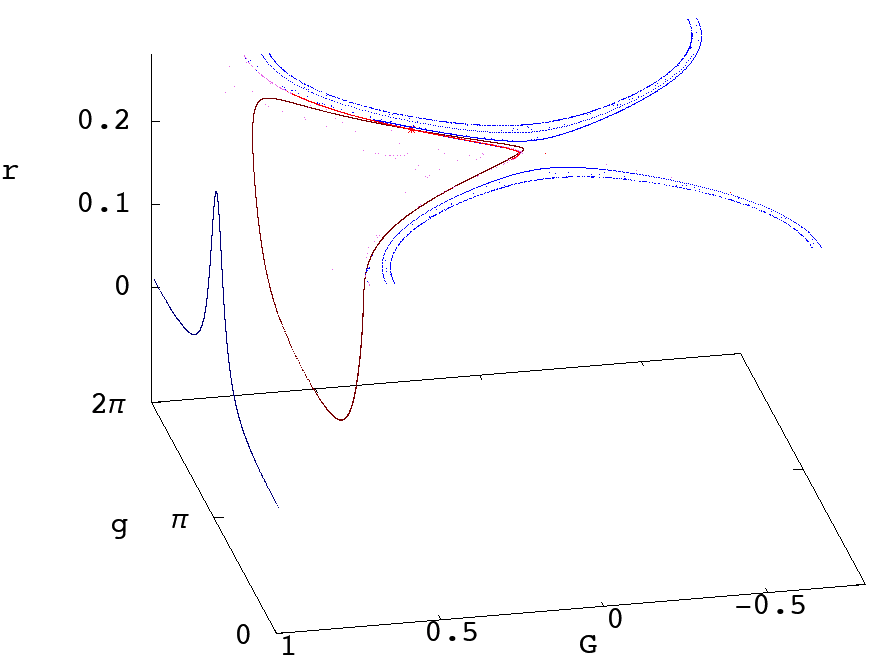}
	\caption{First return map on the plane $\Pi_*$ orthogonal to $\Gamma_{\rm u}$ at  $\rm (r_*, G_*, g_*)$.} \label{st_unst_2} 
\end{figure}

\begin{figure}[H]
\centering
\centering
\centering
\includegraphics[width=7.5cm,height=5cm,draft=false]{figures/fig47_Poi_2_man}
\includegraphics[width=7.5cm,height=5cm,draft=false]{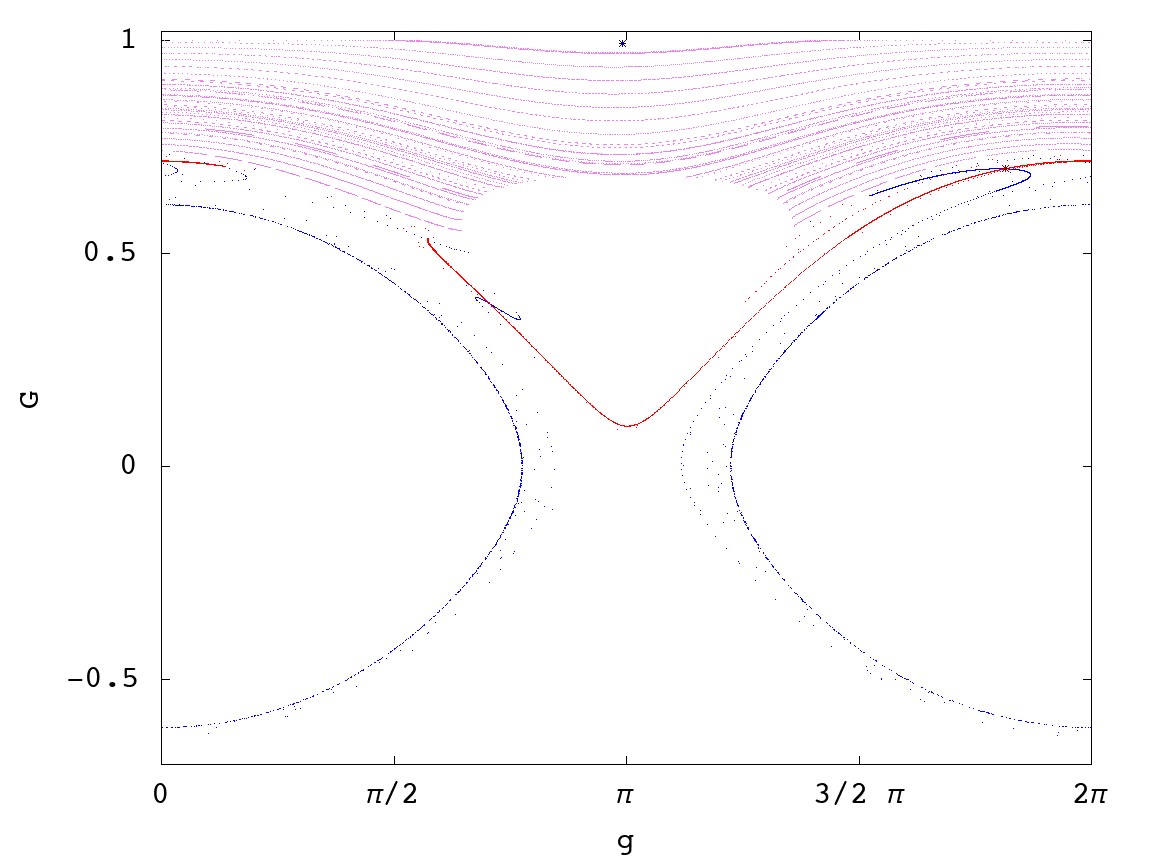}\\
\includegraphics[width=7.5cm,height=5cm,draft=false]{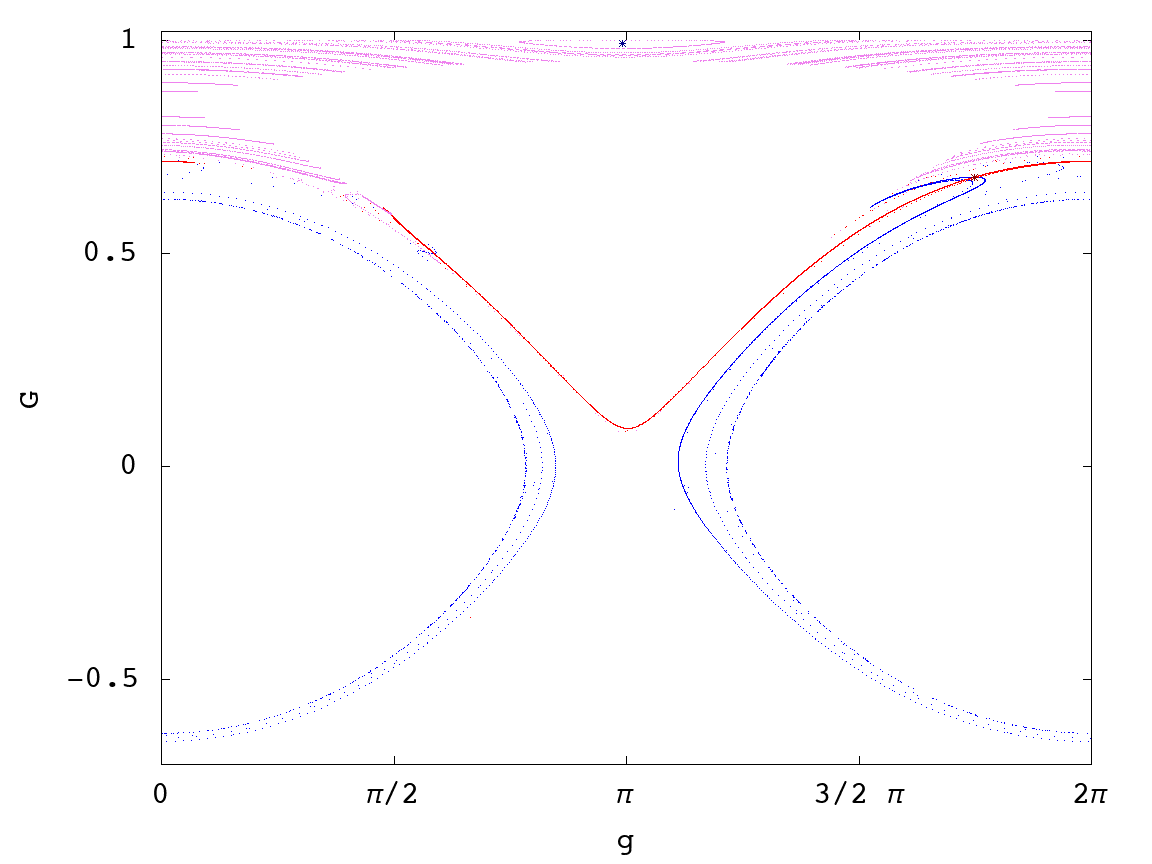}
\includegraphics[width=7.5cm,height=5cm,draft=false]{figures/fig50_Poi_3_man}
	\caption{First return maps on the planes $\Pi_u^i$ orthogonal to $\Gamma_{\rm u}$ at  $\rm (r_i, G_i, g_i)\in \Gamma_{\rm u}$, for $\rm i=1, \dots,4$ where, from upper left to bottom right, respectively, $\rm r_1 =0.13165\, ,\,r_2 = 0.242432 \, ,\,r_3 = 0.252024 \, ,\,r_4 = r_* = 0.26987 $  (we recall that $\rm r_{max}= 0.274496$). } \label{evi5p2} 
\end{figure}

\vskip.2in
\noindent
2.
Let $\Pi_{\rm u}^i$ be different planes  orthogonal to $\Gamma_{\rm u}$ at different points $({\rm r}_{\it i}, {\rm G}_{\it i}, {\rm g}_{\it i})$ of the curve. We consider first return maps 
\begin{equation}\label{poinc2}{\cal P}_{\!\textrm{\tiny$\cal H$,$\Pi_{\rm u}^i$}}:\qquad ({\rm G}, {\rm g})\to ({\rm G}', {\rm g}')
\end{equation} 
along $\Pi_{\rm u}^i$. Incidentally, this procedure provides us with an help to control numerical errors, as we check the invariance of the Lyapunov exponents at $\rm (G^{\it i}_*, g^{\it i}_*)$, for different choices of $\Pi_{\rm u}^i$.\\
 We denote as $\Pi_*$ the plane orthogonal to $\Gamma_{\rm u}$ at \begin{equation}\label{saddlepoint}\rm(r_*, G_*, g_*)=(0.270, 0.346, 1.27\,\pi)\end{equation}
 This point of  $\Gamma_{\rm u}$ has been chosen for being  ``close'' to  ${\cal M}_0$.
With this choice, we detect a homoclinic tangency and absence of splitting:

\begin{numevid}[quasi--homoclinic tangency]
Consider first return maps 
\eqref{poinc2}
on $\Pi_{\rm u}^i$. As soon as  $\rm r^{\it i}_*$ is chosen closer and closer to ${\rm r}_{\rm max}$ the stable tori zone becomes smaller and smaller.  For $\Pi_{\rm u}^i=\Pi_*$, stable motions are not numerically detected, and the unstable, stable manifolds have a homoclinic tangency at $\rm (G_*, g_*)$. In other words a splitting of such manifolds (which have the shape of ${\cal S}_0(\rm r_*)$)  is not numerically detected. 
$\quad\square$\end{numevid}
The normalized stable and unstable eigenvectors in $\rm (G_*, g_*)$ are, respectively, 
\begin{equation*}
	v^s_*  = ( 0.312937 \, ,\,  0.949774)  \,, \quad 
	v^u_*  =  (0.320019 \, , \, 0.947411,) 
	\end{equation*}
and the angle between them is
\begin{equation*}
	\alpha =  0.0074651.
\end{equation*}
\noindent
The results are visualised in Figure~\ref{st_unst_2} and~\ref{evi5p2}. 

\begin{figure}[H]
\centering
\includegraphics[width=5cm,height=4cm,draft=false]{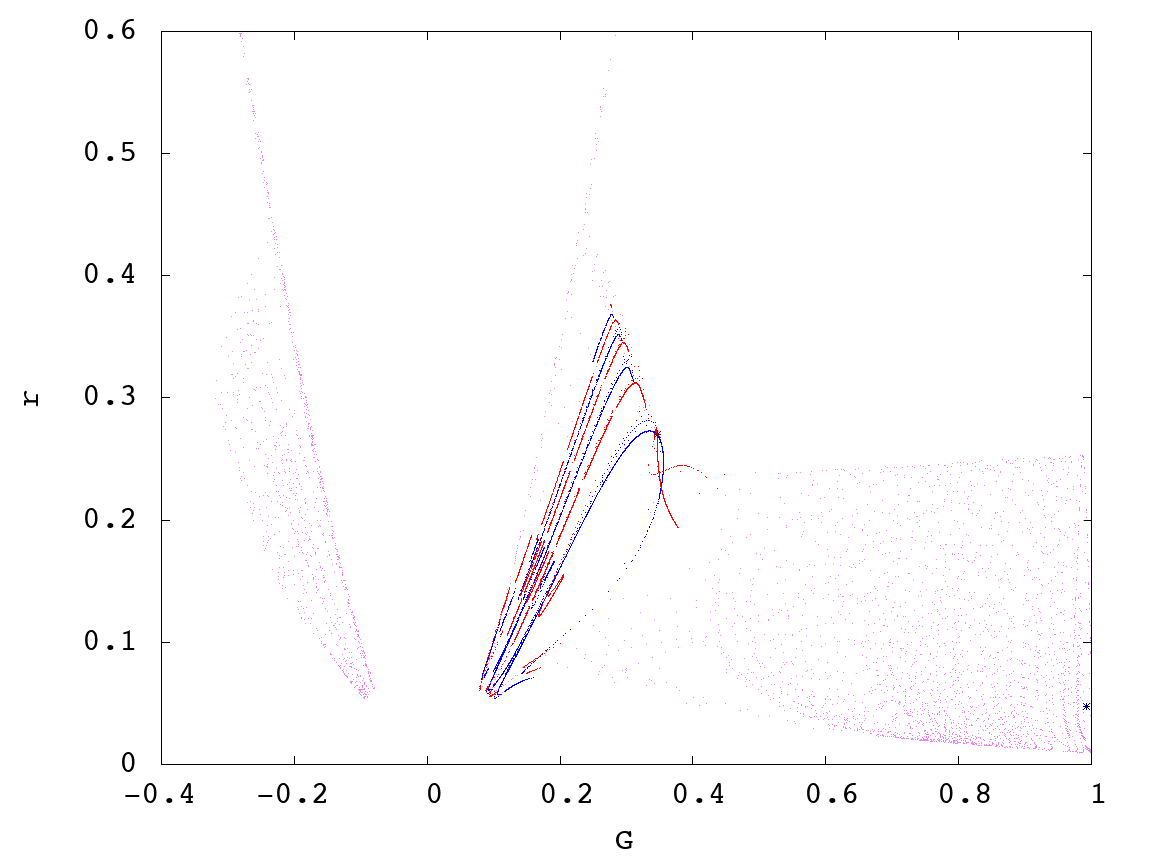}
\includegraphics[width=5cm,height=4cm,draft=false]{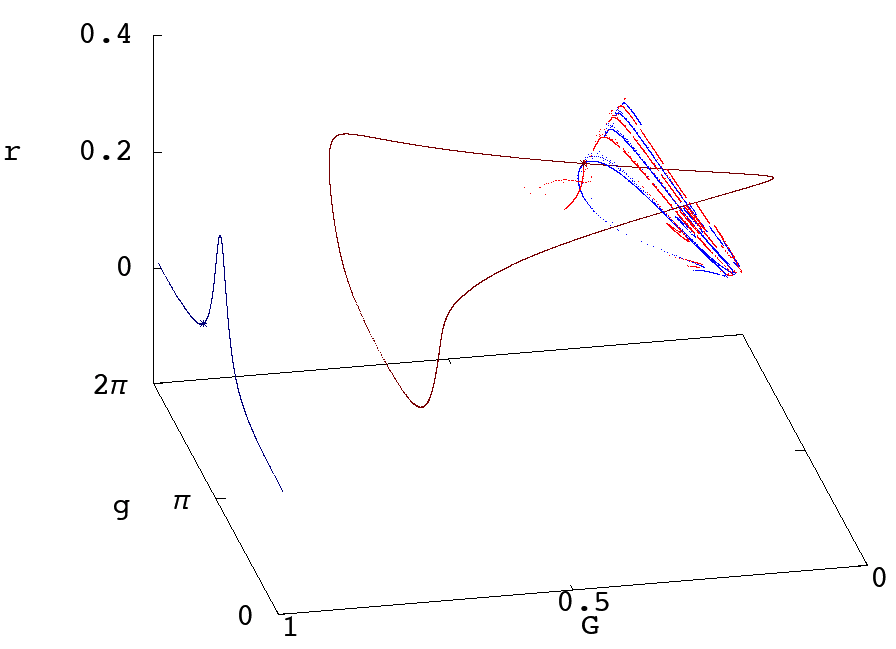}
\includegraphics[width=5cm,height=4cm,draft=false]{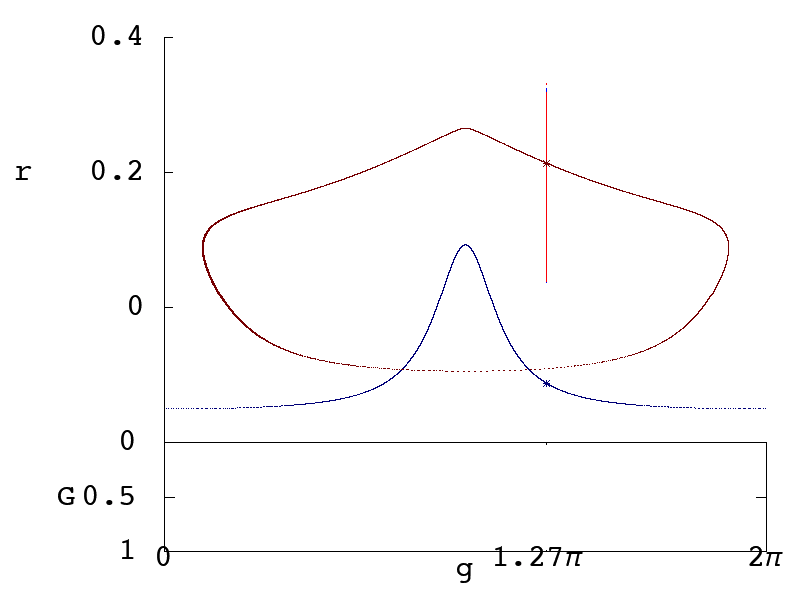}
	\caption{First return map on the plane $\Pi^*=\{\rm g=g_*\}$. Note the aspect of the stable and unstable manifolds in the first panel.} \label{st_unst_2bis} 
\end{figure}

\vskip.2in
\noindent
3. We finally fix the plane $\Pi^*=\{\rm g=g_*\}$. The two--dimensional first return map
\begin{equation}\label{poinc3}{\cal P}_{\!\textrm{\tiny$\cal H$,$\Pi^*$}}:\qquad ({\rm r}, {\rm G})\to ({\rm r}', {\rm G}')
\end{equation}
is depicted in Figure~\ref{st_unst_2bis}. The aspect of the stable and unstable manifolds changes drastically, but  homoclinic intersections are present.

\begin{figure}[H]
\centering
\includegraphics[width=7.5cm,height=5cm,draft=false]{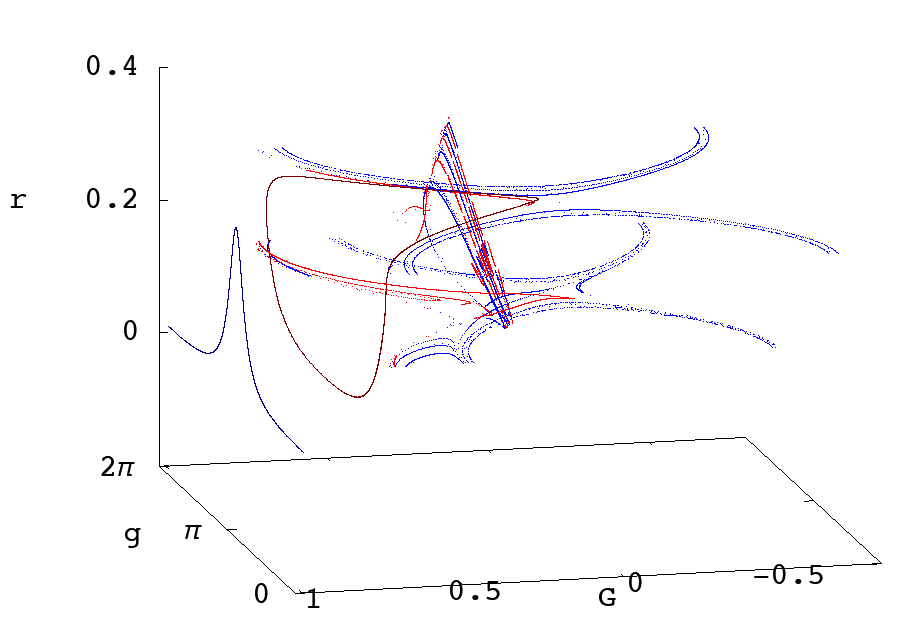}
\includegraphics[width=7.5cm,height=5cm,draft=false]{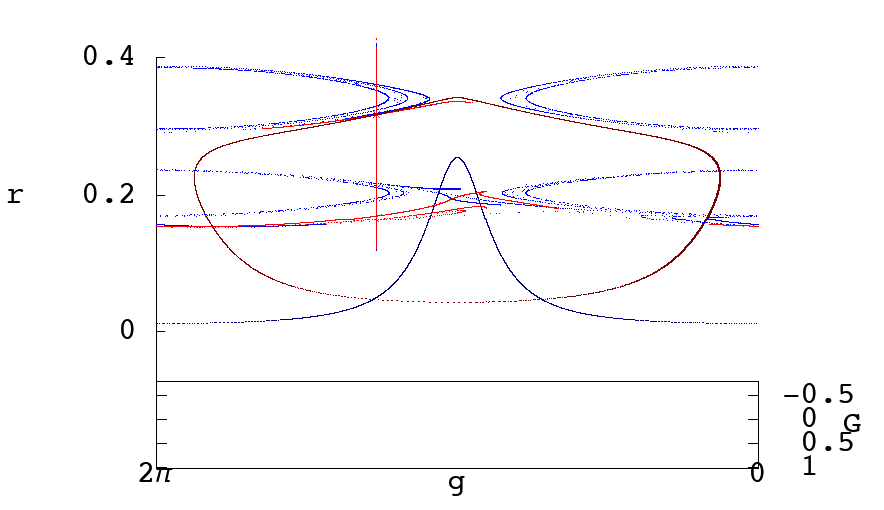}
	\caption{Stable (in blue) and unstable (in red) manifolds constructed using different maps.} \label{st_unst_3} 
\end{figure}

\vskip.2in
\noindent
A comprehensive visualisation including the surface sections and returns maps herewith considered is in Figure~\ref{st_unst_3}.

\begin{figure}[H]
\centering
\includegraphics[width=12cm,height=9cm,draft=false]{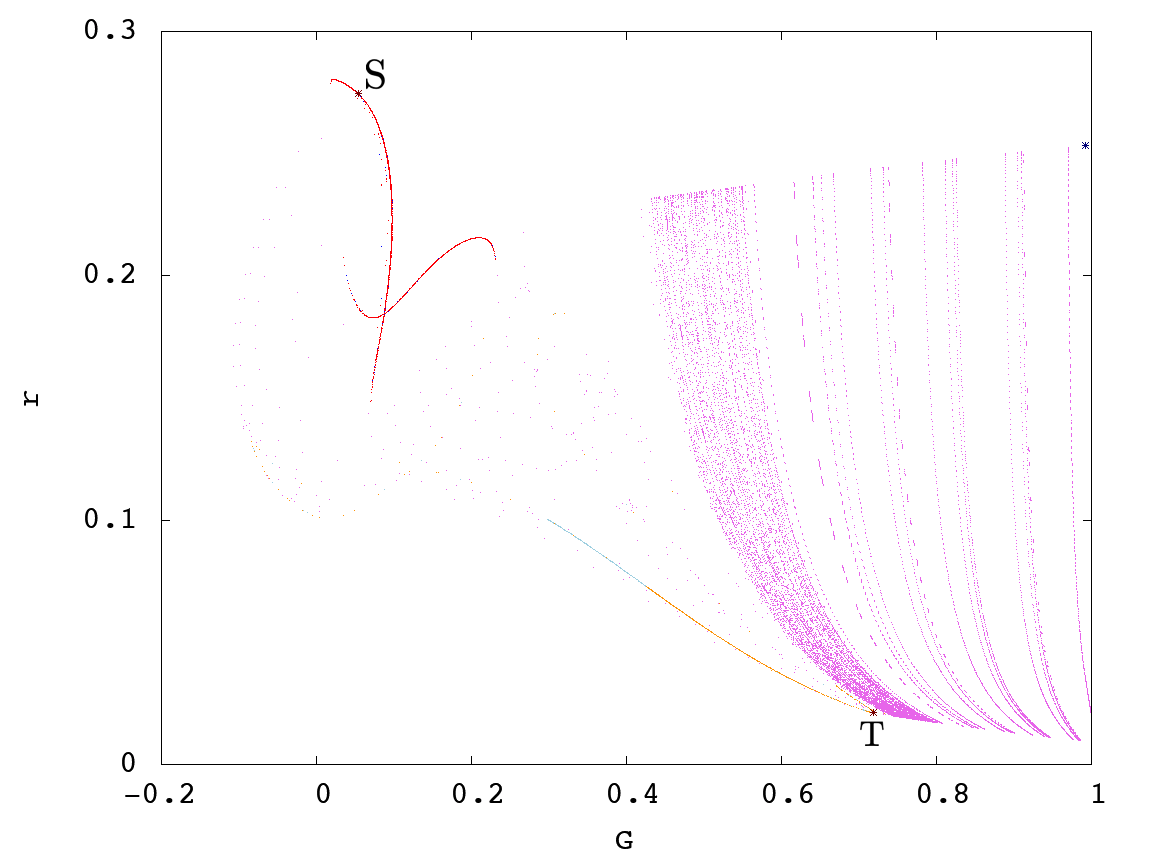}
	\caption{First return map on the plane $\Pi^{\bigstar}$; in purple we plot the orbit sections (regular orbits on the right and chaotic see in the left); in blue and red the stable and unstable manifolds of fixed point $\rm S$; in light--blue and orange the stable and unstable manifolds of fixed point $\rm T$.} \label{poi_pi} 
\end{figure}

\vskip.2in
\noindent
4. As special case of previous point, we fix the plane $\Pi^{\bigstar}=\{\rm g=\pi\}$. Due to the geometrical shape and properties of the curve $\Gamma_{\rm u}$, it is not possible to fix the orthogonal planes at the points $\rm T = (G_{\star},r_{min},\pi) $ and $\rm S =(G^{\star},r_{max},\pi) $. For this reason, we study the structure of the first return map on the plane $\Pi^{\bigstar}$. The points $\rm T,S$ are both hyperbolic fixed points depending on the orbit is run in one direction or in the opposite. We construct stable and stable manifolds for both fixed points and we obtain a complete overlapping of the stable and unstable manifolds of each point. In Figure~\ref{poi_pi}, we can see in red and blue the manifolds related to $\rm S $  (and the blue manifold is completely hidden by the red one) and in orange and light--blue the manifolds related to $\rm T$.

 \section{Symbolic dynamics}\label{Symbolic dynamics}
 \label{section_horseshoe}
 In this section, we discuss numerical evidence of symbolic dynamics for the map ${\cal P}_{\!\textrm{\tiny$\cal H$,$\Pi^*$}}$ in~\eqref{poinc3}, in the sense of the following 
 
 \begin{definition}[Symbolic dynamics; horseshoe]\label{def: symb dyn}\rm
 Let $D\subset{\mathbb R}^2$
 $${f }:\qquad D\to {\mathbb R}^2 \, , $$ 
 we say that ${f }$ has $m$--{\it symbolic dynamics} if
 there exist compact subsets with non--empty interior $N_0$, $N_1\subset D$
 such that,
 for every $n\in \mathbb N$  and any finite sequence $(\sigma_0$, $\ldots$, $\sigma_n)$  of symbols $\sigma_i\in\{0\,,\ 1\}$ having length $n+1$, one can find $x_0\in N_{\sigma_0}$ such that the orbit of $x_0$ under ${f }$, namely, $x_j:={f }^j(x_0)$ is well defined for $j=0$, $\ldots$, $nm$, and $x_{mj}\in N_{\sigma_j}\quad \forall\ j=0\,,\ldots\,,\quad n$.\\
 $1$--symbolic dynamics in $N_0\cup N_1$ is also called {\it horseshoe}. $\quad \square$
 \end{definition}
 \begin{remark}\rm
 Observe that $m$--symbolic dynamics implies $p$--symbolic dynamics for any $p\in \mathbb N$ such that $m|p$.  So, in presence of an horseshoe, $p$--symbolic dynamics holds for any $p\in \mathbb N$. It is also known that a map with a horseshoe is semi--conjugated to a shift 
 $(\sigma_{-1}, \sigma_0, \sigma_1, \cdots)\to (\sigma_{-2}, \sigma_{-1}, \sigma_0, \cdots)$; see, e.g.,~\cite{ZgliczynskiG2004}.
 \end{remark}
 In fact, we have the following
 \begin{numevid}\label{numvid: symbdyn}
 The map ${\cal P}_{\!\textrm{\tiny$\cal H$,$\Pi^*$}}$ in~\eqref{poinc3} has a $3$--symbolic dynamics. Moreover, an orbit $\{x_j\}_{j=0\,\cdots, 3n}$ corresponding to a given sequence $\sigma_0$, $\ldots$, $\sigma_n$, can be chosen to be extendible for $j=0\,\cdots, 3(n+1)$
 and
 periodic, with period $N\in \{1\,,\cdots\,,3(n+1)\}$.
 \end{numevid}

 \noindent
To understand why we assert the Numerical Evidence~\ref{numvid: symbdyn}, we need to recall, below, the  method of {\it covering relations} developed in~\cite{ZgliczynskiG2004}, and already used in~\cite{GierzkiewiczZ19} and, recently, in~\cite{diruzzaDP20} (of course, the interested reader is invited to consult the mentioned literature for more details). It is to be recalled that in~\cite{diruzzaDP20} the method was used to find an horseshoe, while in this paper we  obtain a weaker result (3--symbolic dynamics), seemingly due to the non--existence of heteroclinic connections, as per Numerical Evidence~\ref{fixed points poinc}.

  \subsection*{Covering relations and symbolic dynamics}
 \label{horseshoe}

We simplify the material of~\cite{ZgliczynskiG2004} to the case that the dimension of the space is $2$, as this is needed in our application.

\begin{definition}[$h$--sets,~\cite{GierzkiewiczZ19, ZgliczynskiG2004}]\label{hset}\rm Let $N \subset \R^2$ be a compact set and let $$c_{N} : \R^2 \rightarrow \R^2$$ be an homeomorphism such that $c_N(N) = [-1,1]^2$. 
\item[\rm(i)]  The couple $(N, c_N)$ is called a {\it h--set};  $N$ is called {\it support} of the $h$--set.  
\item[\rm(ii)] Put $$N_c:=[-1,1]^2\,,\quad N_c^-:=\{-1,1\} \times [-1,1]\,,\quad N_c^+:= [-1,1]\times \{-1,1\}$$
and
$$S(N)_c^l :=(-\infty, -1) \times \R\,,\quad S(N)_c^r :=(1,\infty) \times \R\,,\quad N_{c}^{\rm le}=\{-1\}\times [-1,1]\,,\quad N_{c}^{\rm ri}=\{1\}\times [-1,1]$$
 The  sets $$N^-= c_N^{-1}(N_c^-)\,,\quad N^+= c_N^{-1}(N_c^+)\,,\quad N^{\rm le}=c_N^{-1}(N_c^{\rm le})\,,\quad N^{\rm ri}=c_N^{-1}(N_c^{\rm ri}) $$ are called, respectively, the {\it exit set} and the {\it entry set}, while the sets
$$S(N)^l := c_N^{-1}(S(N)_c^l)\,,\quad S(N)^r = c_N^{-1}(S(N)_c^r) $$ are called, respectively, the  {\it left side}, {\it right side}, {\it left edge}, {\it right edge} of $N$.
 $\quad\square$\end{definition}
 The following definition is fitted to the special case (realised in our study) that the unstable manifold has dimension $1$. The simplification compared to the general definition in~\cite{ZgliczynskiG2004, GierzkiewiczZ19} is based on~\cite[Theorem 16]{ZgliczynskiG2004}.
\begin{definition}[Covering relation,~\cite{ZgliczynskiG2004, GierzkiewiczZ19}]\rm
\label{def_covering}
	Let $f : \R^2 \rightarrow \R^2$ be a continuous map and $N$ and $M$ the supports of two $h$--sets. We say that $M$ $f$--covers $N$ and we denote it by $M \stackrel{f} \Longrightarrow N $ if:
	\begin{itemize}
	\item[(1)\phantom{'}] $\exists\, q_0\in [-1, 1]$ such that $f(c^{-1}_{M}([-1, 1]\times \{q_0\}))\subset {\rm int}( S(N)^l \bigcup N \bigcup S(N)^r)$;
	\item[(2)\phantom{'}] $f(M) \bigcap N^+  = \emptyset$;
	\item[(3)\phantom{'}] $f(M^{\rm le})\subset S(N)^{\rm l}\quad {\rm and}\quad f(M^{\rm ri})\subset S(N)^{\rm r}$ or
	\item[(3)'] $f(M^{\rm le})\subset S(N)^{\rm r}\quad {\rm and}\quad f(M^{\rm ri})\subset S(N)^{\rm l}$
\end{itemize}
If $M=N$, we say that $f$ {\it self--covers} $N$.\\	
Conditions (2) and (3) are called, respectively, {\it exit} and {\it entry condition}.
$\quad\square$\end{definition}

\noindent
As in~\cite{GierzkiewiczZ19}, we write
$A  \stackrel{f} \implica B \stackrel{g} \implica C$, etc, $\ldots$,
if $A  \stackrel{f} \implica B$ and $B  \stackrel{g} \implica C$, etc.

\begin{theorem}[\cite{WilczakZ2003}]\label{symbolic dynamics} Let $N_i$, $i=0$, $\ldots$, $k$, be  $h$--sets 	such that $N_0=N_k$. Let $$f_i:\ N_{i-1}\to{\mathbb R} \qquad \forall\ i=1\,,\ldots\,, k\,.$$ be a continuous map such that
$$N_{0} \stackrel{f_1 } \Longrightarrow N_1\stackrel{f_2 }\Longrightarrow \cdots  \stackrel{f_k } \Longrightarrow N_{k}=N_0$$
Then there exists $x_0\in N_0$ such that
\begin{itemize}
\item[{\rm (i)}] $f_i\circ f_{i-1}\circ \cdots \circ f_1(x_0)\in N_{i}\quad \forall\ i=1\,,\ldots\,, k\,;$
\item[{\rm (ii)}] $f_k\circ f_{k-1}\circ \cdots \circ f_1(x_0)= x_0\,.$
\end{itemize}
\end{theorem}

\noindent
We shall use Theorem~\ref{symbolic dynamics} in the following form.

\begin{corollary}\label{cor: symb dyn}
 Let $D\subset{\mathbb R}^2$
 $${f }:\qquad D\to {\mathbb R}^2$$
 and let $N_0$, $N_1$ be $h$--sets in $D$.
Assume that there exist $h$--sets
$M^{(\sigma, \sigma')}_i$, with $i=1$, $\ldots$, $m-1$ and $\sigma$, $\sigma'\in\{0\,, 1\}$,
such that
\begin{eqnarray}\label{ciclo}
N_\sigma \stackrel{f } \Longrightarrow M^{(\sigma, \sigma')}_1 \stackrel{f } \Longrightarrow M^{(\sigma, \sigma')}_2\cdots \stackrel{f } \Longrightarrow M^{(\sigma, \sigma')}_{m-1}\stackrel{f } \Longrightarrow N_{\sigma'}
\quad \forall\ \sigma\,,\ \sigma'\in\{0\,,1\}\,.
\end{eqnarray} 
Then ${f }$ has $m$--symbolic dynamics in $N_0\cup N_1$. In addition, 
 an orbit $x_k$ corresponding, as in Definition~\ref{def: symb dyn}, to a given sequence $\sigma_0$, $\ldots$, $\sigma_n$, can be chosen so that it is well defined for $i=0$, $\ldots$, $(n+1)m$ and, moreover, $x_{(n+1)m}=x_0$.
\end{corollary}

\noindent
{\bf Proof}
Let $n\in \mathbb N$ and $(\sigma_0$, $\ldots$, $\sigma_n)$ a finite sequence of symbols $\sigma_i\in\{0\,,\ 1\}$ having length $n+1$.
Put: 
\begin{eqnarray}\label{ciclo1}
\begin{array}{lllclll}
&&N^0:=N_{\sigma_0}&&&\\\\
&&N^1:=M^{(\sigma_0, \sigma_1)}_1\,,\  &\ldots& N^{m-1}:=M^{(\sigma_0, \sigma_1)}_{m-1}&\,,\  &N^m:=M^{(\sigma_0, \sigma_1)}_m=N_{\sigma_1}\\\\
&&N^{m+1}:=M^{(\sigma_1, \sigma_2)}_1&\ldots& N^{2m-1}:=M^{(\sigma_1, \sigma_2)}_{m-1} && N^{2m}:=M^{(\sigma_1, \sigma_2)}_m=N_{\sigma_2}\\\\
&&&\vdots&&\\\\
&&N^{(n-1)m+1}:=M^{(\sigma_{n-1}, \sigma_n)}_1&\ldots& N^{nm-1}:=M^{(\sigma_{n-1}, \sigma_n)}_{m-1} && N^{nm}:=M^{(\sigma_{n-1}, \sigma_n)}_m=N_{\sigma_n}\\\\
&&N^{nm+1}:=M^{(\sigma_{n}, \sigma_0)}_1&\ldots& N^{(n+1)m-1}:=M^{(\sigma_{n}, \sigma_0)}_{m-1} && N^{(n+1)m}:=M^{(\sigma_{n}, \sigma_0)}_m=N_{\sigma_0}\\\\
\end{array}
\end{eqnarray}
By~\eqref{ciclo}, we have
$$N^{0} \stackrel{f } \Longrightarrow N^1\stackrel{f }\Longrightarrow \cdots  \stackrel{f } \Longrightarrow N^{(n+1)m}\,.$$
Moreover,  $N^0$ and $N^{(n+1)m}$ are defined in~\eqref{ciclo1} so as to verify
\begin{equation}\label{ciclo2} N^0=N_{\sigma_0}=N^{(n+1)m}\end{equation}
(the last row of definitions in~\eqref{ciclo1} has precisely the r\^ole of making~\eqref{ciclo2} true). Applying Theorem~\ref{symbolic dynamics} with
$$k=(n+1)m\,,\quad f_i={f }\,,\quad N_i=N^i\quad   \forall\ i=1\,,\ldots\,, (n+1)m\,,\ N_0:=N^0$$
we infer
the existence of $x_0\in N^0=N_{\sigma_0}$ such that
\begin{itemize}
\item[{\rm (i)}] ${f }^i(x_0)\in N^{i}\quad \forall\ i=1\,,\ldots\,, {(n+1)m}\,;$
\item[{\rm (ii)}] ${f }^{(n+1)m}(x_0)= x_0\,.$
\end{itemize}
Taking, in (i), $i=m$, $2m$, $\cdots$, $nm$, we have the thesis. $\quad \square$

 \begin{remark}\rm As also remarked in~\cite{GierzkiewiczZ19}, if $A  \stackrel{f} \implica B \stackrel{f} \implica C$, not necessarily $A  \stackrel{f^2} \implica C$. Therefore, under conditions of Corollary~\ref{cor: symb dyn}, we cannot conclude that $f^k$ has an horseshoe.
\end{remark}

 \subsection*{Symbolic dynamics for ${\cal P}_{\!\textrm{\tiny$\cal H$,$\Pi^*$}}$} 
\begin{figure}[H]
\centering
\includegraphics[width=7.5cm,height=5cm,draft=false]{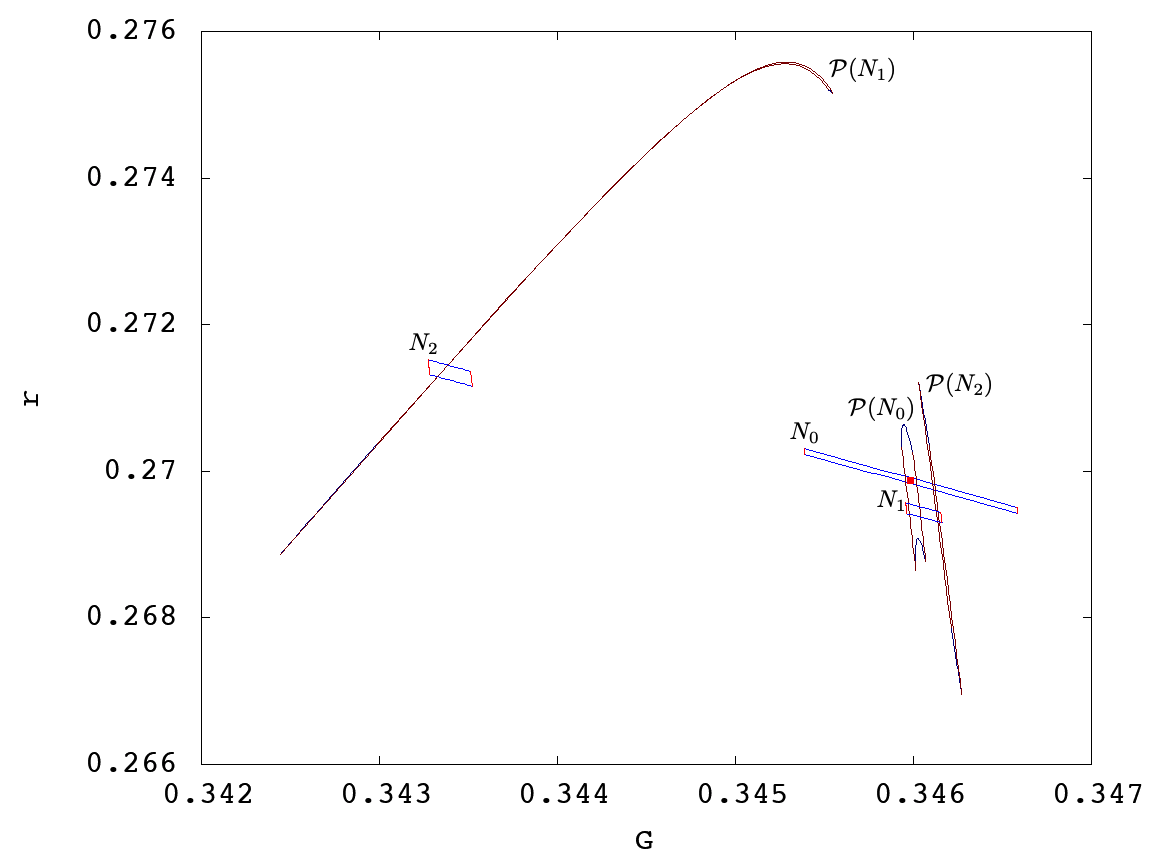}
\includegraphics[width=7.5cm,height=5cm,draft=false]{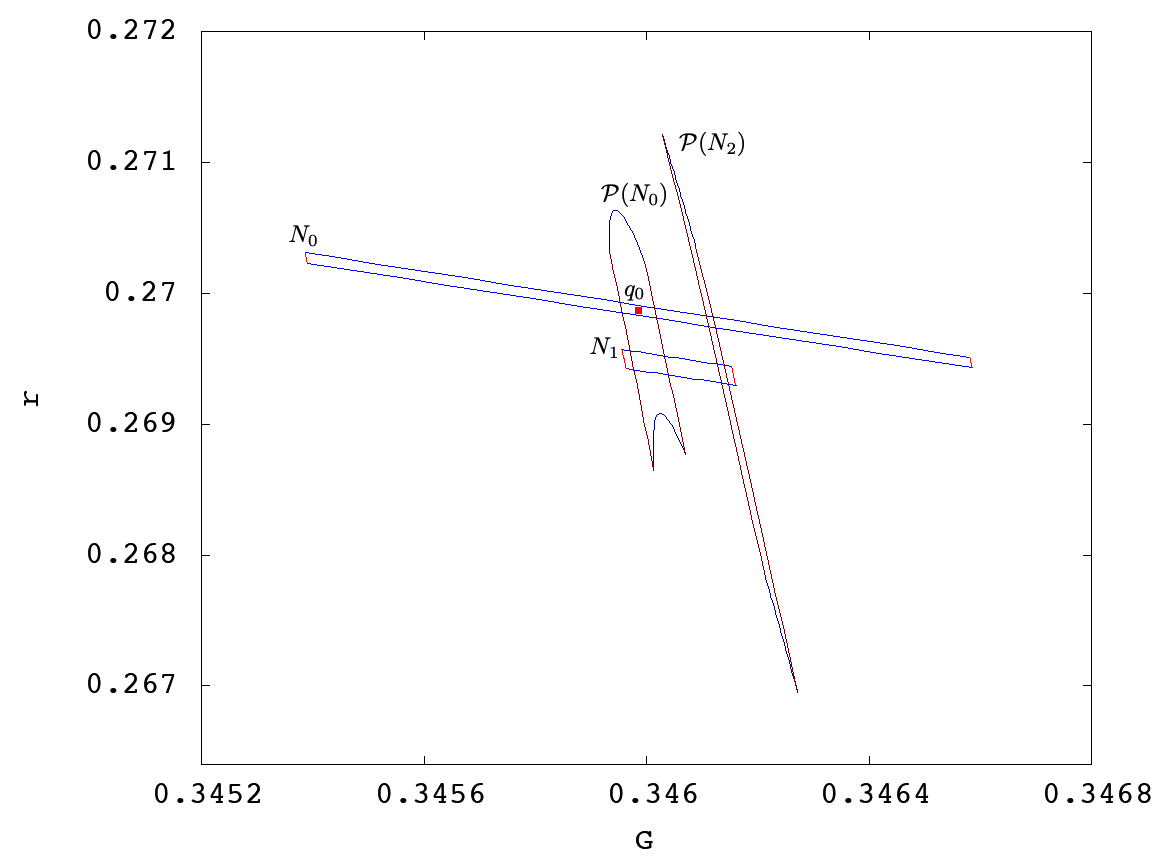}
	\caption{Numerical Evidence~\ref{num ev: cov rel}. Red represents the entry sets and their images and blue the exit sets and their images. The fixed point $\rm(r_0, G_0)$ in~\eqref{q0} is marked in red.  } \label{fig_horseshoe} 
\end{figure}

Let us consider the map ${\cal P}_{\!\textrm{\tiny$\cal H$,$\Pi^*$}}$ in~\eqref{poinc3}. The stable and unstable eigenvectors related to $D{\cal P}_{\!\textrm{\tiny$\cal H$,$\Pi^*$}}$ at 
 \begin{equation}\label{q0}q_0=({\rm r}_0, {\rm G}_0)=(0.26987, 0.345986)
 \end{equation}
have directions, respectively,
$$v^s= (-0.556268\,, \, 0.831003) \, , \quad v^u=(-0.998774 \, , \, 0.0495113),$$
and the angle between them is $\alpha = 0.296467 \, \pi$.
Observe that $q_0$ is the projection of the point~\eqref{saddlepoint} on the plane $\rm (r, G)$.
We denote as $N_0$ the parallelogram through $q_0$ with edges parallel to $v^s$ and $v^u$, namely: \begin{equation}\label{N0}
 	N_0 = q_0 +A_0 v^s +B_0 v^u \, ,
 \end{equation}
 where $A_0$, $B_0$ are the real intevals 
 $$A_0=[-0.000719075, 0.000719075]\,,\qquad B_0=[-0.0000400491, 0.0000400491] \, .$$
We define two analogous  parallelograms:
 \begin{equation}\label{N1}
 		N_1 =  q_1 + A_1 \widetilde v^s + B_1  v^u\,,\qquad
 		N_2 =  q_2 + A_2 v^s + B_2 v^u 
 \end{equation}
where
 \begin{equation*}
 	\left\{
 	\begin{array}{l}
 		q_1 =({\rm r}_1, {\rm G}_1)  =(0.269552,0.34598)\\
 		q_2 =({\rm r}_2, {\rm G}_2) =(0.27124,0.343432)
	\end{array}
 \right.	 \, ,
 \end{equation*}
 with
   \begin{equation*}
 		A_1=[-0.000028763,0.000208532] \,,\qquad  B_1 = [-0.000144177,0.00000400491] \, ,
 \end{equation*}
  \begin{equation*}
 		A_2 = [-0.000179769,0.000107861]\, , \qquad  B_2 = [-0.00000400491,0.000200246] 
 \end{equation*}
 and		 
 $$\widetilde v^s=(-0.556143,0.831003) \, .$$ 
Then we have the following (see Figure~\ref{fig_horseshoe})
\begin{numevid}\label{num ev: cov rel}
 \begin{equation*}
 	N_0  \stackrel{{\cal P}_{\!\textrm{\tiny$\cal H$,$\Pi^*$}}} \implica N_0 \stackrel{{\cal P}_{\!\textrm{\tiny$\cal H$,$\Pi^*$}}} \implica N_1 \stackrel{{\cal P}_{\!\textrm{\tiny$\cal H$,$\Pi^*$}}} \implica N_2 \stackrel{{\cal P}_{\!\textrm{\tiny$\cal H$,$\Pi^*$}}} \implica N_0 \quad , \quad  N_2  \stackrel{{\cal P}_{\!\textrm{\tiny$\cal H$,$\Pi^*$}}} \implica N_1 \, .
 \end{equation*}
 \end{numevid}
 Splitting such relations as
\begin{eqnarray*}
\left\{
\begin{array}{lll}
\displaystyle N_0 \stackrel{{\cal P}_{\!\textrm{\tiny$\cal H$,$\Pi^*$}}} \implica N_1 \stackrel{{\cal P}_{\!\textrm{\tiny$\cal H$,$\Pi^*$}}} \implica N_2 \stackrel{{\cal P}_{\!\textrm{\tiny$\cal H$,$\Pi^*$}}} \implica N_0\\\\
\displaystyle N_0 \stackrel{{\cal P}_{\!\textrm{\tiny$\cal H$,$\Pi^*$}}} \implica N_1 \stackrel{{\cal P}_{\!\textrm{\tiny$\cal H$,$\Pi^*$}}} \implica N_2 \stackrel{{\cal P}_{\!\textrm{\tiny$\cal H$,$\Pi^*$}}} \implica N_1\\\\
\displaystyle N_1 \stackrel{{\cal P}_{\!\textrm{\tiny$\cal H$,$\Pi^*$}}} \implica N_2 \stackrel{{\cal P}_{\!\textrm{\tiny$\cal H$,$\Pi^*$}}} \implica N_0 \stackrel{{\cal P}_{\!\textrm{\tiny$\cal H$,$\Pi^*$}}} \implica N_0\\\\
\displaystyle N_1 \stackrel{{\cal P}_{\!\textrm{\tiny$\cal H$,$\Pi^*$}}} \implica N_2 \stackrel{{\cal P}_{\!\textrm{\tiny$\cal H$,$\Pi^*$}}} \implica N_0 \stackrel{{\cal P}_{\!\textrm{\tiny$\cal H$,$\Pi^*$}}} \implica N_1
\end{array}
\right.
\end{eqnarray*}
 and in view of Corollary~\ref{cor: symb dyn}, the Numerical Evidence~\ref{numvid: symbdyn} follows, with $N_0$, $N_1$, $N_2$ as in~\eqref{N0},~\eqref{N1}.

\begin{figure}[H]
\centering
\includegraphics[width=7.5cm,height=5cm,draft=false]{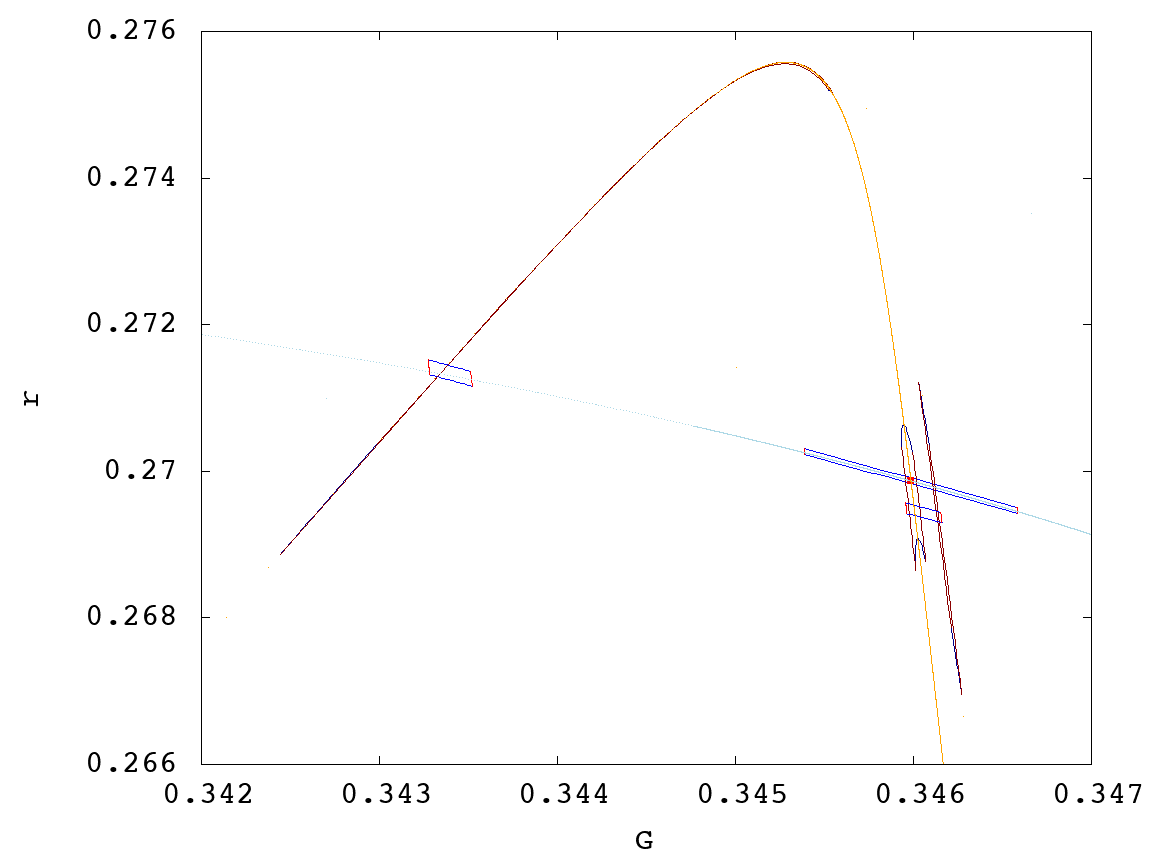}
\includegraphics[width=7.5cm,height=5cm,draft=false]{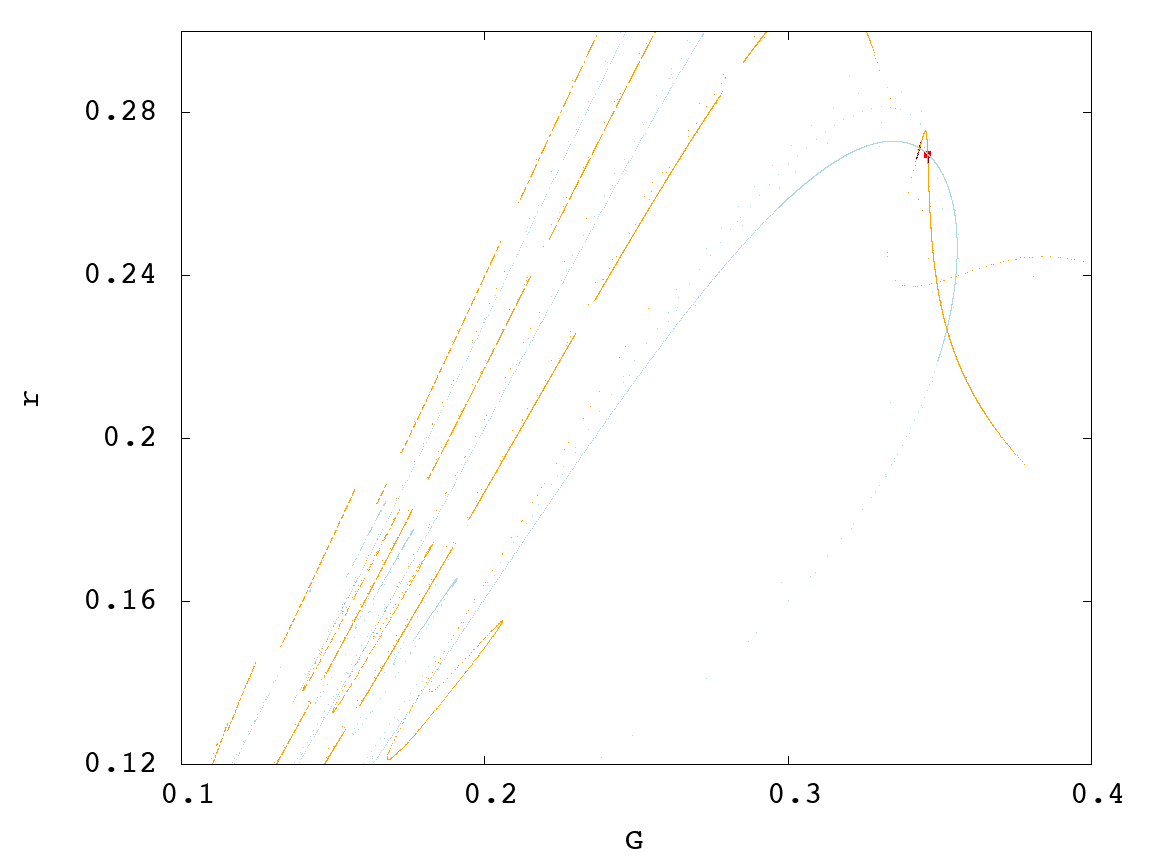}
	\caption{Details on the construction of  of the sets~\eqref{N0} and~\eqref{N1}.  Light blue and orange denote, respectively, the stable and unstable manifolds. The right figure represents a wider region.} \label{fig_horseshoe2} 
\end{figure}

\noindent
We conclude this section providing some detail on the construction of the sets~\eqref{N0} and~\eqref{N1}. As highlighted in Figure~\ref{fig_horseshoe2}, such sets are obtained inspecting the homoclinic intersections of the stable and unstable manifolds through $q_0$.

 \section{Control of errors and conclusions}\label{Control of errors and conclusions}
 In this section we describe how we controlled numerical errors and draw some conclusions.

\begin{figure}[H]
\centering
\includegraphics[width=7.5cm,height=3cm,draft=false]{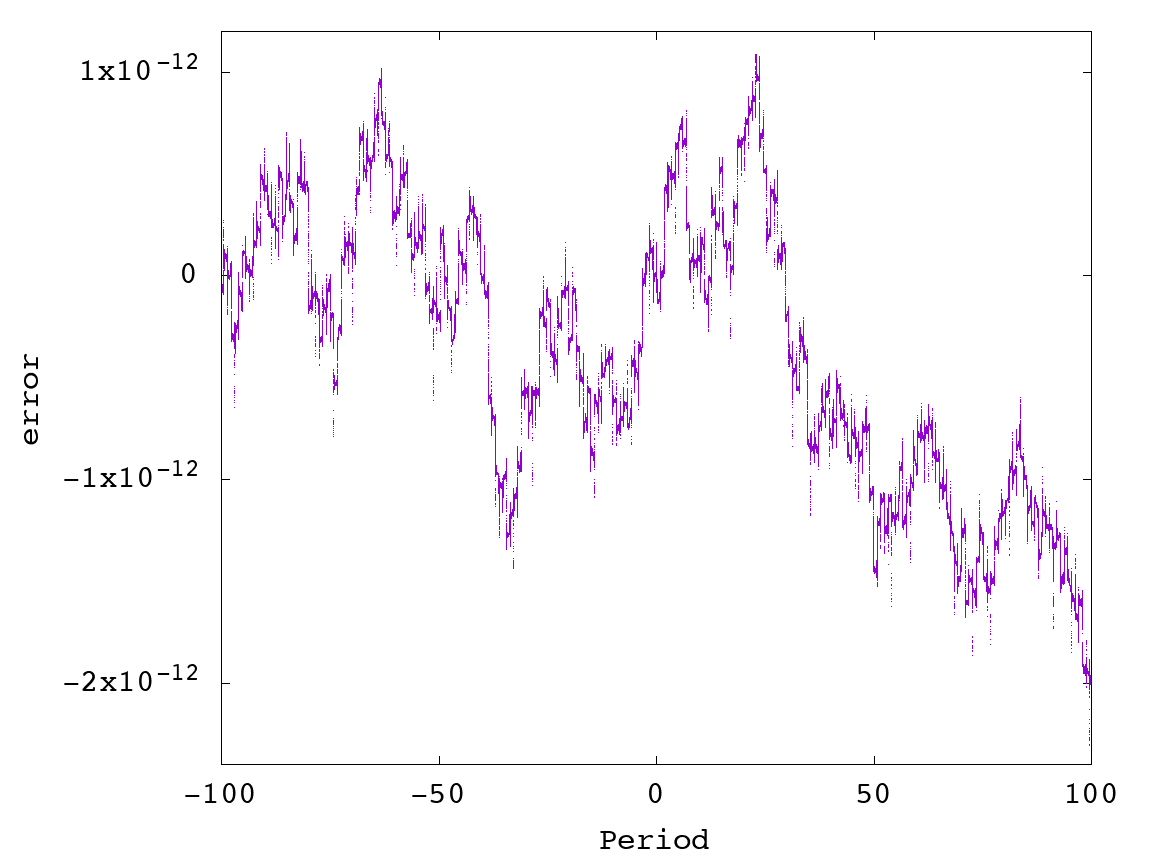}
\caption{Relative energy error in the propagation of the periodic orbit $\Gamma_{\rm s}$ versus number of its period $T_{\rm s}$. } \label{ham_error} 
\end{figure}

 \noindent
 In our computations, we used a double precision. One check of errors was performed by the control of energy which, being a first integral of motion should be constant. Its relative variation was required not to exceed $10^{-10}$, but the error we obtain in our simulations is much smaller. For the orbits we deal with in Section~\ref{eul_integral},~\ref{Neighbourhoods},~\ref{Symbolic dynamics} the relative error is comparable, so we choose to show orbit $\Gamma_{\rm s}$. In Figure~\ref{ham_error}, it can be seen that the relative error remains less that $2.5 \cdot 10^{-12}$ in 200 iterations of the map~\eqref{poinc}.  
 
\begin{figure}[H]
\centering
\includegraphics[width=10cm,height=7cm,draft=false]{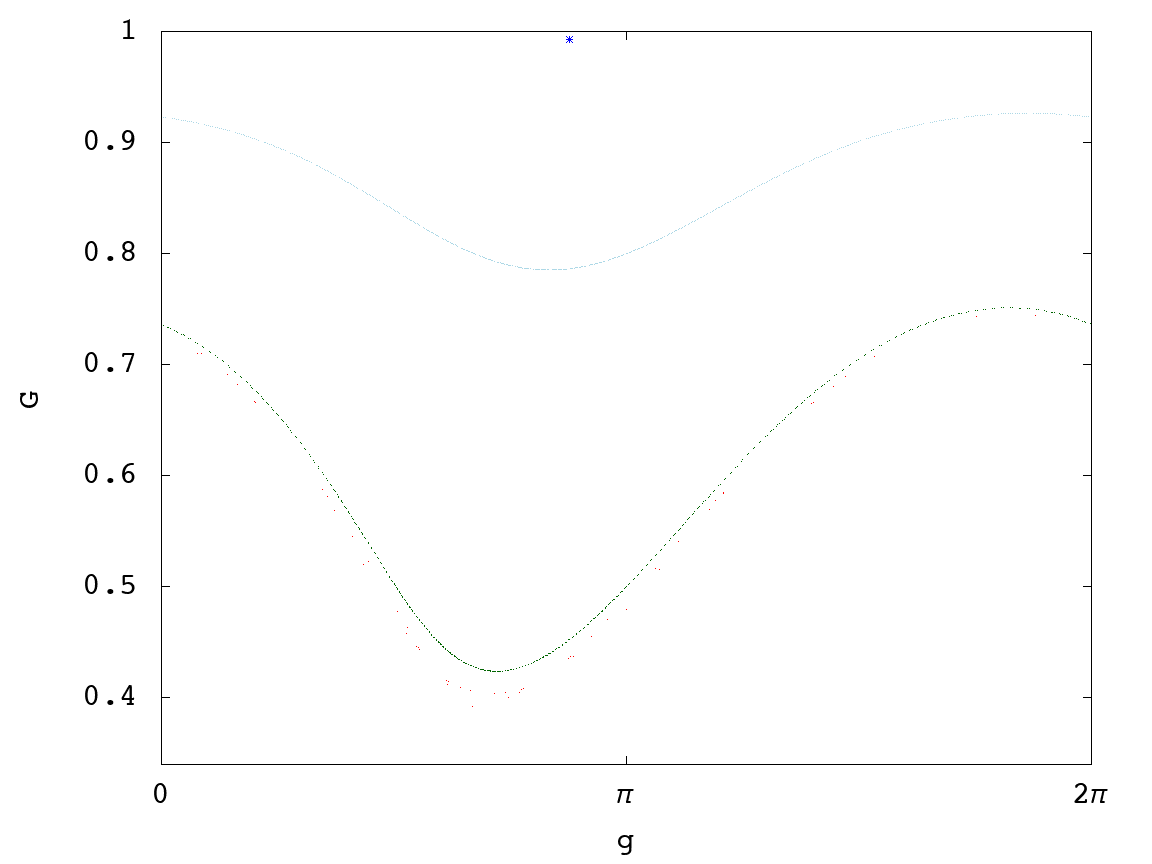}\\
	\caption{Sections of map~\eqref{poinc}  of the 4 orbits starting with initial conditions~\eqref{eq: errors}.} \label{section_errors} 
\end{figure}

\begin{figure}[H]
\centering
\includegraphics[width=7.5cm,height=3cm,draft=false]{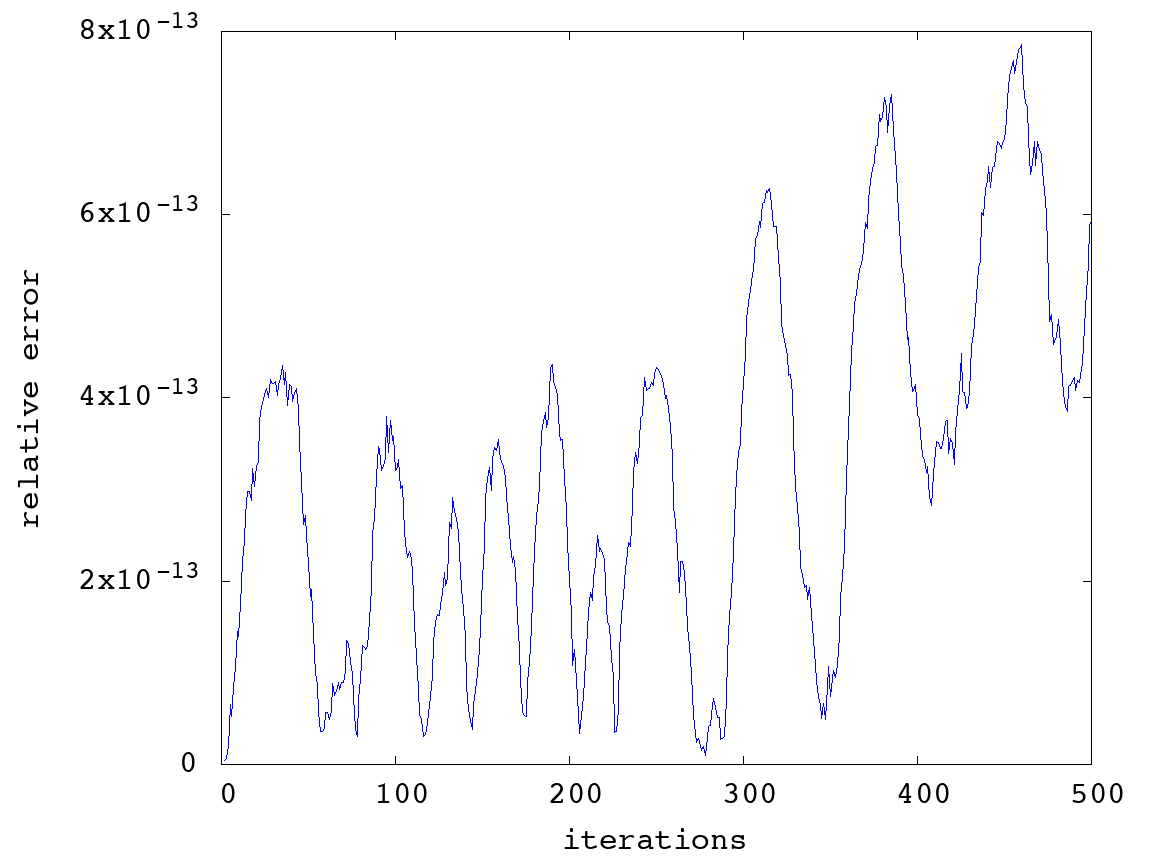}
\includegraphics[width=7.5cm,height=3cm,draft=false]{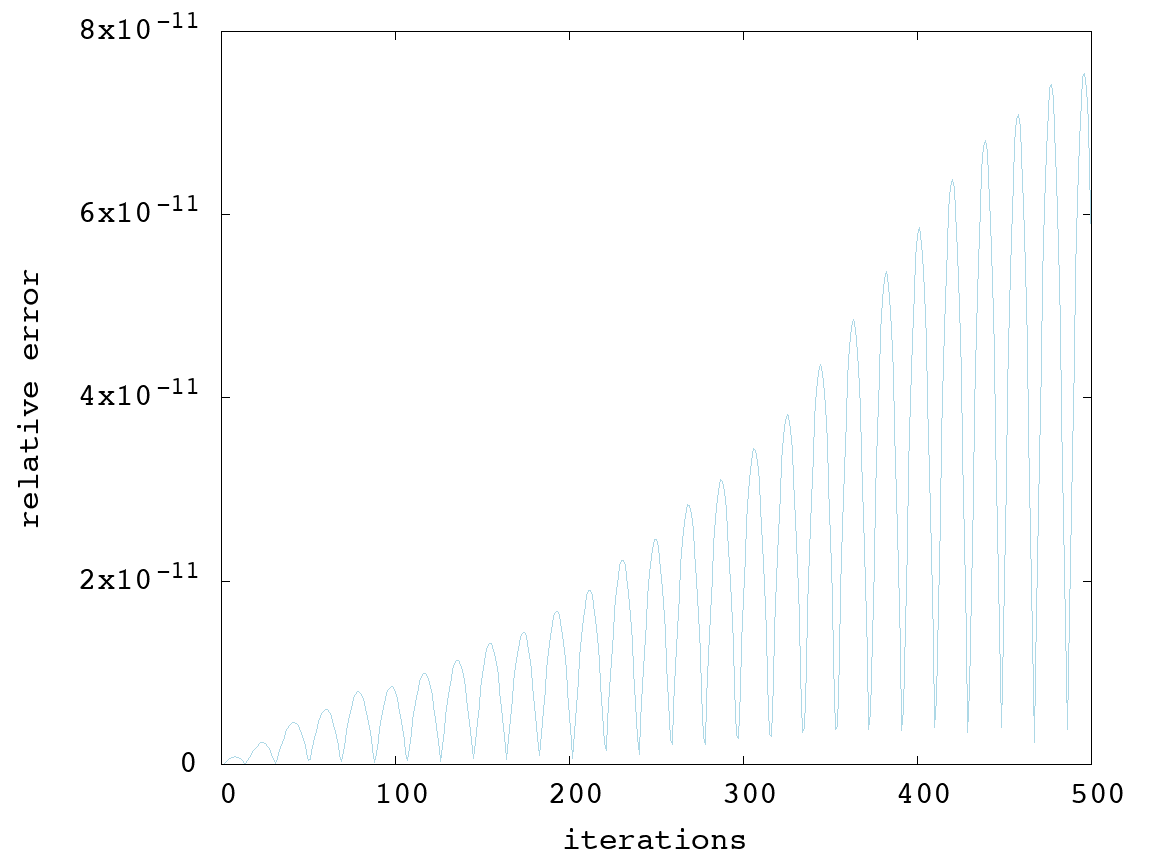}\\
\includegraphics[width=7.5cm,height=3cm,draft=false]{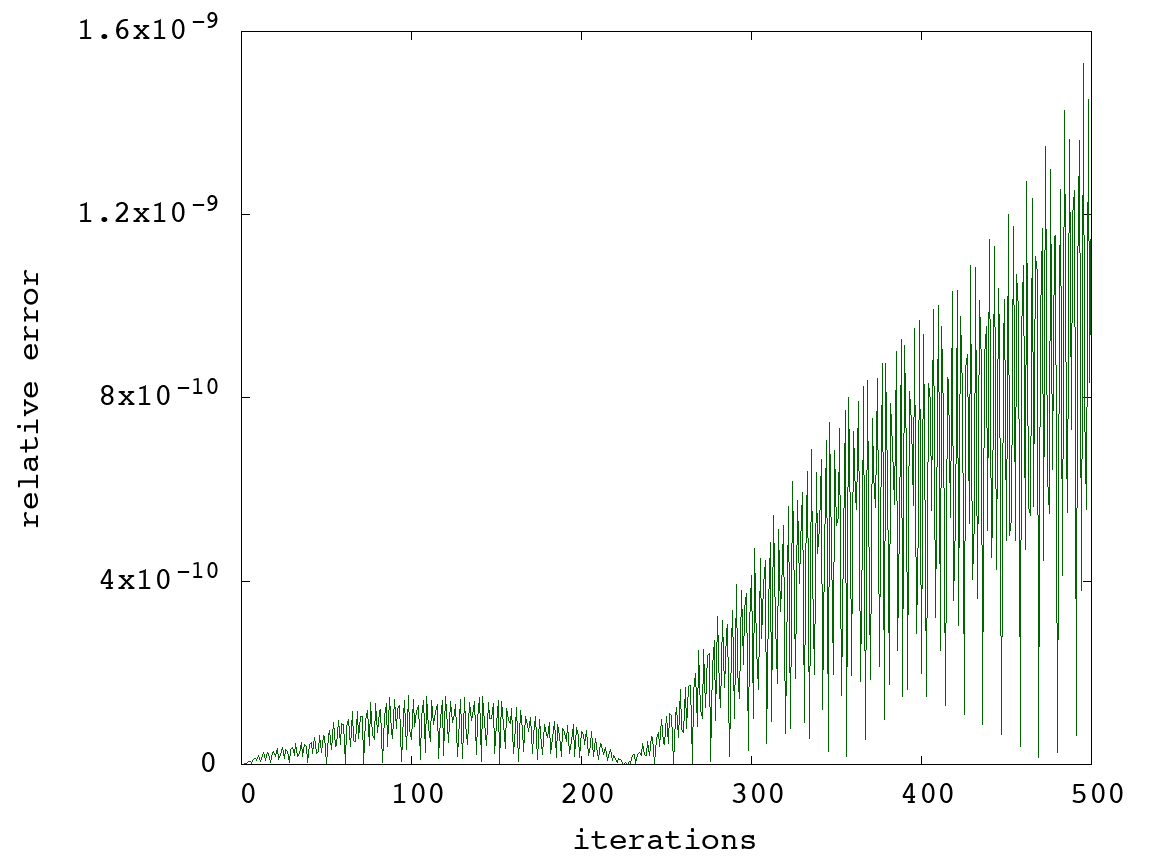}
\includegraphics[width=7.5cm,height=3cm,draft=false]{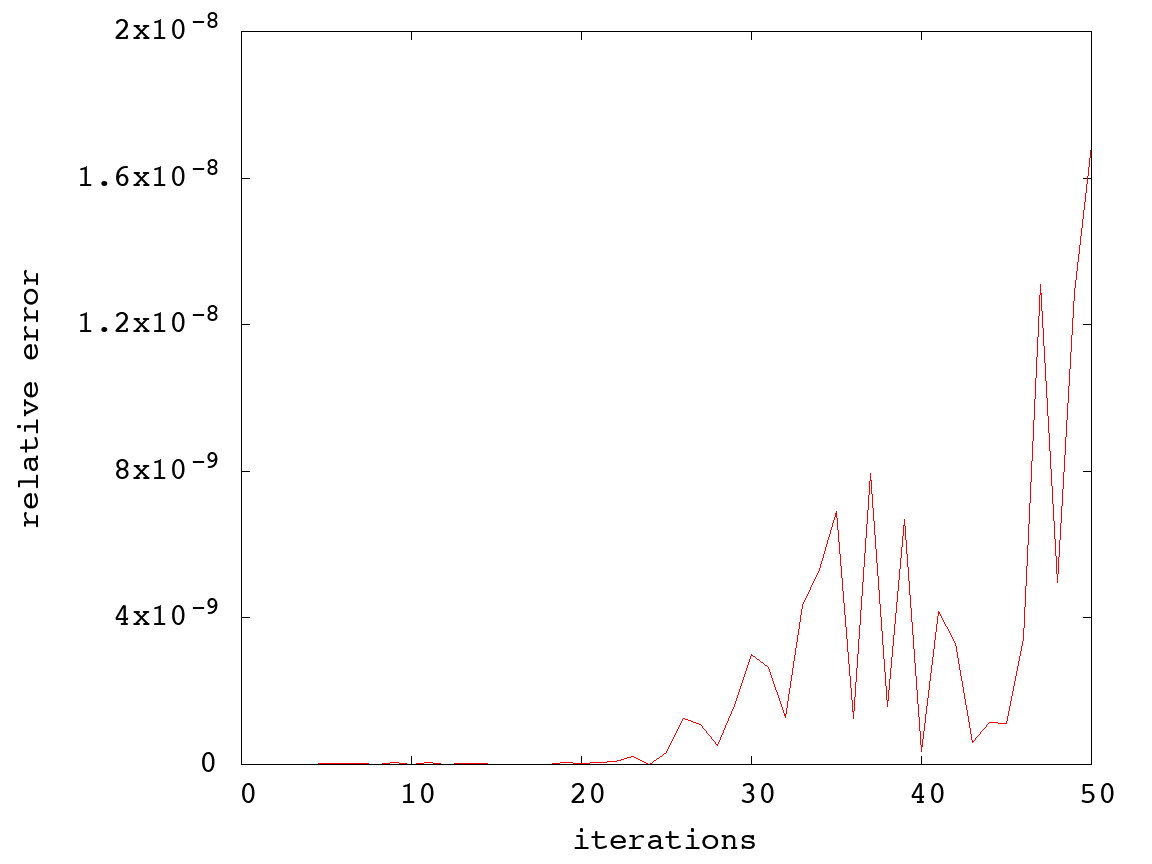}
	\caption{Onward and backwards integration of orbits of map~\eqref{poinc} starting with initial conditions~\eqref{eq: errors} (respectively, blue, light-blue, green and red): relative errors of the coordinates $\rm(G,g)$ of the onward and backward integrations versus the number of iteration of the map are plotted.}\label{errors} 
\end{figure}

 \noindent
 As a further test, we performed onward and backward integrations of orbits of the map~\eqref{poinc} starting with different initial conditions; as example cases, we show 4 orbits with the following initial conditions:
 \begin{equation}\label{eq: errors}
  	\left\{
 	\begin{array}{l}
\rm R =  -11.3668 \, , \quad G = 0.992515 \, ,\quad  r  = 0.13165\, ,\quad  g = 0.878179\,  \pi  \quad (blue)\\
\rm R  = -10.6704 \, , \quad G = 0.8 \, ,\quad  r = 0.13165 \, ,\quad  g = \pi \quad \qquad\qquad\qquad \,\, (light-blue) \\
\rm R = -9.07533 \, , \quad G = 0.5\, ,\quad  r = 0.13165\, ,\quad  g = \pi   \quad \qquad\qquad\qquad \,\, (green) \\
\rm R = -8.94348\, , \quad G = 0.48 \, ,\quad  r = 0.13165 \, ,\quad  g = \pi \quad \qquad\qquad\qquad (red) 
  \end{array}
  \right.  
 \end{equation}
 In Figure~\ref{section_errors}, we plot the sections map~\eqref{poinc}  of the 4 orbits with initial conditions~\eqref{eq: errors}. In Figure~\ref{errors}, we show the errors performed after a number of iterations onward and backward of the map~\eqref{poinc}; the 4 panels show the errors of the 4 orbits (respectively with the same colors) versus the number of iterations of the map~\eqref{poinc}. The number of iterations (500) is chosen as a reference because our simulations do not exceed this number.

\noindent
We are now ready for the conclusions. 

\noindent
In this paper, we discussed about the effects of the level sets of the function~\eqref{euler} on the dynamics of the Hamiltonian~\eqref{ham}. Specifically, in a range where the energy of the averaged, reduced 2--degrees of freedom system~\eqref{ham} has three different scales, one expects that the motions of the system obey to  Conjectures~\ref{picture of motion} and~\ref{conj}. In particular, due to the non--integrability of the system, chaos is expected closely to the envelope ${\cal M}_0$ in~\eqref{M0M1} of the separatrices of $\EE$. After  fixing the  energy level~\eqref{energy}, we computed a Poincar\'e map~\eqref{poinc}, which showed the existence of only two fixed points, having elliptic, hyperbolic character. The level sets of $\EE$ turn to vary a little along the orbit $\Gamma_{\rm s}$ generated by the elliptic fixed point, while it varies more along the orbit $\Gamma_{\rm u}$ generated by the hyperbolic fixed point. However, it turns out that $\Gamma_{\rm u}$ spends most of its time close to the saddle of ${\cal M}_0$, and we investigated the phase space around $\Gamma_{\rm u}$. We  used various 2--dimensional first return maps, and we found  a homoclinic tangency using one of them;  some heteroclinic intersection using another one. Applying the analysis developed in~\cite{WilczakZ2003, ZgliczynskiG2004, GierzkiewiczZ19}, we found $3$--symbolic dynamics in the sense of Definition~\ref{def: symb dyn}.
Our results are so in complete agreement with Conjecture~\ref{conj}, while, as remarked in the introduction, Conjecture~\ref{conjOLD} is still open.

\vskip.2in
\noindent
{\bf Acknowledgments}
We thank the anonymous reviewers for their stimulating remarks, which helped to improve the presentation of the results in the paper.
This paper is supported by the the ERC project 677793 Stable and Chaotic Motions in the Planetary Problem (2016--2022).\\
Figures~\ref{figure1}--\ref{figure3} have been produced with Mathematica${}^\textrm{\textregistered}$. Figure~\ref{model} has been produced with Vectornator and Figures~\ref{energia2}--\ref{errors} have been produced with Gnuplot.

\newpage
\addcontentsline{toc}{section}{References}
\def\cprime{$'$} \def\cprime{$'$}

\end{document}